\documentclass[12pt,a4paper]{article}
\usepackage{amssymb}
\usepackage{amsthm}
\usepackage{amsmath}
\usepackage{pstcol, pst-node}

\newtheorem{thm}{Theorem}[section]
\newtheorem{assum}[thm]{Assumption}
\newtheorem{cor}[thm]{Corollary}
\newtheorem{prop}[thm]{Proposition}
\newtheorem{lem}[thm]{Lemma}
\newtheorem{rem}[thm]{Remark}
\newtheorem{claim}[thm]{Claim}

\theoremstyle{definition}
\newtheorem{defn}[thm]{Definition}
\newtheorem{q}[thm]{Question}
\newtheorem{convention}[thm]{Convention}
\newtheorem{prop-def}[thm]{Proposition-Definition}
\newtheorem{note}[thm]{Note}

\theoremstyle{remark}
\newcommand{\ra}{\rightarrow}

\newcommand{\g}{\mathfrak{g}}

\newcommand{\T}{\mathfrak{t}}
\newcommand{\be}{\begin{equation}}
\newcommand{\bc}{\begin{cor}}
\newcommand{\bt}{\begin{thm}}
\newcommand{\bl}{\begin{lem}}
\newcommand{\bpr}{\begin{prop}}
\newcommand{\br}{\begin{rem}}
\newcommand{\bd}{\begin{defn}}
\newcommand{\ee}{\end{equation}}
\newcommand{\et}{\end{thm}}
\newcommand{\el}{\end{lem}}
\newcommand{\epr}{\end{prop}}
\newcommand{\er}{\end{rem}}
\newcommand{\ed}{\end{defn}}
\newcommand{\ec}{\end{cor}}
\newcommand{\R}{\Bbb{R}}
\newcommand{\C}{\Bbb{C}}

\newcommand{\pll}{\parallel}

\newcommand{\hs}{\hspace*}

\newcommand{\ep}{\epsilon}
\newcommand{\nnn}{\noindent}
\newcommand{\lam}{\lambda}
\newcommand{\Lam}{\Lambda}
\newcommand{\del}{\partial}

\newcommand{\Z}{\Bbb{Z}}
\begin{document}
\date{} 
\title{Convergence of Hermitian-Yang-Mills Connections on K\"ahler Surfaces
 and Mirror Symmetry}

\author{Takeo Nishinou}

\maketitle

\section{Introduction}

The purpose of this paper is to study the limiting behavior of
  a family of Hermitian-Yang-Mills (HYM) connections on K\"ahler
 $T^4$ with an appropriate affine structure when the metric
 on $T^4$ goes to the adiabatic limit.
Through this analysis we exhibit a natural construction
 of (special) Lagrangian manifolds based on the idea of mirror symmetry.
The main point in our analysis is that we allow reducibility 
 of the connections
 (more precisely, reducibility of them along the fibers).
On the other hand, our methods are also applicable to 
 irreducible cases treated in
 \cite{C}, \cite{DS}, \cite{F1}.
In these cases, the discussion
 applies to bundles over products of Riemannian
 surfaces.
We give
 sketches of proofs of parts of results of \cite{DS}, 
 and give a stronger version of the theorem of \cite{C}.

We consider a K\"ahler torus $T^4$ where there is a 
 Lagrangian fibration structure:
$$
T^4 \to T^2_B,
$$
 and some affine structure, called Hessian geometry (see section 2).
Then we take the limit where the diameter of fibers goes to zero.
Such a limit was considered by Dostoglou-Salamon \cite{DS}
 in relation to the Atiyah-Floer conjecture and also appears
 in mirror symmetry as a description of large structure limit
 \cite{F3}, \cite{GW}, \cite{SYZ}.

In fact this paper was motivated by both of them.
On one hand, we are going to try to extend the analysis of \cite{DS}
 to the cases where reducible connections on fibers appear.
On the other hand, we want to give a basic piece for the gauge theory
 in the large structure limit.
Since the latter is related to holomorphic bundles through
 Kobayashi-Hitchin correspondence,
 it is also related, through mirror symmetry, 
 to the symplectic geometry of the mirror manifold.
Namely, through our analysis, we can construct a
 (special) Lagrangian submanifold (called type A D-brane)
 from a family of stable holomorphic bundles (type B D-branes).

We first describe the terms which appeared here with some results
 associated to them.\\

\noindent
$\bold{Mirror \; symmetry}.$
Mirror symmetry originated from
 a (rather apparent) symmetry in superconformal field theory (SCFT). 
Since there are SCFT's constructed by geometrical means
 ($\sigma$-model, roughly speaking a theory of maps to Riemannian
 manifolds (especially Calabi-Yau manifolds)
 from Riemannian surfaces), it has implications to geometry.
In a naive form, it anticipates that given a Calabi-Yau manifold,
 there should be another Calabi-Yau manifold (the mirror pair).
In $\sigma$-model, it is required that theories constructed from
 each piece of a mirror pair coincide. 
More precisely, A-model on a Calabi-Yau $X$ should be isomorphic
 to B-model on the mirror $\check X$ (see below for A- and B-models).
The mathematical implications from this statement are, for example,
 mirror relation of Hodge diamonds and correspondence of 
 Yukawa-cuplings.
The latter means the Gromov-Witten potential (roughly, generating function
 for the number of rational curves) in A-model and
 the solution of a Picard-Lefschetz differential equation in B-model.
The seminal paper \cite{COGP} predicted the number of rational curves
 on the quintic in $\mathbb{P}^4$ using this correspondence. \\

\noindent
$\bold{D-branes}.$
Although the physics flavored name, they first appeared
 in mathematical work by Kontsevich \cite{K}.
There he claimed mirror symmetry should be formulated
 as an equivalence of triangulated categories 
 in symplectic and complex geometry.
D-branes are objects of these categories and
 the D-branes in the symplectic side (called A-model in physics literature)
 are called type A (or A-brane)
 and those in complex side (B-model) are called type B (or B-brane).
Mathematically, a type A D-brane is, roughly, a Lagrangian 
 submanifold with a line bundle and a connection on it,
 and a type B D-brane is an object of the derived category
 of coherent sheaves.
As we have noted, we will construct an A-brane from a family of 
 B-branes on the mirror, that is, 
 we construct a Lagrangian submanifold (with singularities) on a
 K\"ahler torus
 from a family of stable bundles on the mirror.
We note that we have not yet succeeded
 to attach a line bundle with a connection 
 to this Lagrangian submanifold
 (but we can attach some flat bundle.
 See section 8).\\

\noindent
$\bold{Strominger-Yau-Zaslow \; picture} \cite{SYZ}.$
Strominger-Yau-Zaslow (SYZ) picture is one of formulations of
 mirror symmetry, based on so-called T-duality in string theory.
According to this picture, Calabi-Yau manifolds should be equipped with special
 Lagrangian torus fibrations (with singular fibers) if the complex structures
 are near `large structure limit'.
Moreover, the mirror of them should be given by suitable
 compactifications of the dual torus fibrations.

The `large structure limit' (sometimes called the maximal degeneration limit
 in view of Hodge theory) has (conjectural) metric description.
Namely, the Calabi-Yau metric (i.e, Ricci flat K\"ahler metric)
 should degenerate in such a way
 that the diameter of the fiber of special Lagrangian
 fibration shrinks, as the complex structures come close to
 the limit.
This description was recently established in the K3 case, by
 Gross and Wilson \cite{GW}.\\

\noindent
$\bold{Adiabatic \; limit}.$
In \cite{DS}, Dostoglou and Salamon solved the Atiyah-Floer
 conjecture (here meaning the isomorphism of instanton
 Floer homology of a three manifold 
 and symplectomorphism Floer homology of the moduli space of
 flat bundles on a Riemannian surface)
 for $P_f \times \R$, where $P_f$ is a mapping torus
 of a Riemannian surface $\Sigma$, with the diffeomorphism $f$,
 in the case of those $SO(3)$ bundles which when restricted to 
 $\Sigma$, are nontrivial.
They analyzed the behavior of anti self dual (ASD) connections
 on the bundle when the diameter of $\Sigma$ goes to zero.
This limit is called the adiabatic limit.
They showed the isomorphism of the moduli spaces,
 one is the space of ASD connections over $P_f \times \R$ and the
 other is the space of pseudo-holomorphic strips on the 
 moduli space of flat bundles connecting the fixed
 points of the symplectomorphism induced by $f$,
 when the diameter of the fibers is sufficiently small.\\

\noindent
$\bold{Relation \; between \; SYZ \; and \; adiabatic \; limit}.$
Given two topics both concerned with the metric degeneration 
 shrinking the fibers, it is natural to seek for relations
 between them.
In fact, this is the starting point of this paper.
By the result of Gross and Wilson \cite{GW}, it is natural to expect that
 the adiabatic limit is the metric description of large structure
 limit.
At this limit, it is believed that the string theory degenerates
 to the classical geometry (without quantum corrections).
In fact, our analysis gives,
 through going to this limit, the correspondence of mirror objects,
 stable bundles on B-model and Lagrangian submanifolds on A-model.\\

\noindent
$\bold{Doubly \; periodic \; instantons} \cite{BJ, J, J2, J3}.$
Aside from mirror symmetry, there is a topic which we will be concerned with.
A doubly periodic instanton is the popular name of
 an ASD connection on $\R^2 \times T^2$.
There has been large progresses about this topic recently.
This is concerned to us through the analysis of bubbles,
 and our analysis is of different nature from those references.
Namely our aim is the energy quantization theorem for doubly
 periodic instantons, without any assumption on the behavior of
 the connection at infinity.
This point was relatively easy in those references, under the assumption 
 of curvature decay at infinity.\\

Now we describe the main ideas of our construction.
Suppose we are given a family of pairs  $(\widehat M_{\ep}, E_{\ep})$,
 $\ep \in (0, 1]$, of
 2-dimensional K\"ahler tori and 
 stable holomorphic vector bundles (i.e. B-branes) of fixed
 topological type on them
 and each torus has a
 structure of a Lagrangian torus fibration
 $\pi:\widehat M_{\ep} \ra T_B^2$ whose fibers are of diameter $O(\ep)$.
Let $(M, \omega)$ be the SYZ mirror symplectic manifold of $\widehat M_1$
 (see definition 2.1).
Let $\Xi_{\ep}$ be a family of hermitian
 Yang-Mills(HYM) connections on $E_{\ep}$,
 which exist by the theorem of Donaldson, Uhlenbeck and Yau \cite{D, UY}.

Our main result says,
 as $\ep$ goes to zero, $\Xi_{\ep}$ will, modulo possible bubbles,
 converge to a connection on 
 $E_1 \to \widehat{M_1}$ which is flat on each fiber.
Then, since each fiber is a torus, the limit connection will determine
 elements of the dual torus, which are by definition
 points of the fiber of the mirror
 $M_1$.
Then, the geometric version of our theorem says
 these points gather to make up a
 (special) Lagrangian variety (i.e. A-brane).\\
%


The contents of this paper are as follows.
In the next section, we review some geometrical facts about mirror symmetry,
 especially in the context of SYZ conjecture in semi-flat case
 and describe our
 main geometric result (theorem 2.6).
In section 3,
 we describe HYM equations near the adiabatic limit.
Then we state the analytic version of the main result (theorem 3.1)
 and prove it assuming the results in later sections.
Along the way, we review Dostoglou-Salamon's classification of bubbles 
 in our situation.
Section 4 is a preparation for sections 5 and 6.
We study the space of connections on the trivial $SO(3)$ 
 (or $SU(2)$) bundle
 over $T^2$, especially the action of complex gauge transformations.
The main result here is proposition 4.19,
 roughly saying if a connection $A$ is 
 close to a flat connection $A_0$, then
 there is a gauge such that
 $\pll A - A_0 \pll_{C^{r+1}}$
 is bounded by $\pll F_A \pll_{C^r}^{\frac{1}{2}}$.
In section 5, we treat doubly periodic instantons.
The main result here is the energy quantization for
 them without any assumption about curvature decay
 (theorems 5.5 and 5.6).
However, the most degenerate cases like trivial connection
 still disturb us from the unified formulation.
Namely, the energy bound we can show for the connections reduces to zero
 as they come close to those cases
 (but in the genuinely degenerate case, we again have a finite lower
 bound of the energy, theorem 5.6).
In section 6, we treat the bubbling of the type which gives a holomorphic map
 to the moduli of flat connections,
 and finally prove the $C^{\infty}$ convergence of ASD connections
 modulo gauge transformation (theorem 6.6).
In fact, one can specify the limit object associated to a family of
 stable bundles $E_{\ep}$ without proving the $C^{\infty}$
 convergence of HYM connections on these bundles,
 from the point of view of minimal surface theory
 (proposition 6.3).
This corresponds to the holomorphic maps from the base to the space
 of flat connections in Dostoglou-Salamon's case.
In fact, in their case too, one can define a family of holomorphic
 maps by taking the complex gauge equivalence classes
 once the existence of the ASD connection is known
 and so can also define the limit holomorphic map
 without showing the convergence of the connections
 (see the discussion at the end of section 6).
However, the relation between these holomorphic maps
 and HYM connections are subtle 
 and we prove the actual convergence of the HYM
 connections to a holomorphic curve.
This is important when we want to compare
 the moduli spaces of these objects (see remark 7.3).
It is a basic part of the analysis
 for the attempt to extend the Dostoglou-Salamon's result to 
 reducible cases.
Of course, the convergence is of independent interest from the point of view of gauge
 theory, in particular, gauge theory over large structure limit.

The methods in sections 5 and 6 are quite transparent.
That is, first construct an appropriate gauge by applying  
 more or less known general gauge fixing results, and then fine tune.
This latter process is needed to handle with the reducible connections.
Throughout the discussion, we systematically exploit the
 complex gauge transformations in these gauge fixing
 processes.
The point is that the construction is of local nature
 (see the paragraph after corollary 5.9).
This would also enable us to apply the same methods to other
 situations (see section 8).
In particular, our method can be applicable also to non-reducible cases
 like Dostoglou-Salamon \cite{DS} and Fukaya \cite{F1},
 and we can reproduce parts of their results
 without some parts of fine tuning process (so the argument becomes
 simpler than those in the main text).
As an illustration of how our method will be applied
 to irreducible cases, we reproduce some of
 results of \cite{DS}, especially the parts
 concerning bubbling analysis.

On the other hand, Chen \cite{C} treated a convergence problem of ASD
 connections on a product of Riemannian surfaces.
However, although not mentioned explicitly, it seems that in that 
 paper it is assumed that the moduli problem of flat connections on the
 Riemannian surface corresponding to the fiber is transversal and so the moduli
 space is
 smooth of expected dimension.
Moreover, it seems that there some analysis of bubbles
 (bubbles of type two in our language, see the proof of theorem
 3.1) is missing.
See the last of section 6 for these points.
At the last of section 6, we give our version of a formulation of 
 his theorem (in fact stronger than his, in that
 although in \cite{C} $C^0$-convergence was discussed,
 we give $C^{\infty}$ convergence) and prove it.

In section 7, we complete the proof of theorem 2.5.\\

%

I would heartfully like to thank my adviser Kenji Fukaya 
 for leading to this problem and 
 giving me useful suggestions and advice.
I thank Manabu Akaho for informing me about the
 erratum of $\cite{DS}$.
I would also like to thank referees of the earlier version of this paper
 for useful comments and pointing out some errors.
Part of this paper was written during my stay at 
 Universit\"at Freiburg. 
I am very grateful to Bernd Siebert and Mathematisches Institut
 for warm hospitality.\\

\noindent
$\bold{Conventions}.$
We sometimes write $\frac{\del}{\del x}$ as $\del_x$ 
 here $x$ is a coordinate of some manifold.
In various estimates of $C^r$ norms of functions and sections,
 constants will appear which will depend on $r$.
But we do not write this dependence explicitly,
 unless it is crucial, to avoid too many subscripts.

\section{Review of semi-flat mirror symmetry}

In this section, we review some of mirror symmetry and set up
 the particular situation we will be concerned with.
For more details about the materials in
 this section, refer to $\cite{F3}$, section 1 and references therein.

Let $(\widehat M, \omega, J)$ be a K\"ahler surface and $\pi:\widehat M \ra B$ be
 a Lagrangian torus fibration, $B$ is a surface (compact or not).
%
%
%
We assume $\pi$ has an $SL(2, \Bbb{Z})$ structure, namely,
 there is a $\Bbb{Z}^2$ bundle $\Lambda \ra B$ whose structure group
 is $SL(2, \Bbb{Z})$ and put $E = \Lambda \otimes_{\Bbb{Z}} \Bbb{R}$.
Then $\widehat M$ is isomorphic to the torus bundle $E/\Lambda$.
So $\widehat M$ has a natural flat structure induced from $\Lambda$
 (we call it Gauss-Manin
 connection) and we denote it by $\nabla^{GM}$.
Furthermore, $\widehat M$ has a canonical flat section, the zero section $s_0$.

Let $U \subset B$ be an open subset on which the bundle $\Lambda$
 trivializes.
Then on $\pi^{-1}(U)$, we have an action of $\R^2 / \Z^2 = T^2$ by addition,
 once we fix flat sections $\lam_1$, $\lam_2$ of $\Lambda|_U$ which
 generate $\Lambda|_U$ on each fiber. 
We denote this action by 
$$\phi_U : T^2 \times \pi^{-1}(U) \to \pi^{-1}(U),$$
 explicitly, it is given by
$$\phi_U(a, b; x) = x + a\lam_1 + b\lam_2,$$
 here $(a, b) \in \R^2 / \Z^2$ and $x \in \pi^{-1}(U)$.
Note that although this action is not globally defined in general,
 when we take a two trivializing neighbourhood $U, V \subset B$
 with $U \cap V \neq \varnothing$, $\phi_U$ and $\phi_V$
 are related by
$$
\phi_U(t, x) = \phi_V(\psi_{UV} t, x),
$$
 where $x \in \pi^{-1}(U \cap V)$, $t \in \R^2 / \Z^2$ and
 $\psi_{UV}$ is the transformation matrix between $\pi^{-1}(U)$ and
 $\pi^{-1}(V)$ for some fixed bases of fibers
 on $\pi^{-1}(U)$ and $\pi^{-1}(V)$.
In particular, $T^2$-action is transformed to $T^2$-action.
Our assumption, which is often assumed in mirror symmetry is that
 this flat structure on $\widehat M$ is compatible with the K\"ahler structure.
Namely, we assume the following conditions.\\

\nnn
$\bold{1}.$ The symplectic structure $\omega$ and the complex structure $J$ are
 invariant under the local $T^2$ action to fibers.
That is, for an open subset $U$ of $B$ as above, we have  
$$
\phi_U(t)^* \omega = \omega
$$
and
$$
\phi_U(t)_* \circ J (v) = J \circ \phi_U(t)_* (v),
$$
where $t \in T^2$ and $v$ is a tangent vector at some point
 $x \in \pi^{-1}(U)$.\\

\noindent
$\bold{2}.$ Every fiber and the zero section are Lagrangian submanifolds.\\

\noindent
$\bold{3}.$ $J(TB) = TF$, where $B$ is identified with 
 the image of the zero section and 
 $TF$ is the bundle on it whose fibers are tangent spaces to the fibers
 of $\pi$.\\

We call these conditions as K\"ahler $T^2$-structure.
\begin{defn}
We put $M = E^{\ast} / \Lam^{\ast}$,
 where $E^*$ and $\Lambda^*$ are duals of $E$ and $\Lambda$,
 respectively, and call it the mirror of $\widehat M$.
\end{defn}
Note that the symplectic form $\omega$ identifies $T^* F$ with 
 $TB$.
Similarly, the complex structure $J$ identifies $TF$ with $TB$.
With these identifications, flat affine connections
 are induced on $B$ by pulling back $\nabla^{GM}$ or its dual.
The closedness of $\omega$ and the integrability of $J$ imply
 torsion freeness of these connections (see \cite{F3}, section 1).
From this observation, we can prove the following facts
 for manifolds with K\"ahler $T^2$ structure.\\

\noindent
$\bold{I}.$ 
There is an isomorphism $\varphi$ between 
 the tangent bundle $TB$ of the base and $E^*$,
 such that the induced connection $\varphi^* (\nabla^{GM})^{\vee}$
 is a torsion free flat connection on $TB$,
 here $(\nabla^{GM})^{\vee}$ is the dual flat connection of
 $\nabla^{GM}$.
Thus, the base $B$ has an integral affine structure and locally we have
 affine coordinates $\{ s, t \}$
 such that for some linear (the origin is given by the
 zero section $s_0$) coordinates $\{ x, y \}$ of the torus,
 defined on fibers over trivializing neighbourhood $U \subset B$
 (necessarily multi-valued, but it does not matter to our purpose),
 $\{ s, t, x, y \}$ are Darboux coordinates for
 the symplectic structure of $\widehat M$ (symplectic form
 $\omega = ds \wedge dx + dt \wedge dy$).
 We call $\{ s, t \}$  symplectic coordinates for the base and
  $\{ s, t, x, y \}$  symplectic coordinates for $\widehat M$. \\

\noindent
$\bold{II}.$ The base $B$ has 
 affine coordinates $\{ \Check{s}, \Check{t} \}$
 such that for the above linear coordinates
 $\{ x, y \}$ of the torus,
 $\{ \Check{s} + i x, \Check{t} + i y \}$ are complex coordinates for
 the complex structure of $\widehat M$, namely
 $J(\frac{\del}{\del x}) = -\frac{\del}{\del \Check{s}}$ and
 $J(\frac{\del}{\del y}) = -\frac{\del}{\del \Check{t}}$ hold.
We call $\{ \Check{s}, \Check{t} \}$ complex coordinates for $B$ and
  $\{ \Check{s}, \Check{t}, x, y \}$  complex coordinates for $\widehat M$.\\

\noindent
$\bold{III}.$ The mirror manifold $M$ has a K\"ahler structure
 which also satisfies the conditions $1, 2, 3$.\\

\noindent
$\bold{IV}.$ Moreover, complex coordinates for the base $B$ of $\widehat M$
 become symplectic coordinates for the base $B$ of $M$
 and vice versa.
This symplectic coordinates and
 the dual coordinates of $\{ x, y \}$ of the fiber torus will be 
 local symplectic coordinates for the total space of $M$.\\
 
%

%
%
\nnn
$\bold{V}.$  There is an open cover $\{U_a\}_{a \in I}$ of the base $B$
 ($I$ is an index set) and smooth functions $h_a$ on them
 such that their second derivatives $g_{pq} = \frac{\del^2 h_a}{\del p 
 \del q}$, where $p, q$ is $s$ or $t$,
 form a Riemannian metric $g_a = \sum_{p, q} g_{pq} dpdq$ on $U_a$.
These metrics are compatible on the intersections of the
 coverings. \\

The metric on the whole space $\widehat M$
 is defined from this base metric using the 
 zero section, the invariance under the local $T^2$
 action and the invariance of the metric under the action of $J$.
This coincides with the metric defined by the K\"ahler structure
 $\omega, J$.
The following is known. 
\begin{prop}
This metric is Calabi-Yau if and
 only if the determinant $det(g_{pq})$ is 1. \qed
\end{prop} 
In this case, the mirror is also Calabi-Yau.
\begin{assum}
In the following, we will concentrate on the case when $\widehat M$
 has the Calabi-Yau metric.
We denote the condition that $\widehat M$ has the Calabi-Yau metric as 
 ($\bold{CY}$).
\end{assum}

%
%

We can identify the holomorphic volume form.
The holomorphic coordinates on $M$ are given by
 $z_1 = s + ix^*$, $z_2 = t + iy^*$. 
\bpr
The form $dz_1 \wedge dz_2$ determines a globally
 defined holomorphic volume form on $M$.
\epr
\proof
Take another affine coordinates $S, T$ on the base 
 which are related to $s, t$ by $SL(2, \Z)$ transformation.
Take fiber coordinates $X^*, Y^*$ so that
 $(S, T, X^*, Y^*)$ constitute complex coordinates for $M$. 
Namely, $J(\del_S)  =\del_{X^*}$ and $J(\del_T)  =\del_{Y^*}$. 
These satisfy 
$$
\left( \begin{array}{c}
         X^* \\
         Y^* \end{array} \right)
 = \left( \begin{array}{cc}
         a & b \\
         c & d \end{array} \right)
   \left( \begin{array}{c}
         x^* \\
         y^* \end{array} \right)
$$
 for some $\left( \begin{array}{cc}
         a & b \\
         c & d \end{array} \right) \in SL(2, \Z)$. 
We also have 
$$
S = as + bt, T = cs + dt.
$$ 
Putting $Z_1 = S + iX^*$ and $Z_2 = T + iY^*$,
 it is easy to see that
$$
dZ_1 \wedge dZ_2 = dz_1 \wedge dz_2. 
$$\qed

\begin{rem}
As proved,
 $dz_1 \wedge dz_2$ is a holomorphic volume form 
 whatever the metric is. 
When the metric is Calabi-Yau, this form is parallel with
 respect to the metric connection.
In particular, the norm of it is constant in Calabi-Yau case.
\end{rem}

Now the relation between the symplectic and complex coordinates
 can be described using the metric.
In fact, the complex coordinates of the base are given by
 the Legendre transformation of the symplectic coordinates:
$$
\Check s = \frac{\del h}{\del s},\;\; \Check{t} = \frac{\del h}{\del t}.
$$
Define a function
 $\widetilde h$ by $\widetilde h = h - (s \Check{s} + t \Check{t})$
 which is seen as a function of $\Check s, \Check{t}$ by representing $s, t$
 as functions of $\Check s, \Check{t}$.
The symplectic coordinates are described by the inverse Legendre
 transformation:
$$
- s = \frac{\del \widetilde h}{\del \Check s}, \;\;
- t = \frac{\del \widetilde h}{\del \Check t}.
$$

Now we set up the situation with which we will be concerned:\\
%
Let $(\widehat{M}_1, \omega_1, J_1)$
 be a complex 2-dimensional
 K\"ahler torus, $\pi_{1}:\widehat{M}_{1} \ra T_B^2$ a Lagrangian
 torus fibration which satisfies the conditions $1, 2, 3$ and 
 $(\bold{CY})$ above,
 and use the same letters $s, t$, $\Check s, \Check t$, etc... 
 to describe the corresponding notions we have discussed.

We define $J_{\ep}$, $\ep \in (0, 1]$ by requiring
$$J_{\ep}(\del_{\Check s}) = \ep^{-1} \del_{x},
  \;\; J_{\ep}(\del_{\Check t}) = \ep^{-1} \del_{y}$$
 and 
$$J_{\ep}^2 = -1$$
 hold in the complex coordinates of $(\widehat{M}_{1}, \omega_{1}, J_{1})$
 and
$$\omega_{\ep} = \ep \omega.$$ 
These $(J_{\ep}, \omega_{\ep})$ still define
 K\"ahler structures, and satisfy the conditions $1, 2, 3$ and
 $(\bold{CY})$.
Moreover, $s, t$ and $\Check s$, $\Check t$ are symplectic
 and complex coordinates of the base, respectively.
But, the Darboux and complex coordinates of the total space
 are changed, namely, the fiber coordinates must be multiplicated
 by $\ep$.
%
%

We assume we have a family $f_{\ep_{\nu}}:E_{\ep_{\nu}}
                             \ra \widehat{M}_{\ep_{\nu}}$
 of rank 2 stable vector bundles, $\nu \in \Bbb{N}$, which are 
 topologically isomorphic.
We fix a hermitian metric $h$ on them independent of $\nu$.
Then, after complex gauge transformations, each stable bundle admits
 a unique (up to unitary gauge transformations) HYM connection.
We denote this connection by $\Xi_{\ep_{\nu}}$.
%
%
Our main result is the following.
\bt
Let $M$ be the mirror symplectic manifold of $\widehat{M}_1$.
There is a double valued Lagrangian multisection, possibly non-reduced
 and possibly with
 ramifications,
 for $M \ra T_B^2$, determined by the family
 $(E_{\ep_{\nu}}, \Xi_{\ep_{\nu}})$.
The ramifications occur at finitely many points.
%
Moreover, if 
 the first chern class of the bundle is $0$, the
 multisection satisfies the special Lagrangian condition
 on the smooth part. 
\et
\br
As shown in section 7 (see the paragraph before theorem 7.2),
 the ramification
 loci mentioned in this theorem
 are locally diffeomorphic to those of complex curves.
\er
\br 
The multisection may not be unique in general.
In the proof of the theorem, we will consider a converging subsequence.
Since the limit connection may depend on this subsequence,
 the multisection can also depend on them.
\er
The idea is that in the limit $\ep_{\nu} \ra 0$ as $\nu \ra \infty$,
 the curvatures of HYM connections
 will converge to $0$ on each fiber,
 so that flat connections are induced.
Since each fiber is a torus, a flat connection corresponds to a point of
 the symmetric product of the dual torus.
These points constitute the desired Lagrangian multisection.
In the next section, we discuss
 the behavior of HYM equations in the adiabatic limit($\ep_{\nu} \ra 0$). 
\section{Set up and outline}
Let $E$ be a $U(2)$ bundle on a (complex)
 two dimensional K\"ahler torus
 $\pi : \widehat M_1 \ra T^2_B$, with Lagrangian torus fibration
 and K\"ahler $T^2$-structure and Calabi-Yau metric as before.
%
Here $T^2_B$ is the base torus.
As in the last section, we consider a family of such manifolds
 parameterized by $\ep$: $(\widehat{M}_{\ep}, \omega_{\ep}, J_{\ep})$.
Also, we use the same notations 
 $s, t$, $\Check s, \Check t$, etc... as in the last section.
Let
 $\g_E$ be the adjoint bundle of $E$.
Given a hermitian connection on $E$, locally it is written in
 the following form.
\begin{equation}
 \Xi = A + \Phi ds + \Psi dt,
\end{equation}
 where
 $A \in \Gamma(U, \Gamma(T^{\ast}T^2_F \otimes \g_E|_{T^2_F}))$ and
 $\Phi, \Psi \in \Gamma(U, \Gamma(\g_E|_{T^2_F}))$.
$U$ is an open disc of the base $T^2_B$.
$T^2_F$ means a fiber over some $x \in U$.
We can identify all the fibers over $U$ by the flat connection.

Let $F_{\Xi} \in
    \Gamma(\wedge^2 T^{\ast}(\pi^{-1}(U)) \otimes \g_E)$
 be its curvature.
Write the self dual (1, 1) part of $F_{\Xi}$
 with respect to the metric associated to the 
 K\"ahler form $\omega_{\ep}$ as $i\widehat{F_{\Xi}}
 = \Lambda_{\ep} F_{\Xi}$, where
 $\Lambda_{\ep}$ means contraction by $\frac{1}{2}\omega_{\ep}$. 
Then, the HYM equations for the connection $\Xi$
 are
$$
\widehat{F_{\Xi}} = c_{\ep} E_2
$$
 and 
$$
F_{\Xi}^{0, 2} = F_{\Xi}^{2, 0} = 0,
$$
 here $c_{\ep}$ is a constant given by
$$
c_{\ep} = \frac{2\pi deg(E)}{Vol(\widehat M_1, g_{\ep}) rank(E)},
$$
 here
 $deg(E) = \int_M c_1(E) \cup \omega_{\ep}$ and
 $E_2$ is the $2 \times 2$ identity matrix.
%
%
%
%
%
%
%
%
%

Since $\omega_{\ep}$ is linear in $\ep$ and $Vol(M, g_{\ep})$
 is quadratic in $\ep$, $c_{\ep}$ is linear in $\ep^{-1}$,  
 so $c_{\ep} \omega_{\ep}$ does not depend on $\ep$.
We define a constant $c_0$ by $c_0 \omega = c_{\ep} \omega_{\ep}$.  
Therefore, for a HYM connection, $F_{\Xi} - ic_0\omega$ is ASD.

We denote by $F_{mix}$ the terms of the curvature whose differential forms
 have both the fiber and the base directions.
More concretely,
$$
F_{mix} = (d_A \Phi - \del_s A) ds + (d_A \Psi - \del_t A) dt,
$$
 where $d_A$ is the covariant derivatives along fibers.
Note that this is globally well-defined because of the existence
 of the flat structure on the bundle $\widehat{M_1} \to T_B^2$. 

With these notations, HYM connections on $E \to \widehat M_{\ep}$
 satisfy the following equations.
\begin{eqnarray}
\del_t \Phi - \del_s \Psi - [\Phi, \Psi]
             - (det g)^{-\frac{1}{2}}\ep^{-2}
            \ast_{T_F^2}F_A = 0,\\
\ast_{\widehat M_1}(F_{mix} - c_0 \omega) = - (F_{mix} - c_0 \omega),
\end{eqnarray}
here $\ast_{T_F^2}$ is the Hodge operator on $(T_F^2, g|_F)$ (the fibers)
 and $\ast_{\widehat M_1}$ is the Hodge operator on $\widehat M_1$.
The coefficient $det g$ is in fact 1 by the Calabi-Yau condition.

We remark about the relation between $U(2)$-bundles
 and $SU(2)$-bundles. 
The Lie algebra of $U(2)$ decomposes as a direct sum like 
$$u(2) = i \R \oplus su(2).$$
Correspondingly, a connection form $D$ on our $U(2)$-bundle decomposes.
We write it as
$$D = a + D_1,$$
 $a$ is an $i \R$-valued 1-form and
 $D_1$ is an $su(2)$-valued 1-form in a local unitary frame.

The $D_1$-part of the connection defines a connection
 on a different bundle, namely, the $SO(3)$-bundle which is obtained
 as the $su(2)$-direct summand of the adjoint bundle of the principal
 $U(2)$ bundle associated to $E$.
We write it as $Ad E$.

In any case, the restriction of $E$ to a fiber is topologically trivial
 because the fibers are Lagrangian submanifolds and so
 $c_1(T_F^2) = c_0 \omega (T_F^2)  = 0$, and so the restrictions of 
 $Ad E$ to the fibers are also trivial.
It can be seen that $Ad E$ lifts to an $SU(2)$ bundle if and only if
 $c_1$ of $E$ is divisible by two.
When the connection $D$ on $E$ is HYM, the connection $D_1$
 on $Ad E$ is ASD.

The equations (2) and (3) decompose into $i \R$ part and $su(2)$ part.
Namely, we can decompose $A, \Phi, \Psi$ into 
 constant multiples of the identity matrix plus traceless matrices:
$$
A = \left( \left( \begin{array}{ll}
            a_1 & 0 \\
             0  & a_1 \end{array} \right) 
        + X_1 \right) dx 
    + \left( \left( \begin{array}{ll}
            a_2 & 0 \\
             0  & a_2 \end{array} \right) 
        + X_2 \right) dy, 
$$
 with $X_1, X_2 \in su(2)$ and similarly for $\Phi$ and $\Psi$.
We write it as 
$A = A_1 + A_2$, $\Phi = \Phi_1 + \Phi_2$ and
 $\Psi = \Psi_1 + \Psi_2$, where $A_1, \Phi_1, \Psi_1$
 are $i \R$ parts.

Then the equation (2) becomes
\begin{equation}
\del_t \Phi_1 - \del_s \Psi_1 - (det g)^{-\frac{1}{2}}
                \ep^{-2} *_{T_F^2} F_{A_1} = 0 
\end{equation}
and
\begin{equation}
\del_t \Phi_2 - \del_s \Psi_2 - [\Phi_2, \Psi_2]
 - (det g)^{-\frac{1}{2}} \ep^{-2} *_{T_F^2} F_{A_2} = 0.
\end{equation}
(3) becomes
\begin{equation}
\begin{array}{l}
 \ast_{\widehat{M_1}} (F_{mix})_1 + (F_{mix})_1 \\
 = \ast_{\widehat{M_1}} ((d_F \Phi_1 - \del_s A_1)ds
         + (d_F \Psi_1 - \del_t A_1)dt) \\
   \hs{.1in}  + (d_F \Phi_1 - \del_s A_1)ds
         + (d_F \Psi_1 - \del_t A_1)dt \\
 =   c_0 E_2 \omega, 

\end{array}
\end{equation}
and
\begin{equation}
\begin{array}{l}
 *_{\widehat{M_1}}
 ((d_{A_2} \Phi_2 - \del_s A_2)ds + (d_{A_2} \Psi_2 - \del_t A_2)dt) \\
 = 
 - ((d_{A_2} \Phi_2 - \del_s A_2)ds + (d_{A_2} \Psi_2 - \del_t A_2)dt).
\end{array}
\end{equation}
Here in (6) $d_F$ means the exterior differential in the fiber direction.

%

We want to know the behavior of the connections 
 $\Xi_{\ep_{\nu}}$ as $\nu \to \infty$. 
The above equations are used for this analysis.
The conclusion is that after removing bubbles and taking subsequences,
 the connections converge to some limit connection
 which is flat in the direction of the fibers.
We treat the convergence of the $su(2)$-part of the connections
 in the next three sections and we treat the $i \R$-part
 in section 7.

%
%

%
%
%
Our main analytic result (which is for the $su(2)$-part)
 is summarized as follows.
We consider the bundles $Ad E_{\ep_{\nu}}$ and we rewrite
 $\Xi_{\ep_{\nu}}$ as a sequence of ASD connections on them satisfying
 (5) and (7).
From the following, the main theorem 2.6 will be proved.
\bt
Suppose we are given a family $(Ad E_{\ep_{\nu}}, \Xi_{\ep_{\nu}})$ as above,
 that is, a family of ASD connections on $SO(3)$ bundles
 $Ad E_{\ep_{\nu}} \to \widehat{M}_{\ep_{\nu}}$ of fixed topological
 type and hermitian metric,and trivial over the fibers.
As a bundle, $Ad E_{\ep_{\nu}}$ does not depend on $\nu$,
 and we denote it by $Ad E$.
Then there is a countable subset $S \subset T_B^2$ with 
 finite accumulation points satisfying the following property.
There is a subsequence of $\Xi_{\ep_{\nu}}$ (which we denote by the same
 letter to avoid too many subscripts),
 a sequence of gauge transformations $g_{\ep_{\nu}}$
 and a unitary connection $\Xi_0$ on the bundle
 $Ad E|_{\pi^{-1}({T_B^2 \backslash S})}$,
 here $\pi: \widehat{M}_1 \to T_B^2$,
 such that
 on any compact subset $K \subset \widehat{M}_1$
 of $\pi^{-1}(T^2_B \backslash S)$, 
 $g^*_{\ep_{\nu}}(\Xi_{\ep_{\nu}})$ converges to $\Xi_0$
 in $C^{\infty}$ sense.

$\Xi_0$ is flat in the direction of fibers and
 $\Xi_0$ satisfies the equation 
$$
\ast_{\widehat M_1}F_{mix} = - F_{mix}.
$$
\et
\proof{}
%
%
First we recall the classification of bubbles in our situation,
 by Dostoglou and Salamon in the proof of theorem 9.1, $\cite{DS}$
 or proof of theorem 3.1 of its erratum.
In the same way as in $\cite{DS}$, there are three types of bubbles:
\begin{quotation}
\nnn
(1). Instantons on $S^4$\\
(2). Instantons on $\Bbb{R}^2 \times T^2$\\
(3). Holomorphic spheres on $Rep_{T^2}(SO(3))$
\end{quotation}
Here $Rep_{T^2}(SO(3))$ is the representation variety of $\pi_1(T^2)$
 in the gauge group $SO(3)$.
For reader's convenience, we briefly recall the nature of these bubbles.
%

These three types of bubbles are classified due to the
 degree of divergence of the curvature.
Namely, let $\Xi_{\ep_{\nu}}$ be a sequence of ASD connections on 
 the bundles $Ad E_{\ep_{\nu}} \ra 
 \widehat M_{\ep_{\nu}}$, $\ep_{\nu} \ra 0$ as
 $\nu \ra \infty, \nu \in \Bbb{N}$.
Let $w_{\ep_{\nu}}$ be a sequence of points on $T^2_B$,
 converging to $w_0 \in T^2_B$ 
 such that 
$$
c_{\nu} = c_{\nu}(w_{\nu}) 
  = \ep_{\nu}^{-1} \pll F_{\ep_{\nu}} \pll_{L^{\infty}(T_F^2)}
       + \pll \del_t A_{\ep_{\nu}} - d_{A_{\ep_{\nu}}}
          \Psi_{\ep_{\nu}} \pll_{L^{\infty}(T_F^2)}
$$
  diverges.
%

The case (1) appears when the sequence $\ep_{\nu} c_{\nu}$
 is unbounded. In this case, a usual instanton on $S^4$ splits off.
The case (2) is the case when there exists a sequence $w_{\nu} \ra w_0$
 such that  $\ep_{\nu} c_{\nu} (w_{\nu}) \geq \delta > 0$ and not
 of the case (1). Here, an anti self dual connection on $\Bbb{R}^2 \times T^2$
 appears when the rescaling argument is applied.
Finally, the case (3) is the case when there exists a sequence
 $w_{\nu} \ra w_0$ on $T_B^2$ such that $c_{\nu}(w_{\nu})$
 diverges, but $\ep_{\nu} c_{\nu} (w_{\nu})$ tends to zero.
%
 
The case (1) is treated exactly in the same way as in $\cite{DS}$.
In this case, the energy of the bubble is
 bounded from below by the Yang-Mills energy of 1-instanton on $S^4$.
%
In the case (2), the exponential decay estimate proved in theorem 7.4 of
 $\cite{DS}$ is not applicable here
 (by the existence of reducible connections). 
However, we prove that this type of bubbles also brings away
 finite amount of
 energy bounded from below by a positive constant in section 5.
But this bound depends on the boundary condition.
This causes the countableness of the set $S$ defined in section 6. 

On the other hand, the sequence $\Xi_{\ep_{\nu}}$ induces a
 sequence of maps $\phi_{\nu}$ from $T^2_B$ to
 $Rep_{T^2}(SO(3))$.
See the paragraph before remark 6.2.
This is holomorphic when we apply hyperK\"ahler rotation.
By minimal surface theory, the sequence $\phi_{\nu}$
 converges to a limit map $\phi$ modulo finite number of
 points where the sequence develops bubbles (proposition 6.3).
The set $S$ is, roughly (and incorrectly),
 the union of the points on the base
 over which the first and the second types of bubbles or
 the bubbles of the sequence $\phi_{\nu}$ occur
 (precise definition is given measure theoretically
 (definitions 6.4 and 6.8).
The precise relation between measure theoretic divergence
 and bubbles of connections or holomorphic maps
 is interesting and left unclear.
 See propositions 6.27, 6.28 and the paragraph after them).

In section 6, the $C^{\infty}$ convergence
 (modulo gauge transformation and taking subsequences)
 of the connections $\Xi_{\ep_{\nu}}$ is proved (theorem 6.6). 
Then the latter assertion of theorem 3.1 is clear from this $C^{\infty}$
 convergence. \qed\\

Moreover, the map defined by 
$$
x \mapsto [A_x],
$$
 where $A_x$ is the fiber component of the limit connection $\Xi$
 on the fiber $\pi^{-1}(x)$ and $[A_x] \in Rep_{T^2}(SO(3))$
 denotes the complex gauge equivalence class (see the next section for
 this terminology),
 coincides with $\phi|_{T^2_B \setminus S}$.

\section{The space of unitary connections on $T^2$}
In this section, we study the analytic nature of the space of 
 traceless unitary connections on a rank 2, topologically
 trivial, hermitian vector bundle $E$ on $T^2$.
This is used in sections 5 and 6 to deal with the bubbling analysis
 in the reducible case.

First we remark about the relations between $SU(2)$ and $SO(3)$ bundles
 in our situation.
Our $SO(3)$ bundles in the last section and the whole of
 (the main stream of) this paper
 are trivial on each fiber and so lift to $SU(2)$-bundles there.
In this section we will work with this $SU(2)$ bundle.
There are very few differences between them, and
 the space of connections and the action of
 the complexified gauge group (see the paragraph before proposition 4.6)
 are in fact the same.
Moreover, results in this section are applied later to situations
 (in four dimension)
 where we can actually lift our $SO(3)$ bundle to an $SU(2)$ bundle.

The only reason we choose $SU(2)$ bundle is that there is
 convenient language (from complex geometry)
 on this side.

We will later be concerned with small energy
 ASD connections on a rank two bundle on $S \times T^2_F$.
$S$ is a smooth surface with a metric (compact or non compact)
 considered as the base of the fibration. 
While the K\"ahler structure in the last section was such that
 the base $S$ and $T^2_F$ are Lagrangians, we take in this section a
 complex structure such that $T^2_F$ are complex submanifolds
 (however, we do not need the total space to be a complex manifold,
 in particular, this family of $T^2_F$ is only supposed to be
 a smooth family).
This is again only for terminological convenience,
 to give an outlook of the relation between the moduli of flat
 connections and the holomorphic classification
 of bundles (these argument appear until proposition 4.5).
After that, complex structure is not required in this section at all.

It is known that (see theorem 5.1 of the next section)
 small energy ASD connections have gauges in which
 the connection matrices have small $C^r$-norms, so in particular,
 the curvature of the connection restricted on any of the fibers
 is also small.
This means the holomorphic bundle  
 given by restricting $E$ to any of the fiber is semi-stable
 (see Fukaya $\cite{F1}$, or it is easy to see directly in this case,
 because an unstable degree zero bundle on $T^2$ is a direct sum of
 a line bundle with nonzero $c_1$ and its inverse
 (see Friedman \cite{Fr}, chapter 2, theorem 6), so it cannot have
 a connection of small curvature inducing
 the required holomorphic structure). 
So we concentrate on semi-stable bundles on $T^2$ here.

Every (small energy) ASD connection
 determines a holomorphic structure on $E$
 (restricted to fibers),
 and so the holomorphic classification of such bundles
 will be helpful to give some idea how the space of flat connections 
 looks like.
We first recall this. For the proof, see Friedman $\cite{Fr}$, theorem 25
 in chapter 8.
\bt
Any semi-stable holomorphic structure on a degree zero rank two bundle $E$ 
 on $T^2$ is one of the following form: \\
1. $E = \lam \oplus \lam^{-1}$, $\lam$ is a line bundle of
 degree zero. \\
2. $E = \mathcal{E} \otimes \lam$, where $\mathcal{E}$ is the
 unique nontrivial extension of $\mathcal{O}_{T^2}$ by itself,
 and $\lam$ is  one of the four line bundles on $T^2$ such that
 $\lam^{\otimes 2} = \mathcal{O}_{T^2}$.\\
The moduli space of S-equivalence classes of such bundles
 is isomorphic to the Riemann sphere $S^2$.
\et
Here we do not give the definition of S-equivalence,
 because it is only subsidiary to our argument.
But we give a picture below what it means in our cases
 (remark 4.3).
See \cite{Fr}, page 154 for precise definition.
\begin{rem}
Degree zero, rank two holomorphic bundles are always 
 topologically trivial, since they allow $SU(2)$-structure
 and any $SU(2)$ bundle over $T^2$ is topologically trivial.
There is no stable bundle and unstable bundles are 
 the direct sum of line bundles of non-zero degree 
 which are mutually inverse, as mentioned
 above.
\end{rem} 
The last assertion of the theorem can be interpreted from the point
 of view of connections.
We denote the space of unitary connections by
$$
\mathcal{A} = \left\{ \left( \begin{array}{cc}
                          a & b \\
                          c & -a \end{array} \right)d\overline{z} -
                      \left( \begin{array}{cc}
                          a & b \\
                          c & -a \end{array} \right)^{\dag}dz
              \big| a, b, c \in C^{\infty}(T^2, \C)\right\}
$$
In particular, we fix a gauge in which the product connection
 is specified by the connection represented
 by 0-matrix in $\mathcal A$.
All the norms are estimated using this gauge (see also remark 5.2).
There is an affine subspace consisting of connections with 
 constant components:
$$
sl_2 \C = \left\{ \left( \begin{array}{cc}
                         a & b \\
                          c & -a \end{array} \right)d\overline{z} -
                  \left( \begin{array}{cc}
                         a & b \\
                          c & -a \end{array} \right)^{\dag}dz
               \big| a, b, c \in \C \right\}
$$ 
Any flat connection on $E$ is gauge equivalent to a connection 
 of diagonal constant form.
We denote the space of these connections by $\T$ :
$$
\T = \left\{ \left( \begin{array}{cc}
                        a & 0 \\
                        0 & -a \end{array} \right)d\overline{z} -
             \left( \begin{array}{cc}
                        a & 0 \\
                        0 & -a \end{array} \right)^{\dag}dz
                     \big| a \in \C \right\}
$$
Connections in $\T$ which equal modulo
$$
\left\{ \left( \begin{array}{cc}
                  in\pi + m\pi & 0 \\
                      0        & -in\pi -m\pi \end{array}\right)d\overline{z}
                            -
        \left( \begin{array}{cc}
                  in\pi + m\pi & 0 \\
                      0        & -in\pi -m\pi \end{array}\right)^{\dag}dz
                 \big| m, n\in \Z \right\}
$$
 are gauge equivalent
 by gauge transformations 
$$
g = \left( \begin{array}{cc}
         e^{2\pi i(nx - my)} & 0 \\
                0            & e^{-2\pi i(nx - my)}
          \end{array} \right),
$$
 and thus, we have a family of gauge
  equivalent connections parametrized by $T^2$.
Further, the exchange of the diagonal
 components (this is given by the matrix 
 $\left( \begin{array}{cc}
            0 & 1 \\
           -1 & 0 \end{array} \right)
 $) also gives a pair of gauge equivalent connections.
This last action has four fixed points on $T^2$ and the resulting space is
 the Riemannian sphere.

Comparing with the holomorphic classification, these
 connections define the type 1 of the above theorem.
Then, how can we interpret the type 2 bundles? 
In fact, these do not permit flat connections, and connections
 of the form 
$$
\left\{ \left( \begin{array}{cc}
                  \alpha & b \\
                      0  & \alpha \end{array}\right)d\overline{z} -
           \left( \begin{array}{cc}
                  \alpha & b \\
                      0  & \alpha \end{array}\right)^{\dag}dz                               \big| b \in \C^*, \alpha = \frac{\pi}{2},
             \frac{i\pi}{2}, \frac{(1+i)\pi}{2} \right\}
$$
 or
$$
\left\{ \left( \begin{array}{cc}
                      0 & b \\
                      0  & 0 \end{array}\right)d\overline{z} -
           \left( \begin{array}{cc}
                      0 & b \\
                      0  & 0 \end{array}\right)^{\dag}dz                               \big| b \in \C^* \right\}
$$
 give the corresponding holomorphic structures (any $b$ gives the
 same holomorphic structure and in particular, these ($\alpha$ fixed)
 are all complex gauge equivalent).
The former three cases are not trace free and do not
 appear in our analysis.
\begin{rem}
The S-equivalence identifies these with the limit ($b \to 0$) flat
 connections, and so the moduli of S-equivalent classes
 is isomorphic to $S^2$.
 \end{rem}
\br
In the case of trivial $SO(3)$ bundle,
 since the map $su(2) \to so(3)$ corresponding to the double
 cover $SU(2) \to SO(3)$ is given by multiplying by two on $\T$,
 everything is rescaled by $\frac{1}{2}$ and
 the moduli of flat connections is still $S^2$, with four singular points.
\er
We summarize above argument in the next proposition.
\bpr
There is a positive constant $C$ such that
 if $A \in \mathcal{A}$ has curvature
 $\pll F_A \pll_{L^2(T^2)} \leq C$, 
 then $A$ is complex gauge equivalent (see below)
 to either 
 $$
1. \left( \begin{array}{cc}
              a & 0 \\
              0 & -a \end{array} \right)d\overline{z} -
   \left( \begin{array}{cc}
              a & 0 \\
              0 & -a \end{array} \right)^{\dag}dz
         , \, a \in [0, \pi) \times i[0, \pi)
 $$
or
$$
2. \left( \begin{array}{cc}
           0 & 1 \\
           0    & 0 \end{array} \right)d\overline{z} -
      \left( \begin{array}{cc}
           0 & 0 \\
           1    & 0 \end{array} \right)dz.    
$$
\epr
\proof
Recall that any unitary connection of $E$ defines a holomorphic structure.
It follows from remark 4.2 and the previous argument that
 if the curvature is small, the induced
 holomorphic structure is isomorphic to the one which is induced
 by some of the connections of the statement.
Given any connection with small curvature, take a biholomorphic bundle map from it
 to one of the above. 
This map gives the desired complex gauge transformation. \qed \\

Our main objective is the nature of the action of complex gauge
 group on flat connections. 
We are mainly (exclusively in this and the next sections)
 concentrate to connections satisfying the assumption of proposition 4.5
We first recall the complex gauge transformation
$$
\mathcal G^{\C} = \left\{ g \in C^{\infty}(T^2, SL(2, \C)) \right\},
$$ 
 here
          $g$  acts  by
$$
 g^*A = 
          g^{-1} \overline{\del}g + g^{-1}A^{0, 1}g -
             (g^{-1} \overline{\del}g + g^{-1}A^{0, 1}g)^{\dag}.
$$ 
The Lie algebra of $\mathcal{G}^{\C}$ is given by
$$
\mathcal U_{\C} = \left\{ u \in C^{\infty}(T^2, sl_2 \C) \right\}
$$
We remark the following fact.
\bpr
The isotropy groups in the group of complex gauge transformations
 of the connections of the proposition 4.5 are
 described as follows.\\
For the case 1 and $a \neq 0, \frac{\pi}{2}, \frac{i\pi}{2},
           \frac{(1 + i)\pi}{2}$, then
            $ \left( \begin{array}{cc}
             \alpha & 0        \\                        
                0   & \alpha^{-1} \end{array}
                             \right), \alpha \in \C^* $  \\
For the case 1 and $a = 0, \frac{\pi}{2}, \frac{i\pi}{2},
           \frac{(1 + i)\pi}{2}$,
 then $SL(2, \C)$\\
For the case 2, then 
            $ \left( \begin{array}{cc}
             \pm 1 & \alpha        \\                        
                0   & \pm 1 \end{array}
                             \right), \alpha \in \C $ 
\epr
\proof
This follows from a straightforward calculation. 
For the second case, the isotopy
 group of, for example, $a = \frac{\pi}{2}$ case is given by
$$
g = \left( \begin{array}{cc}
                    p & q e^{2\pi iy} \\
      r e^{-2\pi iy}  & s \end{array} \right),
      \left( \begin{array}{cc}
                    p & q \\
                    r & s \end{array} \right) \in SL(2, \C),  
 $$
 where $y$ is one of the standard coordinates on $T^2$.  
 \qed \\

In particular, the whole locus of the moduli space of 
 flat connections is degenerate
 in the sense that each point has positive dimensional isotropy group.
In other words, the moduli problem is non-transversal. In fact, the index
 of the complex
$$
0 \ra \wedge^0 T^*T^2 \otimes su(2) \stackrel{d^0_A}{\ra}
    \wedge^1 T^*T^2 \otimes su(2)
       \stackrel{d^1_A}{\ra} \wedge^2 T^*T^2 \otimes su(2) \ra 0
$$
here $A$ is a flat connection,
 is zero, and so the virtual dimension of
 the moduli of flat connections is zero.
%
On this point, we make some remarks.
\br
1. For connections of the form 
      $$ A = \left( \begin{array}{cc}
                a & 0        \\                        
                0 & -a \end{array} \right)d\bar{z}
          - \left( \begin{array}{cc}
                a & 0        \\                        
                0 & -a \end{array} \right)^{\dag}dz \in \T,
  a \neq \frac{m\pi}{2} + \frac{in\pi}{2}, m, n \in \Z$$  
  $ker d^0_A$ and $ker (d^1_A)^*$ are 1-dimensional,
 and correspondingly, the space $\bigwedge^1 T^*T^2 \otimes su(2)$
 has a 2-dimensional subspace of $d_A$-harmonic forms. 
 This 2-dimensional space is precisely the tangent space to the moduli
 of flat connections. \\

\nnn
2. In particular, those complex gauge transformations which are obtained by 
 exponentiating 
$$
\mathcal U_{\C}^{\perp} = \left\{ \left(
          \begin{array}{cc}
          \alpha & \beta \\
          \gamma & -\alpha \end{array} \right) \in \mathcal{U}_{\C} 
          \big| \int_{T^2} \alpha d\mu = 0, \alpha, \beta, \gamma \in 
             C^{\infty}(T^2, \C) \right\}
$$
acts on $A$ without isotropy, at least locally.
\er
We make the second remark in a more precise form.
\bpr
For a connection $$ A = \left( \begin{array}{cc}
                a & 0        \\                        
                0 & -a \end{array} \right)d\bar{z}
          - \left( \begin{array}{cc}
                a & 0        \\                        
                0 & -a \end{array} \right)^{\dag}dz \in \T, 
 a \neq \frac{m\pi}{2} + \frac{in\pi}{2}, m, n \in \Z$$ 
 there is a neighbourhood $U$ in $\T$ and a neighbourhood
 $V$ of $0 \in \mathcal U_{\C}^{\perp}$ (with $C^{r+1}-topology$.
Here we take $\alpha, \beta, \gamma \in C^r(T^2, \C)$.
We neglect these obvious modifications hereafter)
 such that the map
$$
\phi: U \times V \ra \mathcal A,\;\;\; (A', v) \mapsto (\exp(v))^*A'
$$
is a local diffeomorphism ($\mathcal A$ is given the $C^r$-topology).
\epr 
\proof
The tangent space of $\mathcal{A}$ at $A'$ is 
$$
\begin{array}{ll}
T_{A'}\mathcal{A} & = \Omega^1(T^2, su(2)) \\
                  & = Im d_{A'}^0 \oplus Im(d_{A'}^1)^* \oplus \mathcal{H},
\end{array}
$$
 here $\mathcal{H}$ is the space of $su(2)$-valued harmonic one forms
 with respect to $d_{A'}$ and $d_{A'}^*$.
As in remark 4.7.1, we have 
$$
T_{A'}U = \mathcal{H}.
$$
The tangent space of the orbit of the complex gauge transformations
 at $A'$ is given by 
$$
Im d_{A'}^0 \oplus Im (d_{A'}^1)^*.
$$
So $\phi$ is locally surjective, that is,
 for any $p \in Im\phi$, there exists $W$, a neighbourhood of $p$ such that
 $W \subset Im\phi$. 
To prove the injectivity, it suffices to prove $d_A^0$ and $(d_A^1)^*$ do not
 have kernel in $\mathcal U_{\C}^{\perp}$.
However, since $ker d_A^0$ and $ker (d_A^1)^*$ correspond
 to the isotropy groups of $A$ in the complex gauge transformations,
 the proposition follows from the classification of the isotropy.\qed\\

The tangent space of $\mathcal U_{\C}^{\perp}$
 at $0$ is divided to unitary and hermitian parts.
Let $V_1$ be a neighbourhood of $0$ in the space of 
 sections of traceless, skew hermitian $2 \times 2$ matrices over $T^2$
 with the integral of diagonal entries $0$ and
 $V_2$ be the neighbourhood of $0$ of similar
 space of hermitian matrices.
Then the following is clear.

\bc
Let us take $U$ as in the proposition.
The map 
$$
\psi: U \times V_1 \times V_2 \to \mathcal A, \;\; (A', u, v) \mapsto
 (\exp (u))^* (\exp (v))^* A'
$$
 is a local diffeomorphism. \qed
\ec

We will give a refinement of this proposition later (corollary 4.15).
 At this stage, we give some more remarks.
\br
1. Near the points
   $ A = \left( \begin{array}{cc}
                a & 0        \\                        
                0 & -a \end{array} \right)d\overline{z} -
         \left( \begin{array}{cc}
                a & 0        \\                        
                0 & -a \end{array} \right)^{\dag}dz$,
              $a = \frac{im\pi}{2} + \frac{n\pi}{2}, n, m \in \Z$,
 the structure of $\mathcal{A}$ is more complicated than near other
 points of $\T$.
Among these, we will concentrate on the analysis of the neighbourhood of
 $A = 0$, the rest can be handled in the same manner. \\

\nnn
2. As the connections get closer to $0$ in $\T$, the
 differential of the diffeomorphism $\phi$ of the above proposition 
 will come close to having a kernel. \\

\noindent
3. However, for the gauge transformations obtained by exponentiating
 the elements of
$$
\mathcal U_{\C}^{\amalg} = \left\{ \left(
         \begin{array}{cc}
          \alpha & \beta \\
          \gamma & -\alpha \end{array} \right)
        \big| \int_{T^2} \alpha d\mu
        = \int_{T^2} \beta d\mu = \int_{T^2} \gamma d\mu = 0
        \right\},
$$  
 the differential does not degenerate.
\er
We will state this point in the following form.
\bpr
Let $A \in \T$ be a connection near $0$. 
Let $U$ be a small neighbourhood of $A$ in $sl_2\C$. 
Let $V$ be a small neighbourhood of $0$ in $\mathcal U_{\C}^{\amalg}$
 in the $C^{r + 1}$-topology. 
Then, the map 
$$
U \times V \ra \mathcal A, \;\;\; (A', v) \mapsto (\exp(v))^*A'
$$
 is a local diffeomorphism and there is a constant $C > 0$,
 not depending on $A$  such that
$$
\pll (\exp(v))^*A' - A' \pll_{C^r(T^2)} \geq C \pll v \pll_{C^{r + 1}(T^2)}
$$
\epr
\br
In this proposition and the followings, to simplify the notation 
 we write various constants by $C$ when no confusion can occur.
Their important dependences or non-dependences on the parameters
 will be noted in each case. 
\er
\br
The term `near 0' here means `not too large', in view of the above
 remark 4.10.1.
And it does not mean `very small compared with 1'.
\er
\noindent
\it{Proof of proposition 4.11.}
\rm
Modulo the terms containing the quadratics or higher order terms of $v$, 
 $(\exp(v))^*A'$ is given by
$$
(\exp(v))^*A' = A' + d_{A'}^0 v - (d_{A'}^1)^* (\ast v),
$$
where $\ast$ is the Hodge operator of $T^2$.  
If $U$ is not large, $d_{A'}^0$ and $(d_{A'}^1)^*$
 do not have kernels in $\mathcal U_{\C}^{\amalg}$. 
The proposition is a standard result of elliptic estimates. \qed\\

\bc
Let $A$ as above. 
Let $V_1$ be a neighbourhood of zero
 in the space of traceless skew hermitian  $2 \times 2$ matrix
 functions 
 over $T^2$
 with the integrals of the entries are zero.
Let $V_2$ be a similar space of traceless hermitian matrix 
 functions.
There is a constant $C$ not depending on $A$
 with the following property.
For any sufficiently small positive constant $\ep$,
 take an $\ep$-neighbourhood $W$ of $A$ in $\mathcal A$. Then, any element
 of $W$ is transformed into $sl_2 \C$ by a unique complex gauge
 transformation $g$ of the form 
 $\exp (v_1) \exp (v_2)$, $v_1 \in V_1, \; v_2 \in V_2$, such that
$$
\pll v_1 \pll_{C^{r + 1}(T^2)} < C\ep ,\;
 \pll v_2 \pll_{C^{r + 1}(T^2)} < C\ep
$$
 hold.  \qed
\ec

\bc
Let $A' \in \T$ is a connection near 0 but not equal to 0
 and let $W$ be a neighbourhood of $A'$ 
 in $\mathcal{A}$, which is sufficiently small compared with
 the distance between $A'$ and $0$.
Let $A \in W$ and assume it is 
 complex gauge equivalent to $A'$.
Then there exists a complex gauge transformation $g$ satisfying
 $g^* A = A'$ with 
$$
\pll A' \pll_{C^r(T^2)} \pll g - Id \pll_{C^{r+1}(T^2)}
         < C\pll A - A' \pll_{C^r(T^2)}.
$$
\ec
\proof
By proposition 4.11, for the gauge transformations obtained by
 exponentiating elements of $\mathcal U^{\amalg}_{\C}$
 certainly the estimate is satisfied, even without the factor
 $\pll A' \pll_{C^r(T^2)}$ of the left hand side. 
Since we can transform $A$ to $sl_2 \C$ by
 these transformations, we want an estimate
 for gauge transformations which are obtained by exponentiating
 constant matrices
$$
\left( \begin{array}{cc}
        0 & b \\
        c & 0 \end{array} \right), b, c \in \C.
$$
Since 
$$
\big[ \left( \begin{array}{cc}
         a & 0 \\
         0 & -a \end{array} \right), 
 \left( \begin{array}{cc}
         0 & b \\
         c & 0 \end{array} \right) \big]
 = 
  \left( \begin{array}{cc}
          0 & 2ab \\
       -2ac & 0 \end{array} \right),
$$
the estimate follows. \qed 
\bpr
For any constant $\ep$ with $0 < \ep \ll 1$, the following statement holds.
 Let $A$ be a connection on $E$ which has a matrix representation
 in the given gauge such that there is a connection 
$A_0 = \left( \begin{array}{cc}
          a_0 & 0 \\
            0 & -a_0 \end{array} \right)d\overline{z} -
       \left( \begin{array}{cc}
          a_0 & 0 \\
            0 & -a_0 \end{array} \right)^{\dag}dz \in \T$ near $0$
 with $\pll A - A_0 \pll_{C^r(T^2)} < \ep$.
Then, one (and only one) of the followings holds.\\
1. $A$ is not complex gauge equivalent to a flat connection. \\
2. $A$ is complex gauge equivalent to a flat connection
 $A_1 \in \mathfrak{t}$
 which satisfies $\pll A_1 - A_0 \pll_{C^r(T^2)} < C\ep$,
 where $C$ is a constant not depending on $A$ and $\ep$.
(since $A_1$ and $A_0$ have constant components, the norm is
 in fact the same as $C^0$-norm) 
\epr
\proof
Assume that $A$ is complex gauge equivalent to a flat connection.
By corollary 4.14, $A$ is transformed into $sl_2 \C$
 by a complex gauge transformation $g$ such that 
$
\pll g - Id \pll_{C^{r + 1}(T^2)} < C\ep.
$  
 So the resulting connection $A'$
 with constant components can be written in the
 following form:
$$
A' = \left( \begin{array}{cc}
         a + a_0 & b \\
            c    & -a -a_0 \end{array} \right)d\overline{z} -
     \left( \begin{array}{cc}
         a + a_0 & b \\
            c    & -a -a_0 \end{array} \right)^{\dag}dz
$$
 with $|a|, |b|, |c| < C \ep$. 
We want to get a connection which is complex gauge equivalent to $A'$
 and have the form 
$ \left( \begin{array}{cc}
            \alpha & 0 \\
               0   & -\alpha \end{array} \right)d\overline{z} -
  \left( \begin{array}{cc}
            \alpha & 0 \\
               0   & -\alpha \end{array} \right)^{\dag}dz.
$
 $\alpha$ is given by the eigenvalue 
 $\lam = \pm \sqrt{(a + a_0)^2 + bc}$ of the matrix
$ \left( \begin{array}{cc}
         a + a_0 & b \\
            c    & -a -a_0 \end{array} \right).$\\
(i) The case $|a_0| > 10C\ep.$\\

Preferring the plus signature,
 $\lam = a_0 \sqrt{1 + \frac{2a}{a_0} + \frac{(a^2 + bc)}{a_0^2}}$. 
So 
$$|\lam - a_0| < |a_0| (|\sqrt{1 + \frac{2a}{a_0}
           + \frac{a^2 + bc}{a_0^2}} - 1|)
                 < |2a| + \frac{|a^2 + bc|}{|a_0|} < 3C\ep.$$
(ii) The case $|a_0| < 10C\ep.$\\

Obviously, $|\lam| < 12C\ep$. \qed 
\bc
Let $A_s$, $s \in S$, $S$ is a compact smooth manifold diffeomorphic to
 a two dimensional disc with several holes,  
 be a smooth family of connections on $E$
 satisfying the assumptions of proposition 4.16. 
Furthermore, assume that  every $A_s$ is complex gauge
 equivalent to a flat one
 and denote by $A_{0, s} \in \T$
 the possibly double valued smooth family of
 flat connections to which $A_s$ is close to (as in the assumption of 
 proposition 4.16) and assume that no $A_s$ is complex gauge equivalent to zero.
Then, there is a possibly double valued
 smooth family $g_s$ of complex gauge transformations such that
 each $g_s^*A_s$ is in $\mathfrak{t}$ and satisfies the estimate
$$
\pll g^*_sA_s - A_{0,s} \pll_{C^r(T^2)} < C\ep.
$$ 
Here $g_s$ acts fiberwise (that is, not as a transformation
 of a bundle on $S \times T^2$)
 with $C$ a constant not depending on $\{ A_s \}$.
\ec
\proof
Since we can transform $A_s$ into $sl_2 \C$ by a complex gauge transformation 
 of the form $\exp(e)$, $e \in \mathcal{U}_{\C}^{\amalg}$ 
 (it is unique, by the assumption that $A_s$ is near to $\mathfrak{t}$), 
 the problem is whether we can diagonalize a family of matrices
 $A(s) \subset sl_2 \C, \; s \in S$ by a family of constant complex transformations.
Let $H$, $G$ be matrices which diagonalize $A(s)$ to the same diagonal matrix,
$$
HAH^{-1} = GAG^{-1} =
 \left( \begin{array}{ll}
          a & 0 \\
          0 & -a \end{array} \right), a \neq 0.
$$
(Note that we can also diagonalize $A(s)$ to
 $
 \left( \begin{array}{ll}
          -a & 0 \\
          0 & a \end{array} \right).
 $
Namely, multiplying $H$ or $G$ by
  $
 \left( \begin{array}{ll}
          0 & i \\
          i & 0 \end{array} \right),
 $
 we get a diagonalizing matrix which transforms $A(s)$ to
 $
 \left( \begin{array}{ll}
          -a & 0 \\
          0 & a \end{array} \right).
 $
This produces the double valuedness.)
Then $G H^{-1}$ commutes with the diagonal matrix
 $ \left( \begin{array}{ll}
          a & 0 \\
          0 & -a \end{array} \right).$
Denoting the commutator subgroup of 
 $ \left( \begin{array}{ll}
          a & 0 \\
          0 & -a \end{array} \right)$
 by $G_0$,
 the element $G$ can be written as 
$$
G  = g_1 H
$$
 where $g_1 \in G_0$.
$G_0$ is isomorphic to $\C^*$ since $a \neq 0$.
Consider the space $W$ with the map $\pi: W \to S$ whose fiber
 over $s$ is the diagonalizing matrices in $SL(2, \C)$ of $A(s)$.
We give it the subspace topology of 
 $S \times SL(2, \C)$.  
The fiber of it has two components
 each of which is diffeomorphic to $\C^*$.
So, we can find a possibly double valued section of it by the assumption that 
 $S$ is diffeomorphic to a disc with several holes. \qed
\br
Note that the double valuedness occurs when the 
 complex gauge equivalence classes of the family $A(s)$ 
 rotate around $0 \in \T$.
In our application in the following sections,
 we will meet only such situations that this does not occur.
\er
\bpr
There is a constant $\delta_0$ with the following property. 
Let $A$ be a connection on $E$ with
 $\pll F_A \pll_{C^r(T^2)} < \delta <\delta_0$.
 Moreover, assume there is a flat connection 
$A_1 = \left( \begin{array}{cc}
              a_0 & 0  \\
              0 & -a_0 \end{array} \right)d\overline{z} -
         \left( \begin{array}{cc}
              a_0 & 0  \\
              0 & -a_0 \end{array} \right)^{\dag}dz \in \T$ 
 with $\pll A - A_1 \pll_{C^r(T^2)} < \ep$,
 $\ep$ is taken as in proposition 4.16. 
Then, one and only one of the followings holds.\\
1. $A$ is not complex gauge equivalent to a flat one.\\
2. We can transform $A$ by a unitary gauge transformation $h$ so that
 $$\pll h^*A - A_0 \pll_{C^{r+1}(T^2)} < C \sqrt{\delta},$$
here $A_0 = \left( \begin{array}{cc}
              a & 0  \\
              0 & -a \end{array} \right)d\overline{z} -
             \left( \begin{array}{cc}
              a & 0  \\
              0 & -a \end{array} \right)dz$
 is a flat connection to which $A$ is complex gauge equivalent
 and $C$ is a constant not depending on $A$ and $A_1$.
\epr
\proof
Assume $A$ is complex gauge equivalent to a flat one. 
By proposition 4.16, we can replace $A_1$ with $A_0$, that is,
 we can assume from the first that $A$ is complex gauge equivalent to $A_1$
 with estimates  
$$\pll F_A \pll_{C^r(T^2)} < \delta <\delta_0, \,\,
 \pll A - A_1 \pll_{C^{r}(T^2)} < \ep.$$
By corollary 4.14, there is a complex gauge transformation
 $g$
 which maps $A$ into $sl_2 \C$ and satisfies the estimate 
$$\pll g - Id \pll_{C^{r + 1}(T^2)} < C \ep.$$ 
In particular, the estimate
$$\pll g^*A - A_1 \pll_{C^r(T^2)} < C\ep$$ holds.
 
We can decompose $g$ into the form $g = g_1g_2$, where
 $g_1$ is a unitary gauge transformation and $g_2$ is hermitian,
 where there exists a hermitian matrix
 valued function $e \in \mathcal U_{\C}^{\amalg}$ on $T^2$
 with $exp (e) = g_2$. 
For any connection $A'$ in $sl_2 \C$ with 
 $\pll A' - A_1 \pll_{C^r(T^2)} < C\ep$,
 the curvature of $(g^{-1}_2)^*A'$
 is given, modulo the terms containing the quadratics 
 and higher powers of $e$, 
 by 
$$\Delta_{A'} *e + F_{A'},$$
 where $\Delta_{A'}$ is the Laplace
 operator associated to $d_{A'}$.

\begin{claim}
$\pll F_{g^*A} \pll_{C^r(T^2)} < C\delta$,
 $\pll g^* A - g_1^* A \pll_{C^{r+1}(T^2)} < C \delta$.
\end{claim}
\it{Proof of the claim}.
\rm
We know that $\pll F_A \pll_{C^r(T^2)} < \delta$
 (and so $\pll F_{g_1^*A} \pll_{C^r(T^2)} < \delta$) and
 it is transformed to a flat connection by a complex gauge
 transformation.
Moreover, that transformation can be decomposed into two steps,
 first by an element of corollary 4.14 we transform
 $A$ into $sl_2\C$ and then by a constant matrix to a flat connection.
We observe that the directions of
 the first order variations of the curvature 
 during these two processes are perpendicular in $L^2$ sense
 as follows.

Let us consider the first order variation of the curvature when we 
 apply a gauge transformation $g = exp(e)$, $e \in \mathcal U_{\C}^{\amalg}$
 is hermitian matrix valued, to a connection $A' \in sl_2\C$.
It is, as noted above, given by
$$
\Delta_{A'} *e.
$$
In particular, the integration on the torus of each component of it is  
 $0$ by definition of $\mathcal U_{\C}^{\amalg}$.

On the other hand the transformation from 
 $A'$ to a connection in $\T$ is done by a constant gauge transformation
 (it is just the diagonalizaton).
We assumed $A'$ is complex gauge equivalent to $A_1 \in \T$.
The curvature of the connection produced by complex gauge transforming
 $A_1$ by a constant complex gauge transformation
 is represented by a constant component matrix.
In particular, it is orthogonal to 
 $\Delta_{A'} *e$ in $L^2$ sense on the torus.

So, from $\pll F_A \pll_{C^r(T^2)} < \delta$ 
 we can conclude
 that both of $L^2$-norms of
 these (first variations of) curvatures are bounded by $\delta$.
In particular, $\pll F_{g^* A} \pll_{L^2(T^2)} < \delta$.
Moreover, since one part is a constant
 matrix, the $C^r$-norms of the both of these 
 are bounded by $\delta$.

Since the eigenvalues of the operator $\Delta_{A'}$, restricted
 to the direction $\mathcal{U}_{\C}^{\amalg}$ are bounded
 from below, we have $C^{r+2}$ bound of  $e$ by $\delta$.
From this, we have
 $\pll g^* A - g_1^* A \pll_{C^{r+1}(T^2)} < C\delta$.\qed\\


Now our proposition is reduced to the following lemma.
\begin{lem}
There is a constant $\delta_0$ with the following property. 
Suppose $A$ is a connection belonging to $sl_2 \C$ which is complex gauge
 equivalent to $A_1$ and satisfies $\pll F_A \pll_{L^2(T^2)} < \delta < \delta_0$.
Then, there is a unitary gauge transformation $g$ such that
$$
\pll g^* A - A_1 \pll_{C^{r+1}(T^2)} < \sqrt{\delta}
$$
holds. 
\end{lem}
\proof
It suffices to prove that if the constant valued complex gauge transformation
 by a hermitian matrix
 $P = \left( \begin{array}{cc}
          p & q \\
       \overline{q} & r \end{array} \right)$ of determinant one, 
 transforms the connection $A_1$ into a connection with 
 $\pll F_{P^* A_1} \pll_{C^r(T^2)} < \delta$, then the difference from
 $A_1$ to $P^* A_1$ is bounded by $\sqrt{\delta}$. 
The components of $P^* A_1$ are given by
$$
P^* A_1 = \left( \begin{array}{cc}
       a(pr + |q|^2) & 2aqr \\
         -2a\overline{q}p & -a(pr + |q|^2) \end{array} \right)d\overline{z} -
    \left( \begin{array}{cc}
       a(pr + |q|^2) & 2aqr \\
         -2a\overline{q}p & -a(pr + |q|^2) \end{array} \right)^{\dag}dz
$$
Put it as $Bd\overline z - B^{\dag} dz$.
Its curvature is given by 
$$(BB^{\dagger} - B^{\dagger}B) dzd\overline{z}$$
 where
$$
BB^{\dagger} - B^{\dagger}B = a^2 \left(
         \begin{array}{cc}
        -4r^2|q|^2 + 4p^2|q|^2 & -4q(pr + |q|^2)(p + r) \\
         -4\overline{q}(pr + |q|^2)(p + r) &  4r^2|q|^2 - 4p^2|q|^2 
          \end{array} \right).
$$
The assumption is that the absolute values of the
 components of $BB^{\dagger} - B^{\dagger}B$
 are bounded by $\delta$.
Since $pr = 1 + |q|^2$, we have 
$$|p + r| > 2\sqrt{1 + |q|^2}.$$ 
So 
$$\delta > |a^2 4q(pr + |q|^2)(p + r)|
            > |a^2 4q(1 + 2|q|^2)2\sqrt{1 + |q|^2}|.
$$
It follows that 
$$
|a\sqrt q| < \frac{\sqrt{\delta}}{2\sqrt{2}(1 + |q|^2)^{\frac{3}{4}}}.
$$
Then, 
$$
|a(pr + |q|^2 - 1)| = |a(2|q|^2)|
  < \frac{\sqrt{\delta}|q|^{\frac{3}{2}}}{\sqrt{2}(1 + |q|^2)^{\frac{3}{4}}}
  < \sqrt{\delta}.
$$
On the other hand, we have 
$$4a^2|q|^2 |p^2 - r^2| <\delta.$$
If $\sqrt 2 r < p$, then we have
 $r^2 < p^2 - r^2$ and $p^2 < 2(p^2 - r^2)$.
So,
$$
a^2|q|^2r^2,  a^2|q|^2p^2 < \delta.
$$
The case $\sqrt 2 p < r$ is similar.
For the case of $\frac{r}{\sqrt 2} < p < \sqrt 2 r$,
 we have
$$\frac{r^2}{\sqrt 2} < pr = 1 + |q|^2.$$
So, 
$$
\begin{array}{ll}
|aqr| = |a (\sqrt q)^2 r| & <
      \frac{\sqrt{\delta q} |r|}{2\sqrt{2}(1 + |q|^2)^{\frac{3}{4}}}\\
     & < \frac{\sqrt{\delta q}
      (1 + |q|^2)^{\frac{1}{2}}}{2(1 + |q|^2)^{\frac{3}{4}}}\\
                      & < \frac{\sqrt{\delta}}{2}.
\end{array}
$$
Then $|aqp| < |aqr| \sqrt{2} < \sqrt{\delta}$.
The lemma follows. \qed \\

We also have the parametrized version.
\bc
Let $A_s$, $s \in S$ be a smooth family of connections on $E$ 
 which satisfies the assumptions of corollary 4.17.
Let $B_s \in \T$ be the possibly double valued family of
 flat connections to which $A_s$ is
 complex gauge equivalent. 
Then there is a possibly double valued
 smooth family of unitary gauge transformations $h_s$ such that 
 the estimate
$$
\pll h_s^*A_s - B_s \pll_{C^{r+1}(T^2)}
    < C\pll F_{A_s} \pll_{C^r(T^2)}^{\frac{1}{2}}
$$ 
 holds with $C$ a constant not depending on $\{ A_s \}$.
The estimate is done on each fiber and does not contain the estimate
 in the horizontal directions.\qed
\ec
%

\br
Although we dealt with the $SU(2)$ bundles on $T^2$ in this section,
 most of the proofs and results extend in a straightforward manner to
 $SU(N)$ bundles for any $N \geq 2$.
Namely, the reduction to the case of constant connection matrices
 (the first half of the proof of proposition 4.19
 and the results used there) is the same
 (see also remark 5.14). 
\er
So if one can prove the following
 estimate of the components of constant connection matrices 
 by the norm of the curvature, which corresponds to lemma 4.21,
 all of the results in this section
 are generalized to $SU(N)$ gauge groups.
\begin{q}
Let $A_1$ be a real diagonal matrix of trace zero.
Let $P$ be a hermitian matrix of determinant one.
Put $B = PA_1P^{-1}$.
Suppose all of the components of the matrix $BB^{\dag} - B^{\dag}B$ 
 are of order $\delta^2$ ($\delta$ is sufficiently small).
Are the components of B of order $\delta$?
 \end{q}
\begin{rem}
Another direction of generalization is to consider
 the higher genus cases.
In the trivial $SO(3)$ bundle case,
 the argument may get more complicated.

However, in the cases of $SO(3)$ bundles on  
  $\R^2 \times \Sigma_g$ ($g \geq 2$) which restrict to non-trivial bundles on
 fibers are more easy to handle.
Namely, in these cases the moduli spaces of flat connections 
 are smooth and the moduli problems are transversal.
So we can easily deduce the stronger estimate
$$
\pll A - A_0 \pll_{C^r(D \times \Sigma_g)}
 < C \pll F_A \pll_{C^{r-1}(D \times \Sigma_g)},
$$
 here $A_0$ is a family of
 flat connections and $A$ is a family of connections near to $A_0$,
 $D$ is a disc on $\R^2$
 (for norms of sections, see remark 5.2).
Note that this estimate contains the estimates in the horizontal directions.
Using this and the arguments in the following sections, we
 can reprove parts of results (namely, the analysis of bubbles)
 of the paper of Dostoglou-Salamon \cite{DS} and
 the stronger result of the paper of Chen \cite{C}.
\end{rem}

\section{Non-existence of small energy doubly periodic instantons}
Recently, there are large developments in the study of
 finite energy instantons on $T^2 \times \R^2$,
 which are called by the name 
 doubly periodic instantons, by those people including
 Biquard, Jardim, etc., see $\cite{BJ}$, $\cite{J}$.
Although their studies are quite detailed, they assumed
 the decay of the curvature, namely, quadratic decay at the infinity 
 of $\R^2$. 
They have succeeded to classify doubly periodic instantons
 under this assumption.
We deal in this section with the case when this
 decay of the curvature is not assumed.
The result is that even in this
 case the energy quantization holds
 (for the precise statements,
 see theorems 5.5 and 5.6).
This result will be an important step for
 classifying general doubly periodic instantons. 
 
First we show that a small energy doubly periodic instanton
 determines a constant map to the S-equivalence moduli
 of rank 2 bundles on $T^2$. 
Let $\Xi = A + \Phi ds + \Psi dt$ be a doubly periodic instanton
 on the trivial $SU(2)$-bundle $E$ on $T^2 \times \R^2$.
Here we put on $T^2$ and $\R^2$ affine metrics and we put
 on $T^2 \times \R^2$ the product metric.
Denote by $[A(s, t)]$ the S-equivalence class of the holomorphic structure
 of the rank 2 bundle on $T^2$ determined by the fiber (over $(s, t) \in \R^2$)
 direction part $A(s, t)$ of the unitary connection $\Xi$.
Note that the smallness of the energy guarantees that these 
 holomorphic structures are all semi-simple.
Since the ASD connection $\Xi$ determines a holomorphic structure
 on $E$, $[A(s, t)]$ determines a holomorphic map $\rho$ from $\C = \R^2$
 to $S^2$, the S-equivalence moduli of
 rank two semistable bundles with trivial
 determinant over $T^2$:
$$
\rho: (s, t) \in \R^2 \mapsto [A(s, t)].
$$

Here we quote the following result of Fukaya $\cite{F1}$,
 proposition 3.1 and addendum 3.11.
\bt
Let $Y$ be a compact Riemannian 4-manifold with (or without)
 boundary. 
Let $P_G \ra Y$ be a principal $G$-bundle, $G$ is
 a compact Lie group. 
Let $K \subset Y$ be a compact subset such that
 $K \cap \del Y = \phi$. 

Then, there are constants
 $\ep = \ep(K, Y, P_G)$, $C_r = C_r(K, Y, P_G)$
 and $\lam_G = \lam(K, Y, P_G) \in (0, 1]$ 
 which satisfy the following.

Let $A$ be an ASD connection on $Y$ such that
 $\pll F_A \pll_{L^2(Y)} < \ep$. 
Then, there is a flat connection $A_0$ and a gauge transformation
 $g$ which satisfy the following estimate.
$$
\pll g^*A - A_0 \pll_{C^r(K)} < C_r \pll F_A \pll^{\lam_G}_{L^2(Y)}.
$$
Moreover, the flat connection $A_0$ can be taken
 to satisfy
$$
\pll A_0 \pll_{C^r(K)} < C'_r, 
$$ 
 here $C'_r$ is a constant depending on $Y$ and $K$
 (see the remark below for definition of the $C^l$-norms
 of connections and curvatures).\qed
\et
\br
Here we remark about the definitions of $C^l$ norms of connections
 and curvatures.
Generally, they depend on the gauge one takes
 (however, usually these norms are equivalent).
So we cover $Y$ by open balls on which $P_G$ trivialize and
 fix gauges there.
This gauge should have a property such that the connection matrix
 of $A$ has uniformly bounded $C^l$ norm on each open ball.
This can be done, for example, taking usual Uhlenbeck's
 Coulomb gauge.
All the norms of connections and sections 
 in the theorem and the following arguments are defined using this
 fixed gauge.
Norms of connections mean the norm of the matrix valued section 
 given by the difference of the considering connection and the product
 connection in this fixed gauge.
The gauges which will appear in the following are related to these
 fixed gauges by gauge transformations of bounded 
 $C^{\infty}$ norms. 
So actually the norms defined by using occasional (i.e, non-fixed)
 trivialization give the equivalent norms.
\er

\bt
If the Yang-Mills energy of the connection $\Xi$ is sufficiently small,
 the image of the map $\rho$ is a point. 
\et
\proof
We want to compare the Yang-Mills energy of the ASD connection $\Xi$
 and the energy of the holomorphic map defined by $\Xi$.
For this purpose, we first recall the symplectic and complex structures on
 the moduli of flat connections on a trivial bundle on $T^2$. 
The symplectic form $\omega$ at a flat connection $[A]$ 
 (which is not singular, i.e, not
 having isotropy group $SU(2)$) is given by
$$
\omega(v, w) = - 2\int_{T^2} Tr (v \wedge w), 
$$
where $v$ and $w$ are representatives of $H^1_A(T^2, \g_E)$.
The complex structure of the moduli is given by the Hodge $\ast$-operator
 and this and
 the symplectic structure determine the K\"ahler metric
 (on the nonsingular part).

The energy of a holomorphic map from $\C$ to the moduli is 
 given by integrating the pull-back of $\omega$
 (when it is a constant map to the singular points,
 we define the energy to be 0). 
On the other hand, the Yang-Mills energy of $\Xi$ is
 given by 
$$
YM(\Xi) = 2 \int_{\R^2 \times T^2}
         (|\del_sA - d_A \Phi|^2 + |F_A|^2) d \mu.
$$

If the fiber components $A(s, t), (s, t) \in \R^2$
 of $\Xi$ are flat, $[ A(s, t) ]$ gives a holomorphic
 map from $\C$ to $S^2$.
In this case, the form $\del_sA - d_A \Phi$
 is $d_A$ closed and so represents an element of $H^1_A(T^2)$. 
In particular, using the ASD equation, it follows that
 the energy of the map is given by 
$$
 E([A]) = 2 \int_{\R^2 \times T^2}
         |\del_sA - d_A \Phi|^2 d \mu,
$$
and so $YM(\Xi)$ is equal to the energy of the holomorphic map. 

In our case, $A$'s are not necessarily flat, but 
 $YM(\Xi)$ is small. 
Since in any case $\rho$ is holomorphic,
 the image is a point or the whole of $S^2$. 
Assume the image is $S^2$.  
Then, consider a disk $D$ on $Rep_{T^2}(SU(2))$
 remote from the singular four points.

Take a point $x$ on $\R^2$ whose image 
 $\rho(x)$ is contained in $D$.
Let $B_1(x)$ be the unit ball around $x$.
Then since the energy of $\Xi$ is small,
 we can apply the above theorem
 to the restriction of $E$ to $B_1(x) \times T^2$ so that
 there is a flat connection $\Xi_F$ and a unitary gauge transformation $g$
 satisfying 
$$
\pll g^*\Xi - \Xi_F \pll_{C^r(B_1(x) \times T^2)} < C \pll F_{\Xi}
   \pll_{C^{r-1}(\R^2 \times T^2)}^{\lam},
$$ 
 here $\lam$ does not depend on $x$.
We rewrite this $g^*\Xi$ as $\Xi$ to save letters.
Since $B_1(x)$ is simply connected, the gauge equivalence
 class of the fiber part of $\Xi_F$ (so its image in $S^2$)
 is a point.
Because the curvature of $\Xi$ is sufficiently small, we
 can apply corollary 4.22, so that there is a gauge transformation
 $h$ with the fiberwise estimate
$$
\pll h^* (s, t) A(s, t) - A_{0}(s, t) \pll_{C^r(T^2)} < C \pll F_{A(s, t)}
   \pll_{C^r(T^2)}^{\frac{1}{2}}
$$
 $(s, t) \in B_1 (x)$ holds (in this situation double valuedness does not occur
 because $D$ is remote from the singularities).
Here $A_{0}(s, t)$
is the connection in $\T$ which
 is complex gauge equivalent to $A(s, t)$. 
 
Furthermore, exploiting
 the facts that $D$ is remote from the singularities and
 our connection is ASD and
 the smallness of the
 Yang-Mills energy we can deduce stronger estimates
 including estimates for differentials of the horizontal directions.
\begin{prop}
The gauge transformation $h$ above can be taken so that
 the transformed matrix satisfies
$$
\pll h^* A - A_{0}(s, t) \pll_{C^r(B_1(x) \times T^2)} 
    < C_1 \pll F_{\Xi} \pll_{L^2(B_1(x) \times T^2)}.
$$
\end{prop}
\proof
First, we can assume that 
 the $C^r$ norm of $F_{\Xi}$ is bounded by a constant multiple
 of the $L^2$ norm of it.
This is because it is the case for ASD connections
 with usual Uhlenbeck's Coulomb gauge (see remark 5.2).

Recall that,
 on the open disc $D$, the map
$$
D \times V \times W \to \mathcal{A},
$$
 defined by
$$
(a, u, v) \mapsto \exp(v)^* \exp(u)^* a,
$$
 where $V$ is an open neighbourhood of zero in the
 space of sections over $T^2$ of $2 \times 2$ 
 traceless skew-Hermitian matrices
 with the integrals of diagonal entries zero
 and $W$ is an open neighbourhood of zero in the similar
 space of traceless Hermitian matrices,
 is a local diffeomorphism (corollary 4.9, recall that $D$ is remote from
 singular points).

By the proof of proposition 4.19,
 the connection matrices taken in the gauge of that proposition
 can be characterized by the property
 that it is obtained from $A_{0}(s, t)$
 first by a constant Hermitian gauge transformation
 and then by uniquely determined gauge transformation by 
 the exponential of a Hermitian element of $\mathcal{U_{\C}}^{\amalg}$. 

We write this complex gauge transformation $j$ in the form
$$
j = \exp(e_1) \exp(e_2),
$$ 
 here $e_1$ and $e_2$ are sections over $B_1(x)$
 of the bundles whose fibers are above mentioned spaces of
 Hermitian matrices.

Now the fiberwise curvature of $j^* A_{0}(s, t)$ is given by
$$
\ast_{T^2_F} \Delta_{A_{0}(s, t)} (e_1 + e_2),
$$
 modulo higher order terms.
Since the $C^r$ norm of the curvature is bounded by a constant multiple of 
 $\pll F_{\Xi} \pll_{L^2(B_1(x) \times T^2)}$
 and the operator norms of $\Delta_{A_{0}(s, t)}$ are larger than some constant
 $C$ (determined by the distance from $D$ to the singular points),
 and moreover $\pll A_{0}(s, t) \pll_{C^r} \leq const.$, 
 it follows that the $C^r$ norms of $e_i$ are also bounded by
 constant multiple of 
 $\pll F_{\Xi} \pll_{L^2(B_1(x) \times T^2)}$.
From this, the estimate is clear. \qed\\

Now we complete the proof of the theorem.
Viewing now $h$ as a four dimensional gauge transformation over 
 $E \big|_{B_1(x) \times T^2}$, we can get the connection of the form
 $h^* \Xi = h^* A(s, t) + \Phi (s, t) ds + \Psi (s, t) dt$.
Transforming this by the inverse of the
 hermitian transformation $j$ in the proof
 of the proposition, we get
$$
(j^{-1})^* h^* \Xi = A_{0}(s, t) + \Phi_{0}(s, t) ds + \Psi_{0}(s, t) dt.
$$
By the first part of the $HYM$ equation, which is valid
 under complex gauge transformation (see the first part of the next section
 and proposition 6.13), 
 $\Phi_{0}(s, t)$ and $\Psi_{0}(s, t)$ are in $ker d_{A_{0}(s, t)}$.
So the energy of the holomorphic map is given by
$$
2 \int_{B_1(x) \times T^2} | \del_s A_{0}(s, t)|^2 d\mu.
$$

On the other hand, by the proposition,
 we can deduce using Schwarz's inequality the following
 estimate,
$$
\begin{array}{ll}
YM(h^*\Xi |_{B_1(x) \times T^2}) 
 & = 2 \int_{B_1(x) \times T^2}
  (|\del_s h^*A(s, t) - d_{h^*A(s, t)} \Phi (s, t)|^2
                                 + |F_{h^* A(s, t)}|^2) d\mu \\
 & \geq \int_{B_1(x) \times T^2} | \del_s A_{0}(s, t)|^2
               - C \pll F_{\Xi} \pll_{B_1(x) \times T^2}^2,
\end{array}
$$
 with some constant $C$ not depending on
 $\Xi$.

In particular, the energy of the holomorphic map is
 bounded by the constant multiple of the Yang-Mills energy.
However, this Yang-Mills energy is small by assumption.
So the holomorphic map cannot sweep the region around $\rho(x)$.
This contradicts the assumption that 
 the image of $\rho$ is the whole $S^2$. \qed\\

Before we go further, we need to put an assumption.
Recall an ASD connection $\Xi$ induces a holomorphic structure on 
 $E$.
By complex geometry, if the restrictions 
 of $E$ to fibers over some open set of the base
 are isomorphic to the case 2 of
 theorem 4.1, then all fibers except non-accumulating,
 countable number of fibers
 are isomorphic to it.
We need to avoid this case to use the results of section 4.
So we put the following condition.\\

\noindent
($\bold A$) the restrictions of $E$ to fibers over 
 some open subset in $ \R^2$
 are isomorphic to the ones from the case 1 of theorem 4.1.\\

In this case, together with theorem 5.3,
 the restriction of $E$ to any fiber is
 isomorphic to one from the case 1.
For complex geometric aspects of bundles on elliptic fibrations,
 Friedman \cite{Fr}, section 8 is a good reference.
Here we define the open subset $U_{\eta}$ on the moduli, namely,
$$
U_{\eta} = \bigg\{
     \left( \begin{array}{cc}
              a & 0 \\
              0 & -a  \end{array}  \right) d \overline z
           - \left( \begin{array}{cc}
              a & 0 \\
              0 & -a  \end{array}  \right)^{\dag} dz
            | a \in [0, \pi) \times i[0, \pi),  |a - \alpha| < \eta \bigg\}, 
$$
 here $\alpha = n\pi + im\pi$, $n, m = 0, \frac{1}{2}, 1$
 (more precisely, we think $Rep_{T^2}(SU(2))$
 as
 $ \bigg\{ \left( \begin{array}{cc}
              a & 0 \\
              0 & -a  \end{array}  \right) d \overline z
           - \left( \begin{array}{cc}
              a & 0 \\
              0 & -a  \end{array}  \right)^{\dag} dz
   \bigg|a \in [0, \pi) \times i [0, \pi)
   \bigg\} \bigg/ \pm 1$,
   and $U_{\eta}$
 is the image of the quotient by $\pm 1$.
We treat the elements of
 $\bigg\{ \left( \begin{array}{cc}
              a & 0 \\
              0 & -a  \end{array}  \right) d \overline z
           - \left( \begin{array}{cc}
              a & 0 \\
              0 & -a  \end{array}  \right)^{\dag} dz \bigg\}$
 as if it were in the quotient space if no confusion can happen).

Now we begin the proof of the main theorem of this section, 
 the non-existence of small energy 
 doubly periodic instantons. 
Namely we will prove the following theorems.
\bt
Let $\Xi$ be a small energy doubly periodic instanton and let 
 $ \left( \begin{array}{cc}
              a & 0 \\
              0 & -a  \end{array}  \right) d \overline{z} -
       \left( \begin{array}{cc}
              a & 0 \\
              0 & -a  \end{array}  \right)^{\dag} dz
      $, $a \in [0, \pi) \times i[0, \pi)$
 be the element in $\T$ to which $\rho$ of theorem 5.3 maps to,
 (so the condition
 $\left( \bold A \right)$ is automatic).
Suppose $ \left( \begin{array}{cc}
              a & 0 \\
              0 & -a  \end{array}  \right) d \overline{z} -
       \left( \begin{array}{cc}
              a & 0 \\
              0 & -a  \end{array}  \right)^{\dag} dz
      $ is contained in 
 $Rep_{T^2}(SU(2)) \setminus U_{\eta}$.
Then if the energy of $\Xi$ is sufficiently smaller than $\eta$, 
 then $\Xi$ is a flat connection.
\et
\bt
In the case that the image of the map $\rho$ is one of
 the singular points of $Rep_{T^2}(SU(2))$
 and $\Xi$ satisfies the condition $( \bold A )$,
 if the energy of $\Xi$ is smaller than a positive constant
 $\delta_1$, then $\Xi$ is a flat connection.
\et
\noindent
$\it {Proof \; of \; theorem \; 5.5}.$
In the following, as in the last section, the constants
 in various estimates are simply written by $C$.
Their critical dependences are noted whenever necessary.

Let $R$ be a positive integer and $N_R$ be the annulus
 $B_{R+1} - B_R$ on the base.
Here $B_R$ means the ball of radius $R$ with the center at the origin. 
We first prove the following lemma.
\bl
Let $\Xi = A + \Phi ds + \Psi dt$ as above 
 represented in some gauge of $E$ 
 given on 
 $\R^2 \times T^2$ (remember that $E$ is a topologically
 trivial bundle). 
For any positive number $\delta$, there is 
 a positive constant $\ep_0$ independent of $\Xi$ and $R$, depending on $\delta$
 such that
 if $\Xi$ satisfies $\pll F_{\Xi} \pll_{L^2(\R^2 \times T^2)} < \ep_0$, there 
 exists a unitary gauge transformation $g_N$ on $N_R$
 with $\pll g_N \pll_{C^r(N_R \times T^2)} < C$ such that if 
 $g_N^*\Xi = A' + \Phi'ds + \Psi' dt$, then
$$ \pll A'(s, t) - A_0 \pll_{C^r(T^2)} <\delta, \; (s, t) \in N_R $$ holds.
Here
 $A_0 = \left( \begin{array}{cc}
              a & 0 \\
              0 & -a  \end{array}  \right) d \overline{z} -
       \left( \begin{array}{cc}
              a & 0 \\
              0 & -a  \end{array}  \right)^{\dag} dz \in \T$,
 the connection to which $\{A(s, t) \}$ are
 complex gauge equivalent.
\el
\proof
Since we will use this lemma also in the proof of theorem 5.6,
 we do not assume $a \neq 0$ here.
 
Hereafter, we denote the fiber part of the gauge transformed
 connection by 
 some gauge transformation $j$ of $\Xi = A + \Phi ds + \Psi dt$
 by $j^* A(s, t)$
 (this is an abbreviated notation, but actually the fiber
 part is given by $j^* A$ when $j$ is seen as a 
 family of gauge transformations on fibers). 

By the energy smallness, the $C^r$-norm of the connection matrix
 of $\Xi$ is bounded by a constant $C$, not depending on $\Xi$,
 as mentioned in remark 5.2.

Let $U_1$ be the region defined by 
 $$U_1 = \left\{ z \in N_R | 0 < arg z < \frac{3\pi}{R} \right\}.$$
Applying theorem 5.1 to $U_1 \times T^2$,
 we conclude that there exist a unitary gauge transformation $g_1$
 on $E|_{U_1 \times T^2}$
 and a flat connection $\Xi_F$ such that the inequality
$$
\pll g_1^*\Xi - \Xi_F \pll_{C^r(U_1 \times T^2)}
   < C \pll F_{\Xi} \pll_{L^2(U_1' \times T^2)}^{\lam}
$$
 holds, with
$$
\pll \Xi_F \pll_{C^r(U_1 \times T^2)} < C.
$$ 
Here for any open subset $U$, $U'$ is a larger open subset defined by
$$
U' = \{x \in  \R^2 | dist(x, U) <1 \}.
$$

Since $g_1$ satisfies 
$$
g_1^{-1} d g_1 = -g_1^{-1} \Xi g_1 + \Xi_F + \alpha
$$
 with the $C^r$-norm of $\alpha$ small and the $C^r$-norms of 
 $\Xi$ and $\Xi_F$ bounded, we can estimate the $C^r$ norm of 
 $g_1$,
$$
\pll g_1 \pll_{C^r(U_1 \times T^2)} < C
$$
 with some constant $C$ depending only on $r$.

Denoting the fiber component of $\Xi_F$ by $A_F(s, t)$, $(s, t) \in U_1$,
 we have
$$
\pll g_1^*A(s, t) - A_F(s, t) \pll_{C^r(U_1 \times T^2)}
 < C \pll F_{\Xi} \pll_{L^2(U_1' \times T^2)}^{\lam}.
$$

Since $\Xi_F$ is flat and $U_1$ is simply connected,
 $A_F(s, t)$ is unitary gauge equivalent to an 
 element of $\T$ independent of $(s, t)$. 
In particular, there is a unitary gauge transformation 
 on the bundle on $T^2$, parametrized by $(s, t) \in U_1$, $g'(s, t)$, 
 such that $g'(s, t)^*A_F(s, t) \in \T$ holds. 
We can take $g'$ so that the $C^r$-norm of it is also bounded.
This can be shown as follows.
Since $U_1$ is simply connected, there is a gauge in which the flat connection
 $\Xi_F$ is represented by $A_{00} + 0 ds + 0dt$, $A_{00} \in \T$.
We take $g'$ so that it transforms $\Xi_F$ to
 $A_{00} + 0 ds + 0dt$.
Then we have
$$
g'^{-1} d g' = g'^{-1} \Xi_F g' + A_{00},
$$
 and the estimate follows as before:
$$
\pll g' \pll_{C^r(U_1 \times T^2)} < C.
$$
Then we have
$$
\pll g'^*g_1^*A(s, t) - g'^*A_F(s, t) \pll_{C^r(U_1 \times T^2)} 
   < C\pll F_{\Xi} \pll_{L^2(U_1' \times T^2)}^{\lam}.
$$
On the other hand, since $g'^* g_1^* A(s, t)$ is complex gauge
 equivalent to $A_0$, we can put $g'^*A_F(s, t) = A_0$ by
 proposition 4.16. 
Denote $g_1 g' = g$.
As we have seen, the $C^r$-norm of $g$ is bounded.

Doing the same procedure on the overlapping region $U_2 \times T^2$,
 where 
$$
U_2 = \left\{ z \in N | \frac{2\pi}{R} < arg z < \frac{5\pi}{R} \right\},
$$
we have a gauge transformation $h$ such that
$$
\pll h^* A(s, t) - A_0 \pll_{C^r(U_2 \times T^2)}
 < C\pll F_{\Xi} \pll_{L^2(U'_2 \times T^2)}^{\lam}
$$
and 
$$
\pll h \pll_{C^r(U_2 \times T^2)} < C
$$
 hold.

The difference of the gauges $g(s, t)^{-1}h(s, t)$ satisfies the inequality
$$
\pll (g(s, t)^{-1}h(s, t))^*A_0 - A_0 \pll_{C^r(U_1 \cap U_2 \times T^2)}
 < C\pll F_{\Xi} \pll_{L^2(U'_1 \cup U'_2 \times T^2)}^{\lam},
$$
 which follows from the triangle inequality.

By corollary 4.14, there exists
 $v(s, t) \in \mathcal U_{\C}^{\amalg} \times (U_1 \cap U_2)$ such that
$$\pll v(s, t) \pll_{C^r((U_1 \cap U_2) \times T^2)}
 < C\pll F_{\Xi} \pll_{L^2(U'_1 \times T^2)}^{\lam}, \,\,
 \exp(v(s, t))^*(g(s, t)^{-1}h(s, t))^*A_0 \in sl_2\C$$  holds with
 $$\pll \exp(v(s, t))^*(g(s, t)^{-1}h(s, t))^*A_0
 - A_0 \pll_{C^r((U_1 \cap U_2) \times T^2)}
                < C\pll F_{\Xi} \pll_{L^2(U'_1 \times T^2)}^{\lam}.$$ 
This means $\exp(v(s, t))^* (g(s, t)^{-1}h(s, t))^*$ is
 a constant gauge transformation 
 (constant in the direction of the fibers,
 it can depend on the base coordinates and this dependence is
 $C^r$-bounded)
 and so $g(s, t)^{-1}h(s, t)$ is written as a composition of a constant
 gauge transformation and $exp(-v(s, t))$. 

Since the gauge group $SU(2)$ is connected and $exp(-v(s, t))$
 is sufficiently $C^r$-small, we can take a gauge transformation $\widetilde g$
 on $U_1 \times T^2$ such that
$$
\widetilde g = \left\{ \begin{array}{ll}
           g(s, t)^{-1}h(s, t) & on \; (U_1 \cap U_2) \times T^2 \\       
           id      & on \; \{ z \in N | 0
                        < arg z < \frac{\pi}{R} \}
             \end{array} \right.      
$$
in such a way that the norm
 $sup_{(s, t) \in U_1} \pll \widetilde g^* A_0(s, t) - A_0 \pll_{C^r(T^2)}$
 (here $\widetilde g^* A_0$ is considered as a family of connections on $T^2$,
 and the base component is neglected)
 does not exceed
 $sup_{(s, t) \in U_1 \cap U_2}
      \pll (g^{-1} h)^*A_0(s, t) - A_0 \pll_{C^r(T^2)}$.
This can be seen as follows.
Let $\rho$ be a cut off function on $U_1 \times T^2$ which
 depends only on the base and satisfies
$$
\rho = \left\{ \begin{array}{ll}
         1 \; on \; (U_1 \cap U_2) \times T^2 \\
         0 \; on \; \{ z \in N_R | 0 < arg z < \frac{\pi}{R} \}
             \end{array} \right. 
$$
 and 
$$
0 \geq \rho \geq 1, \; \pll \rho \pll_{C^r} \leq 3.
$$
We can extend $v(s, t)$ to $U_1$ so that it satisfies
$$
\pll v(s, t) \pll_{C^r(U_1 \times T^2)} < C \pll F_{\Xi}
 \pll_{L^2(U'_1 \times T^2)}^{\lam}.
$$
As noted above, the transformation $(g^{-1} h)^* = h^* (g^{-1})^*$
 is written as $\exp(-v(s, t))^* exp(a)^*$, where $a$ is 
 constant in the direction of the fibers. 
If the isotropy group of $A_0$ is $SU(2)$, then constant gauge transformations
 do not transform $A_0$.
So $\widetilde g = \exp(\rho a) \exp(-\rho v(s, t))$ suffices our purpose.

Let $\mathfrak{su}(2)$ be the space of flat connections 
 represented by constant matrices in the given gauge.
The adjoint action 
$$
Ad : SU(2) \times \mathfrak{t} \to
 \mathfrak{su}(2), \; Ad(A)(A_0) = A^{-1}A_0A,
$$
 for $A_0 \in  \mathfrak{t}$, whose isotropy is not
 $SU(2)$ is, when we restrict $SU(2)$ to $\exp \mathfrak{h}^{\perp}$
 (here  $\mathfrak{h}^{\perp}$ is the set of elements of $su(2)$
 with zero diagonal components),
 locally a diffeomorphism.
In particular, if we take $exp(\rho a)^* A_0$, the norm 
 $\pll exp(\rho a)^* A_0 - A_0 \pll_{C^r(T^2)}$ on each fiber does not exceed
 that of $\pll exp(a)^* A_0 - A_0 \pll_{C^r(T^2)}$.
So in this case also the definition
 $\widetilde g = \exp(\rho a) \exp(-\rho v(s, t))$ suffices.

Now we have constructed a trivialization on 
 $(U_1 \cup U_2) \times T^2$, which is 
 obtained from the given gauge by transforming by $\widetilde{g}^* g^*$
 on $U_1$ and $h^*$ on $U_2$, so that in this trivialization, 
 the fiber part of the connection matrix, which we rewrite
 by $\Xi = A + \Phi ds + \Psi dt$, 
 satisfies
$$
\pll A(s, t) - A_0 \pll_{C^r(T^2)}
  < C\pll F_{\Xi} \pll_{L^2(N'_R \times T^2)}^{\lam}.
$$
Continuing this process until we go around $N_R$,
 we construct a trivialization on $N_R \times T^2$
 so that the required inequality holds. 
In the final step where we treat the region
$$
U = \left\{ z \in N_R \bigg| \frac{(2R - 2)\pi}{R}
   < arg z < \frac{(2R + 1)\pi}{R} \right\},
$$
 the gauges of $U_1$ and $U$ may not coincide on the overlap.
However, by the same argument as above, we can construct
 a global gauge on $N_R \times T^2$.
\qed\\

Now we go back to the proof of the theorem. 
Apply the previous theorem 5.3 to $\Xi$ and 
let $[A_0]$, $A_0 \in \T \setminus U_{\eta}$
 be the class of the connections representing the
 S-equivalence class to which $A(s, t)$ are mapped. 
By the above lemma 5.7, we can find a gauge on $N_R \times T^2$
 such that the fiber component $A$
 of $\Xi$ satisfies the estimate
$$
\pll A(s, t) - A_0 \pll_{C^r(T^2)} < \delta_0.
$$
 for a small constant $\delta_0$,
 which depends only on $\pll F_{\Xi} \pll_{L^2(\R^2 \times T^2)}$.
Moreover, taking $\pll F_{\Xi} \pll_{L^2(\R^2 \times T^2)}$ sufficiently small,
 we can assume
 that $\delta_0$ can be taken so small
 that it is possible to apply corollary 4.22.

Here we again use proposition 5.4 to obtain the estimate
 in the horizontal direction.
We rewrite it in a form suited to the paragraph here.
\begin{prop}
There is a unitary gauge transformation $a$
 on $N_R \times T^2$ by which the connection matrix
 in the gauge of lemma 5.7 is transformed to a matrix satisfying
$$
\pll a^*A(s, t) - A_0 \pll_{C^r(U_i \times T^2)} 
    < C \pll F_{\Xi} \pll^{\frac{1}{2}}_{L^2(U'_i \times T^2)},
$$
 for every $i$,
 where $a^* A(s, t)$ denotes the fiber component of the gauge
 transformed connection over the fiber on $(s, t) \in \R^2$.
 \qed
\end{prop}
Here we wrote in the right hand side
 $\pll F_{\Xi} \pll^{\frac{1}{2}}_{L^2(U'_i \times T^2)}$
 instead of 
 $\frac{\pll F_{\Xi} \pll_{L^2(U'_i \times T^2)}}{\eta}$,
 noting the assumption that $\pll F_{\Xi} \pll_{L^2(U'_i \times T^2)}$
 is much smaller than $\eta$.
For this point, see the remark 5.14.

Hereafter we write
 $a^* \Xi$ by $\Xi = A + \Phi ds + \Psi dt$ for notational simplicity.
When we complex gauge transform $\Xi$ in a form that the 
 fiber components are in $\T$, the rest of the components
 are also in diagonal constant forms 
 (which is seen from the part of ASD equation which is valid under
 complex gauge transformation, see the first part of the next section
 and lemma 6.13). 
Noting this, the following is clear from the proof of proposition 5.4.
\bc
There are complex traceless diagonal $2 \times 2$ matrix valued
 sections $\Phi_{0}(s, t)$ and $\Psi_{0}(s, t)$ over
 $N_R \times T^2$ such that 
 they are constant in the direction of fibers and the following estimates
$$
\pll \Phi(s, t) - \Phi_{0}(s, t) \pll_{C^r(U_i \times T^2)}
   < C \pll F_{\Xi} \pll^{\frac{1}{2}}_{L^2(U'_i \times T^2)},
$$
$$
\pll \Psi(s, t) - \Psi_{0}(s, t) \pll_{C^r(U_i \times T^2)}
   < C \pll F_{\Xi} \pll^{\frac{1}{2}}_{L^2(U'_i \times T^2)},
$$
 hold. \qed
\ec

The point in these estimates and the following
 construction is that the $C^r$-norms are estimated
 by local quantity ($L^2$-norm of the curvature
 over $U'_i \times T^2$) and not
 by $\delta_0$ appeared before or by
 $\pll F_{\Xi} \pll_{L^2(N_R \times T^2)}$.
This allows us to deduce estimates which do not depend on $R$. 

Note again that in this case the double valuedness of the gauge transformation does
 not occur, since the fiber part of the connections are complex gauge equivalent to
 the unique one in $\T$ (see remark 4.18).
The constant $C$ depends only on $\delta_0$. 
Next we have to estimate the components $\Phi$ and $\Psi$
 under some suitable gauge.\\

Let $U_1$ be the region on $\R^2$ as above.

\begin{prop}
There is a gauge transformation $i$ on $U_1 \times T^2$,
 depending only on the base and
 commutes with $A_0$ when restricted to gauge transformations on bundles
 on fibers, satisfying 
$$
\pll i^*\Xi - A_0 \pll_{C^r(U_1 \times T^2)}
 < C \pll F_{\Xi} \pll_{L^2(U'_1 \times T^2)}^{\frac{1}{2}}. 
$$
Here $A_0 = A_0 + 0 ds + 0dt$ is a connection in four dimension.
\end{prop}
\proof

By proposition 5.8 and corollary 5.9,
 we see that there are diagonal matrix valued functions $\Phi_0$, $\Psi_0$,
 constant along fibers, such that
$$
\Xi = A_0 + \Phi_0 ds + \Psi_0 dt + \kappa, 
$$
 where
$
\pll \kappa \pll_{C^r(U_1 \times T^2)} < C
           \pll F_{\Xi} \pll_{L^2(U'_1 \times T^2)}^{\frac{1}{2}}.
$
Write $ A_0 + \Phi_0 ds + \Psi_0 dt$ as $\Xi_0$.

Then the curvature of $\Xi$ is given by
$$
F_{\Xi} = F_{\Xi_0} + d_{\Xi_0}\kappa + \kappa \wedge \kappa.
$$
From this, we have
$$
\pll F_{\Xi_0} \pll_{C^r(U_1 \times T^2)} < C
           \pll F_{\Xi} \pll_{L^2(U'_1 \times T^2)}^{\frac{1}{2}}.
$$

%
We can see $\Phi_0 ds + \Psi_0 dt$
 as a connection matrix on the rank two bundle on $U_1$.
Then the $C^r$ norm of the curvature of it is bounded by
 $ C \pll F_{\Xi} \pll_{L^2(U'_1 \times T^2)}^{\frac{1}{2}}$.
So by taking two dimensional Coulomb gauge (see \cite{DK}, proposition 4.4.11,
 or in this case as a direct consequence of harmonic
 integral theory),
 we can gauge transform by a diagonal matrix valued function $i$
 (so commutes with $A_0$) to a connection which satisfies
$$
\pll i^* \Phi_0 \pll_{C^r(U_1 \times T^2)} <
  C \pll F_{\Xi} \pll_{L^2(U'_1 \times T^2)}^{\frac{1}{2}}, \:
\pll i^* \Psi_0 \pll_{C^r(U_1 \times T^2)} <
  C \pll F_{\Xi} \pll_{L^2(U'_1 \times T^2)}^{\frac{1}{2}}.
$$
Then it is clear that
$$
\pll i^*\Xi - A_0 \pll_{C^r(U_1 \times T^2)}
 < C \pll F_{\Xi} \pll_{L^2(U'_1 \times T^2)}^{\frac{1}{2}}.
$$
 also holds. \qed\\

Now take another open region $U_2$ as before and
 do the same construction, so we get 
$$
 \pll j^* \Xi  
    - A_0 \pll_{C^r(U_2 \times T^2)} <
   C \pll F_{\Xi} \pll_{L^2(U'_2 \times T^2)}^{\frac{1}{2}}.
$$
 for a gauge transformation $j$.
The difference between the two gauges $h = j^{-1} i$
 depends only on the base and satisfies
$$ \begin{array}{ll}
h^*( j^* \Xi)
  & = h^{-1}d_Uh + h^{-1} j^* \Xi h \\
  & = i^* \Xi , \end{array}  
$$ 
and
$$
h^{-1} A_0 h = A_0.
$$
Here $d_U$ denotes the exterior differential of horizontal direction.
From this, we see
$$
\pll h^{-1} d_U h \pll_{C^r(U_1 \cap U_2 \times T^2)} 
 < C (\pll F_{\Xi} \pll_{L^2(U'_1 \times T^2)}^{\frac{1}{2}}
       + \pll F_{\Xi} \pll_{L^2(U'_2 \times T^2)}^{\frac{1}{2}})
$$
We can extend $h$ to $U_2 \times T^2$
 smoothly depending only on the base, which we denote by $\widetilde{h}$,
 so that it preserves $A_0$ and satisfies,
$$
\widetilde h = \left\{ \begin{array}{ll}
           h & on \; U_1 \cap U_2 \\
           \widetilde h_0 & on \; \frac{4\pi}{R}
           < arg z < \frac{5\pi}{R} \end{array}, \right.
$$
 here $\widetilde h_0$ is a constant matrix
 (constant both in the fiber and in the base directions) 
 which commutes with $A_0$ and
$$
\pll \widetilde h - \widetilde h_0 \pll_{C^r(U_2 \times T^2)}
 < C (\pll F_{\Xi} \pll_{L^2(U'_1 \times T^2)}^{\frac{1}{2}}
       + \pll F_{\Xi} \pll_{L^2(U'_2 \times T^2)}^{\frac{1}{2}}).
$$
Then we get a global gauge transformation $\overline{g}$ on 
 $U_1 \cup U_2$ defined by
$$
\overline{g} = \left\{ \begin{array}{ll}
           i & on \; U_1 \\
           j \widetilde{h} & on \; U_2 \end{array} \right.
$$ 
 which satisfies
$$
\pll \overline{g}^* \Xi
   - A_0 \pll_{C^r(U_1 \times T^2)}
   < C \pll F_{\Xi} \pll_{L^2(U'_1 \times T^2)}^{\frac{1}{2}},
$$
$$
\begin{array}{l}
 \pll \overline{g}^*
 \Xi - A_0 \pll_{C^r((U_2 \cap \{ z | \frac{3\pi}{R} 
 < arg z < \frac{4\pi}{R} \}) \times T^2)} 
 <
   C (\pll F_{\Xi} \pll_{L^2(U'_1 \times T^2)}^{\frac{1}{2}}
       + \pll F_{\Xi}
            \pll_{L^2(U'_2 \times T^2)}^{\frac{1}{2}}),
\end{array}
$$
and
$$
 \pll \overline{g}^*
 \Xi - A_0
 \pll_{C^r((U_2 \cap \{ z | \frac{4\pi}{R} < arg z < \frac{5\pi}{R} \})
     \times T^2)}
     < C \pll F_{\Xi} \pll^{\frac{1}{2}}_{L^2 (U'_2 \times T^2)}.
$$
Continuing this process until we go around $N_R$, 
 we may have, at the last step, 
 a discrepancy of the gauge on the overlap
$$
U_R \cap U_1 = \left\{ z \in N_R \bigg| \frac{(2R - 2)\pi}{R} < arg z
    < \frac{(2R+1)\pi}{R} \right\}
 \cap \left\{ z \in N_R \bigg| 0 < arg z < \frac{3\pi}{R} \right\}.
$$ 
If we denote this discrepancy of the gauge by $h_N$,
 by the above argument, $h_N$ satisfies 
$$
h_N = \exp{(2\pi Q)} h_{N, v},
$$
and
$$
h_N^{-1} A_0 h_N = A_0
$$
 where  $Q$ is a constant matrix which commutes
 with $A_0$ and the gauge transformation $h_{N, v}$ 
 satisfies
$$
\pll h_{N, v} \pll_{C^r((U_R \cap U_1) \times T^2)}
 < C (\pll F_{\Xi} \pll_{L^2(U'_R \times T^2)}^{\frac{1}{2}}
         + \pll F_{\Xi} \pll_{L^2(U'_1 \times T^2)}^{\frac{1}{2}}).   
$$
Then we apply a gauge transformation on the strip 
$
[0, 1] \times [0, \frac{(2R+1)\pi}{R}] \times T^2,
$
 defined by
$$
a = \exp(f(\theta)Q),
$$
 where $f(\theta)$ is a smooth function on $[0, \frac{(2R+1)\pi}{R}]$
 such that
 $f(0) = 0$, $f(\theta) = 2\pi$ on
 $2\pi \leq \theta \leq \frac{(2R+1)\pi}{R}$.
Here $(r, \theta) \in [0, 1] \times [0, \frac{(2R+1)\pi}{R}]$
 is the polar coordinates on (the cover of) $N_R$.

After applying it, the discrepancy of the gauge on the
 overlap is given by $h_{N, v}$, and modifying the gauge on
 $[0, \frac{3\pi}{R}]$ by a gauge transformation $\Check{h}$, 
 which satisfies
$$
\Check{h} = \left\{ \begin{array}{ll}
             h_{N, v} & on \; [0, \frac{\pi}{R}] \\
             id. & on \; [\frac{2\pi}{R}, \frac{3\pi}{R}] \end{array} \right.,
$$
 we get a global gauge transformation $g$ satisfying the following.
\begin{prop}
$$
\pll g^*\Xi - (A_0 + f'(\theta)Q d\theta) \pll_{C^r(U \times T^2)}
  < C \pll F_{\Xi} \pll_{L^2(\widetilde U' \times T^2)}^{\frac{1}{2}}
$$
 on any $U$ with
$$
U = \left\{ z \in N_R \bigg| \frac{M\pi}{R} < arg z < \frac{(M+3)\pi}{R} \right\},
$$
 $M \in \R$ and
$$
\widetilde{U} = \left\{ z \in N_R \bigg|
   \frac{(M-3)\pi}{R} < arg z < \frac{(M+6)\pi}{R} \right\}.
$$
$\widetilde U'$ is the larger open set containing $U$ introduced before. \qed
\end{prop}
Note that $A_0 + f'(\theta)Q d\theta$ is a flat connection,
 since $Q$ and $A_0$ commute.


Now let $\chi$ be a cut-off function depending only on the
 radius on $\R^2$ such that
$$ 
 \chi(r) = \left\{ \begin{array}{ll}
                   0 & on \; B_R \\
                   1 & on \; \R^2 - B_{R+1} \end{array}, \right.
$$
with
$$
\pll \chi(r) \pll_{C^r} < 3.
$$
Set
$$
\widetilde \Xi = g^* \Xi + \chi(r) (\Xi_{00} - g^* \Xi),
$$
where $\Xi_{00}$ is the flat connection $A_0 + f'(\theta)Q d\theta$. 
Then we have
$$
F_{\widetilde \Xi}
  = F_{g^*\Xi}|_{B_{R+1}} + d_{g^*\Xi}(\chi(r)(\Xi_{00} - g^*\Xi))
    + \chi^2 (r) (\Xi_{00} - g^*\Xi) \wedge (\Xi_{00} - g^*\Xi),
$$
 which is zero on $(\R^2 \setminus B_{R+1}) \times T^2$,
 because there $\widetilde{\Xi} = \Xi_{00}$.

Let $\rho$ be a function on $\R^2$ which depends only on the
 radius of $\R^2$ such that
$$
\rho(r) = \left\{ \begin{array}{ll}
                   1 & on \; B_{R+1} \\
                   0 & on \; \R^2 - B_{R+2}. \end{array} \right.
$$
Let $\overline \Xi$ be the connection on $E$ defined by
$$
\overline \Xi = \rho \widetilde{\Xi}.
$$
Then, we calculate that on $(B_{R+2} \setminus B_{R+1}) \times T^2$, 
$$
F_{\overline{\Xi}} = d(\rho(A_0 + f'(\theta)Q d\theta))
  + \rho^2 (A_0 + f'(\theta)Q d\theta) \wedge (A_0 + f'(\theta)Q d\theta).
$$
The second term vanishes because 
 $A_0$ is diagonal and $Q$ also commutes with $A_0$.
So,
$$
 F_{\overline \Xi} \wedge F_{\overline \Xi}
  = 0
$$
on $\R^2 - B_{R+1}$ and 
$$
 F_{\overline \Xi} \wedge F_{\overline \Xi}
  = F_{\widetilde \Xi} \wedge F_{\widetilde \Xi}
$$
 on the whole $\R^2$.  
While, $\overline \Xi$ can be extended to a connection 
 over a bundle on $S^2 \times T^2$, and so 
 the integration of $Tr F_{\overline \Xi} \wedge F_{\overline \Xi}$
 is a topological quantity. 
Combining these remarks, we see that
$$
\int_{\R^2 \times T^2} Tr F_{\widetilde \Xi} \wedge F_{\widetilde \Xi}
 = 8\pi^2 n
$$
 for some $n \in \Z$. 
Now for any sufficiently large $R$, we have,
$$ \begin{array}{lll}
\pll F_{g^*\Xi} \pll^2_{L^2(\R^2 \times T^2)} + 8\pi^2 n & \\
  = - \int_{\R^2 \times T^2} (Tr F_{g^*\Xi} \wedge F_{g^*\Xi} - 
          Tr F_{\widetilde \Xi} \wedge F_{\widetilde \Xi}) & & \\
  = - \int_{(\R^2 - B_{R+1}) \times T^2} Tr  F_{g^*\Xi}
               \wedge F_{g^*\Xi} & & \\
  \hs{.3cm} + \int_{(B_{R+1} - B_R) \times T^2}
              Tr(F_{g^*\Xi} \wedge d_{g^*\Xi}\chi(r)(\Xi_{00} - g^*\Xi)
                +   d_{g^*\Xi}(\chi(r)(\Xi_{00} - g^*\Xi))
              \wedge F_{g^*\Xi}) & & \\
  \hs{.3cm} + \int_{(B_{R+1} - B_R) \times T^2}
              Tr(F_{g^* \Xi} \wedge \chi^2 (r)
                (\Xi_{00} - g^*\Xi) \wedge (\Xi_{00} - g^*\Xi) & & \\
    \hs{4cm}   + \chi^2(r)(\Xi_{00} - g^*\Xi) \wedge (\Xi_{00} - g^*\Xi)
                   \wedge F_{g^* \Xi}) & &\\
  \hs{.3cm} + \int_{(B_{R+1} - B_R) \times T^2}
            Tr d_{g^* \Xi}\chi(r)(\Xi_{00} - g^* \Xi)
    \wedge d_{g^* \Xi}\chi(r)(\Xi_{00} - g^*\Xi) & & \\
  \hs{.3cm} + \int_{(B_{R+1} - B_R) \times T^2}
            Tr(d_{g^*\Xi}(\chi(r)(\Xi_{00} - g^*\Xi))
                \wedge ( \chi^2(r)(\Xi_{00} - g^*\Xi) \wedge (\Xi_{00} - g^*\Xi))
               & & \\
   \hs{4cm} + ( \chi^2(r)(\Xi_{00} - g^*\Xi) \wedge (\Xi_{00} - g^*\Xi))
                 \wedge d_{g^*\Xi}(\chi(r)(\Xi_{00} - g^*\Xi))) & &\\
  \hs{.3cm} + \int_{(B_{R+1} - B_R) \times T^2}
            Tr(\chi^4(r)(\Xi_{00} - g^*\Xi) \wedge (\Xi_{00} - g^*\Xi)
                  \wedge (\Xi_{00} - g^*\Xi) \wedge (\Xi_{00} - g^*\Xi)). & &
   \end{array}
$$
 for some $n \in \Z$.

We first estimate the second term. 
Since 
$$
\begin{array}{ll}
|\int_{(B_{R+1} - B_R) \times T^2}
Tr(F_{g^*\Xi} \wedge d_{g^*\Xi}\chi(r)(\Xi_{00} - g^*\Xi)
                +   d_{g^*\Xi}(\chi(r)(\Xi_{00} - g^*\Xi))
              \wedge F_{g^*\Xi})| & \\
 \hs{3cm} < 2 \int_{N_R \times T^2} |F_{g^*\Xi}|
               |d_{g^*\Xi}\chi(r)(\Xi_{00} - \Xi)| d\mu, &
\end{array}
$$
 we are going to estimate the right hand side.
By Schwarz's inequality and proposition 5.11, we have
$$ \begin{array}{ll}
\int_{N_R \times T^2} |F_{g^*\Xi}| |d_{g^*\Xi}\chi(r)(\Xi_{00} - g^*\Xi)| d\mu
 & < C \sum_i \int_{U'_i \times T^2} |F_{g^*\Xi}|
        \pll F_{g^*\Xi} \pll_{L^2(\widetilde{U'_i}
                    \times T^2)}^{\frac{1}{2}} d\mu \\
 & < C \sum_i \pll F_{g^*\Xi} \pll_{L^2(\widetilde{U_i}
            \times T^2)}^{\frac{3}{2}}.
 \end{array}
$$
Now, since $F_{g^* \Xi}$ is square integrable, 
 there is $R$ such that 
$$
- \int_{N'_R \times T^2} TrF_{g^*\Xi} \wedge F_{g^*\Xi}
    < \frac{\pll F_{g^*\Xi} \pll_{L^2(\R^2 \times T^2)}^8}{R}
$$
 holds 
 since, if not,
$$
- 2 \int_{\R^2 \times T^2} TrF_{g^*\Xi} \wedge F_{g^*\Xi} 
     > \pll F_{g^*\Xi} \pll^8_{L^2( \R^2 \times T^2)} (\frac{1}{R} +
      \frac{1}{R+1} + \cdots),
$$
 holds, which is a contradiction because the right side diverges. 
Now we want to estimate 
 $\sum_i \pll F_{g^*\Xi} \pll_{L^2(\widetilde U'_i \times T^2)}^{\frac{3}{2}}$ 
 under the condition 
$$ \begin{array}{ll}
 \sum_i \pll F_{g^*\Xi} \pll_{L^2(\widetilde U'_i \times T^2)}^2
 & < - 4\int_{N' \times T^2} TrF_{g^*\Xi} \wedge F_{g^*\Xi} \\
 & < \frac{4\pll F_{g^*\Xi} \pll_{L^2(\R^2 \times T^2)}^8}{R}.
 \end{array}
$$
We use the following elementary result.
\bl
Let $a_1, \dots, a_n$ be a sequence of positive numbers.
Let $r_1 < r_2$ be positive numbers.
We vary $a_1, \dots, a_n$ under the condition that
 the sum $\sum_{i = 1}^{n} a_i^{r_2} = C$, $C$ is a fixed constant.
Then, $\sum_{i = 1}^{n} a_i^{r_1}$ takes the maximum when
 all the $a_i$ are equal.  
\el
\proof
First we prove the case of $n = 2$.
Then we have $a_1^{r_2} + a_2^{r_2} = C$.
So $a_1 = (C - a_2^{r_2})^{1/{r_2}}$.
Put $f(a_2) = a_2^{r_1} + (C - a_2^{r_2})^{r_1/{r_2}}$.
It is easy to see that $f(a_2)$ takes its maximum (under the condition
 $a_1$, $a_2$ $> 0$) when $a_2^{r_2} = C/2$, namely, $a_1 = a_2$.
The general case follows from this.
Namely, suppose $a_1, \dots, a_n$ are not equal.
Then, take maximum and minimum $a_i$, $a_j$.
We vary these under the condition $a_i^{r_2} + a_j^{r_2}$ is fixed.
Then, by the above argument, $a_i^{r_1} + a_j^{r_1}$ takes
 its maximum when we take $a_i = a_j$, and this new
 sequence $a_1, \dots, a_n$ has larger $\sum_{i = 1}^{n} a_i^{r_1}$
 than the original one.
Iterating this procedure, $a_1, \dots, a_n$ converges to
 the sequence such that $a_1 = \cdots = a_n$.
\qed \\

\nnn
Using this lemma, We have
$$ \begin{array}{ll}
\sum_i \pll F_{g^*\Xi} \pll_{L^2(\widetilde U'_i \times T^2)}^{\frac{3}{2}}
 & < (2R-1) (\frac{4\pll F_{g^*\Xi}
     \pll_{L^2(\R^2 \times T^2)}^8}{R} \cdot \frac{1}{2R-1})^{\frac{3}{4}} \\
 & < \frac{5 \pll F_{g^*\Xi} \pll_{L^2(\R^2 \times T^2)}^6}{R^{\frac{1}{2}}}.
 \end{array}
$$
We have, by the same calculation as above,
$$\begin{array}{lll}
|\int_{N_R \times T^2}
              Tr(F_{g^* \Xi} \wedge \chi^2 (r)
                (\Xi_{00} - g^*\Xi) \wedge (\Xi_{00} - g^*\Xi) & & \hs{2cm} \\
      \hs{.3cm} + \int_{{N_R} \times T^2}
              Tr(\chi^2 (r)
                (\Xi_{00} - g^*\Xi) \wedge (\Xi_{00} - g^*\Xi)
                   \wedge F_{g^* \Xi}| d\mu & & \\
      \hs{.3cm} < \int_{N_R \times T^2}
            2|Tr( \chi^2 (r) F_{g^* \Xi} \wedge
                    (\Xi_{00} - g^*\Xi) \wedge (\Xi_{00} - g^*\Xi))| & & \\
      \hs{.3cm} < \sum_i \pll F_{g^* \Xi} \pll^2_{L^2(\widetilde U'_i \times T^2)} & & \\
      \hs{.3cm} <   \frac{(2R-1)(4 \pll F_{g^*\Xi}
                 \pll_{L^2(\R^2 \times T^2)}^8)}{R} \cdot \frac{1}{2R-1} & & \\
      \hs{.3cm} < \frac{4 \pll F_{g^*\Xi} \pll_{L^2(\R^2 \times T^2)}^8}{R},
\end{array}
$$

$$ \begin{array}{ll}
|\int_{N_R \times T^2}
     Tr d_{g^*\Xi}\chi(r)(\Xi_{00} - g^*\Xi)
      \wedge d_{g^*\Xi}\chi(r)(\Xi_{00} - g^*\Xi)| & \hs{3cm} \\
 \hs{.5cm}
    < \int_{N_R \times T^2} |Tr d_{g^*\Xi}\chi(r)(\Xi_{00} - g^*\Xi)|^2 d\mu &\\
 \hs{.5cm} < \sum_i \pll F_{g^*\Xi} \pll_{L^2(\widetilde U'_i \times T^2)} & \\
 \hs{.5cm} < (2R - 1) (\frac{4\pll F_{g^*\Xi}
       \pll_{L^2(\R^2 \times T^2)}^8}{R} \cdot \frac{1}{2R-1})^{1/2} & \\
 \hs{.5cm} < 4\pll F_{g^*\Xi} \pll_{L^2(\R^2 \times T^2)}^4,
 \end{array}
$$

$$
\begin{array}{lll}
|\int_{{N_R} \times T^2}
            Tr(d_{g^*\Xi}(\chi(r)(\Xi_{00} - g^*\Xi))
                \wedge ( \chi^2(r)(\Xi_{00} - g^*\Xi) \wedge (\Xi_{00} - g^*\Xi))
               & & \\
   \hs{.3cm} + ( \chi^2(r)(\Xi_{00} - g^*\Xi) \wedge (\Xi_{00} - g^*\Xi))
                 \wedge d_{g^*\Xi}(\chi(r)(\Xi_{00} - g^*\Xi)))| & &\\
    \hs{.3cm}  < \int_{N_R \times T^2}
          |Tr (\chi^2(r) d_{g^*\Xi}(\chi(r)(\Xi_{00} - g^*\Xi))
                   \wedge (\Xi_{00} - g^*\Xi) \wedge (\Xi_{00} - g^*\Xi))| d\mu & & \\
    \hs{.3cm} < \sum_i \pll F_{g^*\Xi}
                 \pll^{\frac{3}{2}}_{L^2(\widetilde U'_i \times T^2)} & & \\
    \hs{.3cm} < (2R - 1) (4 \frac{2\pll F_{g^*\Xi}
       \pll_{L^2(\R^2 \times T^2)}^8}{R}
            \cdot \frac{1}{2R-1})^{\frac{3}{4}} & & \\
    \hs{.3cm} < \frac{5 \pll F_{g^*\Xi} \pll_{L^2(\R^2 \times T^2)}^6}{R^{\frac{1}{2}}},
\end{array}
$$

$$ \begin{array}{ll}
|\int_{{N_R} \times T^2}
            Tr(\chi^4(r)(\Xi_{00} - g^*\Xi) \wedge (\Xi_{00} - g^*\Xi)
                  \wedge (\Xi_{00} - g^*\Xi) \wedge (\Xi_{00} - g^*\Xi))| & \\
 \hs{.3cm} < \sum_i \pll F_{g^* \Xi} \pll^2_{L^2(\widetilde U'_i \times T^2)} &  \\
      \hs{.3cm} <   (2R-1)(\frac{4\pll F_{g^*\Xi}
                 \pll_{L^2(\R^2 \times T^2)}^8}{R} \cdot \frac{1}{2R-1}) &  \\
      \hs{.3cm} < \frac{4 \pll F_{g^*\Xi} \pll_{L^2(\R^2 \times T^2)}^8}{R}.
\end{array}
$$

By these estimates, we have
$$ \begin{array}{ll}
- \pll F_{g^*\Xi} \pll_{L^2(\R^2 \times T^2)}^4
   - & \int_{(\R^2 - B_{R+1}) \times T^2} Tr F_{g^*\Xi} \wedge F_{g^*\Xi} \\
 & < \pll F_{g^*\Xi} \pll_{L^2(\R^2 \times T^2)}^2 + 8\pi^2 n \\
 & < \pll F_{g^*\Xi} \pll_{L^2(\R^2 \times T^2)}^4
   - \int_{(\R^2 - B_{R+1}) \times T^2} Tr F_{g^*\Xi} \wedge F_{g^*\Xi}.
\end{array}
$$
From this inequality, we see $n = 0$.
Moreover, for any large number $M$, there is $R > M$ for which 
 this inequality holds.
So it follows that
$$
\pll F_{g^*{\Xi}} \pll^2_{L^2(\R^2 \times T^2)}
  < \pll F_{g^*{\Xi}} \pll^4_{L^2(\R^2 \times T^2)}.
$$
Since $\pll F_{g^*{\Xi}} \pll_{L^2(\R^2 \times T^2)}$ is small,
 $F_{g^*\Xi} = 0$. \qed\\

Then we turn to the case where the image of $\rho$ is
 $0 \in \T$.
The other cases when the images of $\rho$ is one of the singular points of 
 $Rep_{T^2}(SO(3))$ can be treated in the same manner.\\

\noindent
\it{Proof of theorem 5.6.} 
\rm
First we apply lemma 5.7 to construct a gauge
 in which the fiber component of the connection is
 sufficiently close to $\T$.
In this case, we can transform the fiber component of
 the connection to $0 \in \T$ by a uniquely determined
 gauge transformation as in corollary 4.14.
As noted in the proof of proposition 5.4,
 in this gauge the $\Phi$ and $\Psi$ components
 are also diagonal and constant in the
 fiber direction.
On the other hand, the method of the
 proof of proposition 5.4
 is valid in this case also by replacing $D$ by $0 \in \T$,
 $V$ by a neibourhood of zero in the space of traceless 
 skew hermitian $2 \times 2$ matrix valued sections over
 $T^2$ with the integrals of the entries are zero,
 and $W$ by a similar space of traceless hermitian 
 matrix valued sections.
Then propositions 5.8, 5.9, 5.10 are also valid
 with stronger statements replacing 
 $\pll F_{\Xi} \pll_{L^2}^{\frac{1}{2}}$ by
 $\pll F_{\Xi} \pll_{L^2}$.
The rest of the proof goes in the same way. \qed

\br
It may be plausible that in general, any ASD connection
 whose image by $\rho$ is a point
 is necessarily a product of a flat connection on $T^2$ and
 a trivial connection on $\R^2$.
If this is shown, it will give the desirable energy bound by the constant
 $8 \pi^2$.
\er
\br
In the proof, we used the assumption that
 $\pll F_A \pll_{L^2(\R^2 \times T^2)}$
 is much smaller than $\eta$ only when
 we deduced the estimate in the horizontal direction
 (proposition 5.8).
This proposition is described in a weaker form than possible
 (proposition 5.4 is its stronger form).
The reason for this is that if we can prove
 proposition 5.8 without the above assumption,
 the whole proof is valid and shows the energy quantization 
 for small energy ASD connections (there is no comparison with $\eta$).
Unfortunately, we have this estimate only in the fiber direction yet
 (corollary 4.22).

On the other hand, this means that we did not need
 corollary 4.22 for our proof (only used the technique)
 and so we can extend the energy quantization result
 in this section to higher rank cases
 when we assume the image of the map $\rho$
 to the representation space is a point which
 is remote from highly degenerating loci,
 because there the argument used in proposition 5.4
 is still valid. 
\er

Applying similar argument to the case of
 $SO(3)$ bundles $E$ on $\R^2 \times \Sigma_g$ ($g \geq 2$),
 which restrict on each fiber to a nontrivial bundle,
 we can reproduce the stronger energy quantization result
 sketched in Dostoglou-Salamon \cite{DS}, pages 632-633. 
See also its erratum, page 11.
\begin{prop}
Any ASD connection on $E$ over $\R^2 \times \Sigma_g$ has
 energy $8 \pi^2 n$, $n \in \Bbb{N}$.
\end{prop}
\proof
First we remark that in this case it is not clear that
 the map $\rho$
 from the base space to $Rep_{\Sigma_g}(E)$ (which is the 
 moduli of flat connections on a non-trivial $SO(3)$ bundle $E$
 on $\Sigma_g$, and it is known that $Rep_{\Sigma_g}(E)$
 is a smooth, compact K\"ahler manifold of dimension $6g - 6$) 
 to be defined by
$$
x \mapsto [A(s, t)],
$$
 where $[A(s, t)]$ denote the complex gauge 
 equivalence class of the fiber part of the connection
 over $(s, t) \in \R^2$ is defined over whole $\R^2$,
 because here we do not assume the smallness of the energy.
But at least we know that when the curvature of the 
 connection is small, this map is well defined,
 which is seen from lemma 2.14 of \cite{F1}.

Even so, there is a difference from the main text, namely
 the image of the map need not be a point. 
However, we still can show by the method of the proof of 
 proposition 5.4, the energy of this map is bounded by the Yang-Mills
 energy of the ASD connection.
In particular 
 the image of the annulus $N_R$ bounded by circles of radius
 $R$ and $R+1$ will be included in a small open ball $D$ of 
 $Rep_{\Sigma_g}(E)$ if $R$ is sufficiently large.
Then we can assume that there is a slice consisted of flat connections
 in the space $\mathcal{A}$ of
 connections on $E$ which represents the complex gauge equivalence
 classes of $D$ (see again lemma 2.14 of \cite{F1}).
We denote by $\phi$ the map from $N_R$ to $\mathcal A$
 defined by composing $\rho$
 and the lift to this slice.
We can assume the $C^r$ norm of $\phi$ is small.

Moreover, for any positive number $\delta$,
 there is $R$ such that 
$$
\int_{N_R \times \Sigma_g} |F_{\Xi}|^2 d\mu < \frac{\delta}{R}
$$
 holds, because if not, 
$$
\pll F_{\Xi} \pll_{L^2(\R^2 \times \Sigma_g)}^2
               > \sum_{i = M}^{\infty} \frac{\delta}{i},
$$
 holds, which contradicts to the fact that
 $F_{\Xi}$ is square integrable.

Now take open subsets $U_i$ as in the proof of theorem 5.3.
That is,
$$
U_1 = \{ z \in N_R | 0 < arg z < \frac{3\pi}{R} \},
$$
 and so on.
By theorem 5.1, we can unitary gauge transform $\Xi|_{U_i}$
 so that $g^* \Xi$ is close to a flat connection $C_i$:
$$
\pll g^*\Xi - C_i \pll_{C^r(U_i \times \Sigma_g)}
    < C \pll F_{\Xi} \pll_{L^2(N'_R)}^{\lam}.
$$
Here $C$ depends only on the metric on $\Sigma_g$.
We can assume the fiber component $A_i$ of $C_i$
 is close to  $\phi(s, t)$ and in fact replace
 $A_i$ by $\phi$ keeping the estimate valid
 (but $C_i$ may be no longer flat),
 because the $C^r$-norm of the curvature is sufficiently small and
 $\Xi$ is ASD.
Denote $g^* \Xi = A + \Phi ds + \Psi dt$.

By the last remark 4.26 of section 4, we can further gauge transform 
 so that the fiber component satisfies
$$
\pll h^* A - \phi(x)\pll_{C^r(U_i \times \Sigma_g)}
 < C \pll F_{\Xi} \pll_{L^2(N'_R \times \Sigma_g)} < C \frac{\delta}{R},
$$ 
 $x \in U_i$.
Note in remark 4.26, the right hand side was $\pll F_{\Xi} \pll_{C^r}$,
 but here $\Xi$ is ASD and its $C^r$ norm is estimated by
 the $L^2$ norm.

On the other hand, when we complex gauge transform 
 $A(s, t)$ to $\phi(s, t)$, the $\Phi$ and $\Psi$ components are determined
 uniquely from the part of the ASD equations, which is valid under
 complex gauge transformation.
Denote by $\Xi_0$ the connection determined in this way.

The point here is that since there is no isotropy group 
 of the fiber part of the connection, the whole gauge transformation
 of $\Xi$
 is fixed once we specify the gauge transformation of the fiber direction 
 on all fibers.
In particular, we have an estimate
$$
\pll h^* \Xi - \Xi_0 \pll_{C^r(U_i \times \Sigma_g)}
   < C \frac{\delta}{R}.
$$

Moreover, by the same reason,
 we can globalize the gauge from each $U_i$ to the whole $N_R$,
 keeping the inequality
$$
\pll h^* \Xi - \Xi_0 \pll_{C^r(N_R \times \Sigma_g)}
    < C \frac{\delta}{R},
$$  
 valid.

Now the cut-off construction as in the proof of theorem 5.5 applies
 (this time we cut off to the connection which is determined by the 
 method above from the fiber part of connections
 $Im \phi|_{ \{x \in \R^2 \big| |x| > R\} }$),
 with error term of order $C \delta$.
Thus when $R$ is large enough, we can show by the method of the proof of 
 proposition 5.4, that the difference between $YM(\Xi \big|_{|x| > R})$
 and $YM(\Xi' \big|_{|x| > R})$, where $\Xi'$ is the connection
 made by the cutting-off construction is bounded by
 $C \delta + C' YM(\Xi \big|_{|x| > R})$.
By taking $R$ large enough, this number can be bounded by an arbitrary
 small positive number $\ep$.
The cut-off connection extends to a connection on a bundle
 over the compactified space $S^2 \times \Sigma_g$
 and so has energy $8 \pi^2 n$, $n \in \Bbb{N}$.
It follows that if we write the energy of
 the original ASD connection $\Xi$ over $\R^2 \times \Sigma_g$
 as $YM(\Xi)$, we have
$$
8\pi^2 n - \ep < YM(\Xi) < 8 \pi^2 n + \ep. 
$$ 
Since $\ep$ can be taken to be arbitrary small, we have the
 energy quantization result.\qed


\section{Bubbles as holomorphic maps to the representation space}
In this section we treat the third type of the bubbles,
 namely the case when it induces a holomorphic map to the representation 
 space $Rep_{T^2}(SO(3))$ and prove the main analytic result
 theorem 3.1.
\begin{convention}
In the following, by hyperK\"ahler rotation,  
 we think the base and the fibers of $\widehat{M_{\ep}}$
 are holomorphic submanifolds.
And in this section, we write the $SO(3)$-bundles $Ad E_{\ep_{\nu}}$
 simply as $E_{\ep_{\nu}}$.
\end{convention}

Since an ASD connection on $\widehat{M_{\ep}}$ is still ASD on hyperK\"ahler
 rotated manifolds, there is no problem in considering
 the behavior of a sequence
 of ASD connections on  hyperK\"ahler
 rotated manifolds.
In particular, stable holomorphic bundles on $\widehat{M_{\ep}}$ are still 
 holomorphic and stable on  hyperK\"ahler
 rotated manifolds, because they allow
 non-degenerate ASD connections.

In this convention, the K\"ahler form is given by
$$
\omega_{\ep} = d \check{s} \wedge d \check{t}
                 - \ep^2 dx \wedge dy, 
$$
 here $\check{s}$, $\check{t}$ are complex coordinates of the base,
 introduced in chapter 2.
To bring the equation in a more tractable form,
 we change the coordinates on the 
 base from affine coordinates to the isothermal coordinates.

Recall that the metric on the base is locally given by
 a potential $h$:
$$
g_{ab} = \frac{\del^2 h}{\del a \del b},
$$
$a, b \in  \{ s, t \}$
 (note that the underlying Riemannian metric is not changed under
 hyperK\"ahler rotation).
It is a classical result that 
there exist coordinates $S, T$ such that the metric is written
 in this coordinates as
$$
g = f(S, T)(dS^2 + dT^2),
$$
 here $f$ is a positive function.

Then, in coordinates $(S, T, x, y)$, the 
 ASD equations for a connection 
 $\Xi = A + \Phi dS + \Psi dT$ on $\widehat{M_{\ep}}$ are 
$$
(H_1) \;\; (\del_T A - d_A\Psi) + *_{T^2_F}(\del_S A - d_A \Phi) = 0,
$$
$$
(H_2) \;\; \del_T \Phi - \del_S \Psi - [\Phi, \Psi]
                - f(S, T) \ep^{-2} *_{T^2_F} F_A = 0.
$$
 where $*_{T^2_F}$ is the Hodge operator on the non-scaled fibers
 (it can vary as the fiber varies).

Note that the first of these equations continues to hold
 if we apply a complex gauge transformation to the connection,
 because it is the equation saying that the (0,2) and (2,0) parts
 of the curvature vanish, and this property is 
 conserved under complex gauge transformations.

We are considering the following situation.
Let $(E_{\ep_{\nu}} \to \widehat{M}_{\ep_{\nu}}, \Xi_{\ep_{\nu}})$ be a 
 rescaling family of complex two dimensional K\"ahler
 tori with affine structures, 
 family of $SO(3)$ bundles of fixed topological type on them
 coming from stable holomorphic bundles and
 ASD connections on them, as in section 2.
%

Our result can be stated in a form that
 there is a subsequence of these connections
 which converges to a limit connection
 on $\pi^{-1}(T^2_B \setminus S)$, for some countable subset $S$
 with finite accumulation points.
First we have to specify this subset.

By proposition 4.5,  if $A_{\nu}(s, t)$ induce semisimple holomorphic
 structures on the bundles restricted to fibers and
 are not complex gauge equivalent to
$
 \left( \begin{array}{ll}
          0 & b \\
          0 & 0 \end{array} \right)d \overline z - 
 \left( \begin{array}{ll}
          0 & b \\
          0 & 0 \end{array} \right)^{\dag}dz,
$
 there is a complex gauge transformation $g_{\nu}(s, t)$
 parametrized by the base, which satisfies
$$
g_{\nu}^* A_{\nu} \in \mathfrak{t}.
$$
In any case, the S-equivalence class of $g_{\nu}^* A_{\nu}$,
 $[g_{\nu}^* A_{\nu}]$ gives a sequence of sections
 $\phi_{\nu} : T^2_B \setminus Q_{\nu} \to X$,
 where $Q_{\nu}$ is the finite point set where the bundle restricted to the 
 fibers over $x \in Q_{\nu}$ is unstable.
Here $X$ is the $S^2$-bundle over $T_B^2$ whose fiber over
 $x \in T_B^2$ is $Rep_{T_{F, x}^2}(SO(3))$,
$T_{F, x}$ is the fiber over $x$.

The complex structure of $X$ is defined locally in the trivializing open 
 subset of $\widehat{M_1} \to T_B^2$,
 by the product of the complex structure of the base $T_B^2$
 and the complex structure of $Rep_{T_{F}^2}(SO(3))$.
This defines a global complex structure.
The symplectic form of $X$ is defined by the sum of the pullback of
 that of the base and the singular symplectic form
 on $Rep_{T_{F}^2}(SO(3))$, discussed in the proof of theorem 5.3.
These define a (orbi-) K\"ahler structure on $X$.

Here we used the terms `S-equivalence' and `unstable'
 in a somewhat extended way, because $A_{\nu}$
 are connections on $SO(3)$-bundles.
We simply mean by S-equivalence here
 the identification mentioned in remark 4.3
 and unstability of the connection is defined to be
 the property that the corresponding connection on the lifted
 $SU(2)$ bundle on $T^2$ defines an unstable holomorphic structure.

In general, at the level of holomorphic bundles,
 after a sequence of allowable elementary modifications,
 we can assume the restriction of $E_{\ep_{\nu}}$ to any fiber is semistable,
 though this process changes the topological types of the bundles
 (see \cite{Fr}, pages 148-149).
But this also shows that the section $\phi_{\nu}$ can be extended to $T^2_B$
 because this process does not change the complex gauge equivalence classes
 of the bundles on other fibers.
We denote this extended map also by $\phi_{\nu}$.
\begin{rem}
In general, the complex structure 
 (or the conformal class of the metric) of the fiber torus varies,
 and so does the complex structures of the associated representation spaces.
However, in the case of $SO(3)$ bundles, since $Rep_{T^2}(SO(3))$
 is $S^2$ and so there is a unique complex structure on it.
Moreover, the first part of the ASD equation ($H_1$)
 (for complex gauge transformed 
 connection) just means this section is holomorphic with respect to the 
 complex structure of the base associated to the conformal structure,
 which is of course clear from the point of view
 of hyperK\"ahler rotation.
\end{rem}
We can consider the energies of these sections and
 they are bounded by
 $8\pi^2(c_2(E) - \frac{1}{2}c_1^2(E))$
 of the
 $U(2)$ bundle $E$ from which we have started
 (this is $-2\pi^2 p_1(Ad E)$ of our $SO(3)$ bundle $Ad E$)
 and it does not depend on $\nu$.

Together with the energy bound, we have the following observation
 from the theory of Gromov compactness and minimal surfaces
 (\cite{SU}, \cite{Y}).
\begin{prop-def}
Taking a suitable subsequence, $\phi_{\nu}$
 converges to a limit holomorphic section
 $\phi$ on $T_B^2$ except for finite number of points. 
Let $S_1$ be this finite number of points.
Since the energy of $\phi$ is bounded from above, it can be
 extended to a holomorphic section
$
\phi : T_B^2 \to X. 
$ \qed
\end{prop-def}
Let $S_2$ be the subset of $T_B^2$ defined by\\

$S_2$ = $\{ x \in T_B^2 | \phi(x)$ is one of the points
           on $Rep_{T^2}(SO(3))$
 which have isotropy $SO(3) \}$.\\

$S_2$ is a finite subset of $T_B^2$, unless $\phi$ is a constant section
 to the singular points.
Here we put an assumption, which is somewhat similar to the assumption 
 ($\bold A$) of the previous section.\\

\noindent
($\bold B$) The section $\phi$ is not a constant section to
 the singular points (that is, the points corresponding to the 
 connections with isotropy $SO(3)$)
 on $Rep_{T^2}(SO(3))$. \\

Let us take an open disk $V$ on $T^2_B$ on which
 the bundle
 $\widehat{M_{\ep_{\nu}}} \to T_B^2$ trivializes.
Let $\eta$ be a small positive number and
 let $\phi\big|_V^{-1}(U_{\eta})$ be a set defined by
$$
\phi\big|_V^{-1}(U_{\eta}) =
          \phi\big|_V^{-1} \bigg( p \bigg(\left( \begin{array}{cc} 
                              a & 0 \\
                              0 & -a    \end{array} \right) d\overline{z}
                           - \left( \begin{array}{cc}
                              a & 0 \\
                              0 & -a  \end{array} \right)^{\dag} dz \bigg) 
                                \bigg| a \in [0, \frac{\pi}{2})
                                           \times i[0, \frac{\pi}{2}),
                               |a - \alpha| < \eta \bigg),
$$
 here $\alpha = \frac{m}{4} + \frac{in}{4}$,
 $m, n = 0, 1, 2$ (see remark 4.4)
 and $p : \T \to Rep_{T^2}(SO(3))$ is the quotient map.
Let $T_B^2 = \cup_i V_i$ be a finite covering by
 these disks.
We set
$$
\phi^{-1}(U_{\eta}) = \bigcup_i \phi\big|_{V_i}^{-1}(U_{\eta}).
$$ 
If $\eta$ is sufficiently small,
 it will be a small neighbourhood of $S_2$.

On the other hand, 
$$
\mu_{\nu}(U) = \int_{\pi^{-1}(U)} |F_{\Xi_{\ep_{\nu}}}|^2
 d\mu_{\widehat{M_{\ep}}},
$$
 here $U$ is an open subset of the base,
 and $\pi: \widehat M_{\ep_{\nu}} \to T^2_B$
 gives a sequence of measures on the base.
Since this is bounded in the space of measures,
 there is a limit measure $\mu$ of this sequence.

With this preparation,
 we define $S_{3, \eta}$ as follows.\\
$$
S_{3, \eta} = \{x \in T_B^2 
                 | \mu(x) > \delta_{\eta} \}, 
$$
 here $\delta_{\eta}$ is the constant which gives the energy bound for
 the doubly periodic instantons mapped into 
 $Rep_{T^2}(SU(2)) \setminus
    U_{\eta}$ by the map $\rho$ of theorem 5.3.
Obviously, $S_{3, \eta}$ is a finite subset.
It is also clear that $S_{3, \eta}$ contains all type one bubbles
 and type two bubbles satisfying the above condition.

\begin{defn}
$S_{\eta} = S_1 \cup \phi^{-1}(U_{\eta}) \cup S_{3, \eta}.$
\end{defn}




Let us rescale the metric of $\widehat{M}_{\ep_{\nu}}$ by
 $\frac{1}{\ep_{\nu}}$.
We define the connection
 $\Xi_{\nu} = A_{\nu} + \Phi{_\nu} dS_{\nu} + \Psi_{\nu} dT_{\nu}$
 on the rescaled space by
$$
A_{\nu}(x) = A_{\ep_{\nu}}(x_0 + \ep_{\nu} x),
$$
$$
\Phi_{\nu}(x) = \ep_{\nu} \Phi_{\ep_{\nu}}(x_0 + \ep_{\nu} x), \;\;
\Psi_{\nu}(x) = \ep_{\nu} \Psi_{\ep_{\nu}}(x_0 + \ep_{\nu} x),
$$
 where $x_0$ is some point on  $\widehat{M}_{\ep_{\nu}}$
 and $x$ is a coordinate on $\widehat{M}_{\ep_{\nu}}$.
Here $S_{\nu}$ and $T_{\nu}$ are $\frac {1}{\ep_{\nu}}$-rescaled
 coordinates of ($S, T$) on the base.
We observe the following.
\begin{lem}
The relations between the component of the
 curvatures on rescaled and non-rescaled
 manifolds are given as follows.
$$
(i) \pll F_{A_{\nu}}\pll_{L^{\infty}(\frac{1}{\ep_{\nu}}
 \widehat{M_{\ep_{\nu}}})} 
 = \pll F_{A_{\ep_{\nu}}} (x_0 + \ep_{\nu} x)\pll_{L^{\infty}(\widehat{M_1})},
$$
$$
(ii) \ep_{\nu}^{-1} \pll d_{A_{\nu}} \Psi_{\nu}
   - \del_{T_{\nu}} A_{\nu} \pll_{L^{\infty}(\frac{1}{\ep_{\nu}}
 \widehat{M_{\ep_{\nu}}})}
 = \pll d_{A_{\ep_{\nu}}} \Psi_{\ep_{\nu}}
  - \del_T A_{\ep_{\nu}} \pll_{L^{\infty}(\widehat{M_1})},
$$
and
$$
(iii)  \ep_{\nu}^{-2} \pll \del_{T_{\nu}} \Phi_{\nu} - \del_{S_{\nu}} \Psi_{\nu}
   - [\Phi_{\nu}, \Psi_{\nu}]\pll_{L^{\infty}(\frac{1}{\ep_{\nu}}
     \widehat{M_{\ep_{\nu}}})}
 = \pll \del_T \Phi_{\ep_{\nu}} - \del_S \Psi_{\ep_{\nu}}
    - [\Phi_{\ep_{\nu}}, \Psi_{\ep_{\nu}}] \pll_{L^{\infty}(\widehat{M_1})}.
$$
Here the left hand sides are the norms with respect to the
 rescaled metric (that is,
 $\frac{1}{\ep_{\nu}^2} g_{\widehat{M_{\ep_{\nu}}}}$)
 and the right hand sides are the norms
 with respect to the metric of $\widehat{M}_1$. \qed
\end{lem}
Hereafter we denote $(S_{\nu}, T_{\nu})$ simply by
 $(S, T)$ for notational simplicity.

By the characterization of the third type of the bubbles (see 
 the proof of theorem 3.1),  
 at all points on 
 $\pi^{-1}(T_B^2 \setminus S_{\eta})$,
 the $L^{\infty}$-norm (with respect to the metric on
 $\widehat{M}_1$) of the fiber part of the curvature $F_{\Xi_{\ep_{\nu}}}$ 
 and the term
 $\ep_{\nu}\pll d_{A_{\ep_{\nu}}} \Psi_{\ep_{\nu}}
   - \del_T A_{\ep_{\nu}} \pll_{L^{\infty}}$ converge to zero.  
%
%
%
Now we can state the main theorem of this section.
\bt
There is a subsequence of $\Xi_{\ep_{\nu}}$ which converges to
 a limit connection in $C^{\infty}$ sense (with respect to the
 metric on $\widehat M_1$) locally on $\pi^{-1}(T^2_B \setminus S_{\eta})$. 
\et
\begin{rem}
In the case when the bundle is non-trivial on fibers
 and the four manifold is a product of
 Riemannian surfaces $\Sigma_h \times \Sigma_g$,
 $g \geq 2$, 
 this is essentially done by Fukaya \cite{F1}, mainlemma 4.42.
However, as in the case of the second type
 of bubbles of the previous section, 
 in our case the triviality of the bundle restricted on each fiber
 causes additional difficulty
 coming from the existence of the kernel of the operator $d_{A_{\nu}}$
 when $A_{\nu}$ is flat.
\end{rem}
Once this theorem is proved, taking a sequence 
 $\eta_i$ of positive numbers converging to zero,
 we can see the convergence of the connections on
 $\pi^{-1}(T^2_B \setminus S)$ by the diagonal argument,
 where the set $S$ is defined as follows.
\begin{defn}
$S = S_1 \cup S_2 \cup \bigcup_i S_{3, \eta_i}$.
\end{defn}
Note that because of theorem 6.6, $S_{3, \eta_i} \setminus S_{3, \eta_j}$
 for $\eta_j < \eta_i$ is contained in 
 $\phi^{-1} (U_{\eta_i} \setminus U_{\eta_j})$.
So only the points in $S_2$ can be the accumulation points
 of $S$, as required.

The rest of this section is devoted to the proof of this theorem.

Let $x$ be an arbitrary point on 
 $T_B^2 \setminus S_{\eta}$.
Let $B_{\iota}(x)$, with $\iota$ a small positive number,
 be an $\iota-$neighbourhood of $x$
 in $T_B^2$ and we assume it is contained in
 $T_B^2 \setminus S_{\eta}$.
We will consider the region $\pi^{-1}(B_{\iota}(x))$ with 
 the $\frac{1}{\ep_{\nu}}$-rescaled metric.
We write this rescaled region as
 $(\frac{1}{\ep_{\nu}}B_{\iota}(x) \times T^2,
 \frac{1}{\ep^2_{\nu}}g_{\ep_{\nu}})$,
 and abbreviate it as $\frac{1}{\ep_{\nu}}B_{\iota}(x) \times T^2$,
 including the information of the metric.
Let $x' \in \frac{1}{\ep_{\nu}}B_{\iota}(x)$ and let
 $B_1(x')$ be the unit ball around $x'$.
We assume this is included in $\frac{1}{\ep_{\nu}}B_{\iota}(x)$.
 
By theorem 5.1, for sufficiently large $\nu$,
 there is a flat connection of the form 
$
\Xi_0 = A_0 + 0 dS + 0 dT,
$
 where $A_0 \in \mathfrak{t}$ is a constant matrix valued one form
 which is not a singular point of $Rep_{T^2}(SO(3))$,
 and a unitary gauge transformation $g_{0, \nu}$ such that 
$$
\pll g^*_{0, \nu} \Xi_{\nu} - \Xi_0 \pll_{C^r(B_1(x') \times T^2)} < C \ep_0^{\lam}
$$
 holds, with $\lam$ a positive number depending only on 
 the metric on $B_1(x') \times T^2$.
Moreover, we can assume $A_0$ does not
 depend on $\nu$, because
 $[A_{\ep_{\nu}}]$ is a converging sequence of holomorphic curves and by 
 $\frac{1}{\ep_{\nu}}$ rescaling, the differentials of these maps become
 very small.
Here $\ep_0$ is a constant which bounds the $L^2$ norms of
 the curvature of $\Xi_{\ep_{\nu}}$ over $B_{\iota}(x) \times T^2$.
By  definition of the excepted set $S_{\eta}$, 
 $\ep_0$ can be taken to be smaller than $\delta_{\eta}$.
In particular, it is much smaller than $\eta$.

Put
$$
g^*_{0, \nu}\Xi_{\nu} = A_{\nu}' + \Phi_{\nu}'dS + \Psi_{\nu}'dT.  
$$
We can assume the complex gauge equivalence classes of 
 $A'_{\nu}$ are in $Rep_{T^2}(SO(3)) \setminus \phi^{-1}(U_{\eta})$
 for sufficiently large $\nu$.

Now we will globalize this construction by gluing.
Let $\widetilde{A_{\nu}}(S, T)$ be a $\mathfrak{t}$-valued section
 which is complex gauge equivalent to $A'_{\nu}(S, T)$
 on each fiber. 
\begin{prop}
There is a unitary gauge transformation $g_{\nu}$ on $\frac{1}{\ep_{\nu}} B_{\iota}(x)$
 and a $\mathfrak{t}$-valued section $A_0(S, T)$
 on $\frac{1}{\ep_{\nu}} B_{\iota}(x)$ with
$$
\pll A_0(S, T) - \widetilde{A_{\nu}}(S, T)
 \pll_{C^r(\frac{1}{\ep_{\nu}} B_{\iota}(x))} < C \frac{\ep_0^{\lam}}{\eta},
$$
 such that
$$
\pll g_{\nu}^* \Xi_{\nu} - A_0(S, T) \pll_{C^r(\frac{1}{\ep_{\nu}} B_{\iota}(x))}
   < C' \frac{\ep_0^{\lam}}{\eta}
$$
 holds with $C$, $C$'s and $A_0(s, t)$ do not depend on $\nu$.
\end{prop}
Proof:
Cover $\frac{1}{\ep_{\nu}} B_{\iota}(x)$ by
 finite number of open disks $\{ B_i \}_{i \in I}$ of radius one
 (later we will mention how they should be taken),
 $I$ is an index set (having finite elements).
On each open disk, take a gauge in which the connection
 matrix become close to a flat connection of the form 
 $\Xi_i = A_i + 0 dS + 0 dT$ as above, $A_i \in \mathfrak{t}$
 does not depend on $S, T$ and $\nu$.
We write it as
$$
g_{i, \nu}^* \Xi_{\nu} = A'_{\nu, i} + \Phi'_{\nu, i}dS + \Psi'_{\nu, i}dT
$$
 with
$$
\pll g_{i, \nu}^* \Xi_{\nu} - \Xi_i \pll_{C^r(B_i \times T^2)} < C \ep_0^{\lam}.
$$
\begin{claim}
Denoting the family of connections in $\mathfrak{t}$
 which are complex gauge equivalent to 
 $A'_{\nu, i}(S, T)$ by $\widetilde{A_{\nu, i}}(S, T)$, we have
$$
\pll \widetilde{A_{\nu, i}}(S, T) - A_i \pll_{C^r(B_i \times T^2)}
         < C \frac{\ep_0^{\lam}}{\eta},
$$ 
 $C$ does not depend on $\nu$.
\end{claim}
\noindent
\it{Proof of the claim.} \rm
As noted above, before rescaling $[A_{\ep_{\nu}}]$ are in a 
 sequence converging to a holomorphic curve $\phi$.
So for sufficiently large $\nu$, the norm
 $\pll [A_{\ep_{\nu}}] \pll_{C^r(B_{\iota}(x))}$
 ($[A_{\ep_{\nu}}]$ is seen as an element of $\T$)
 is bounded uniformly with respect to $\nu$.
So $\pll \widetilde{A_{\nu, i}}(S, T) -
 \widetilde{A_{\nu, i}}(S_0, T_0)
           \pll_{C^r(\frac{1}{\ep_{\nu}} B_i \times T^2)}$, 
 $(S_0, T_0) \in B_i$ a fixed point,
 is uniformly bounded with the upper 
 bound of order $\ep_{\nu}$
 (since $\widetilde{A_{\nu, i}}$ is constant in the fiber direction,
 it is an estimate for the horizontal direction).
On the other hand,
 since
 $\pll g_{i, \nu}^* \Xi_{\nu} - \Xi_i\pll_{C^r(B_i \times T^2)} < C \ep_0^{\lam}$
 and $\ep_0$ is small enough, $A'_{\nu, i}(S_0, T_0)$
 is complex gauge equivalent to 
 a connection in $\mathfrak{t}$ (it is $\widetilde{A_{\nu, i}}(S_0, T_0)$
 by definition), 
 whose $C^r$ distance from $A_i$ is less than $C \ep_0^{\lam} /\eta$
 by corollary 4.15.
It follows that 
$$
\pll \widetilde{A_{\nu, i}}(S, T) - A_i \pll_{C^r(B_i \times T^2)}
 < C \frac{\ep_0^{\lam}}{\eta},
$$
 as required. \qed \\
\begin{rem}
From this discussion, it follows that the difference between
 $A_i$ and $A_j$ for $B_i \cap B_j \neq \phi$
 is bounded by $C \frac{\ep_0^{\lam}}{\eta}$.
\end{rem}
Now suppose $B_{i_1}$ and $B_{i_2}$ intersects, and compare
 the gauges.
Let $g(S, T)$ be a gauge transformation on the rank two trivial bundle on
 $T^2$ which transforms 
 $A'_{\nu, i_1}(S, T)$ to $A'_{\nu, i_2}(S, T)$.
If $g'(S, T)$ is another gauge transformation which transforms
 $A'_{\nu, i_1}(S, T)$ to $A'_{\nu, i_2}(S, T)$,
 the 
 difference of these two is given by the isotropies of
 $A'_{\nu, i_1}(S, T)$ or $A'_{\nu, i_2}(S, T)$.
It is trivial when the holonomy is not degenerated 
 on $T^2$ (that is, when the holonomy group generates
 $SO(3)$), and $S^1$ or $SO(3)$ when it is degenerated
 ($SO(3)$ case does not occur since we are working on the complement of 
 $\pi^{-1}(U_\eta)$).  
Because the connection is ASD, the holonomy varies analytically
 when we fix a loop on the fiber.
In particular, if it is not degenerated
 at a point, it is not degenerated
 except measure zero loci.
We treat degenerated and non-degenerated cases separately.\\

\noindent
(i) the non-degenerated case \\
 
In this case, the difference between the gauges $g$
 on $B_{i_1}$ and $B_{i_2}$ are uniquely determined 
 by the gauge transformation which sends
 $A'_{\nu, i_1}(S, T)$ to $A'_{\nu, i_2}(S, T)$
 on the open loci where the holonomies are not degenerate.
This transformation uniquely smoothly extends to the whole intersection.
And the difference between $g$ and $1$
 is $C^0$-small on $(B_{i_1} \cap B_{i_2}) \times T^2$ because on each fiber 
 this difference is estimated 
 as
$$
\pll g - 1\pll_{C^r(T^2)} < C \frac{\ep_0^{\lam}}{\eta},
$$
 using the fact that $A'_{\nu, i_1}(S, T)$ and $A'_{\nu, i_2}(S, T)$
 are near to the same connection on $\mathfrak{t} \setminus U_{\eta}$.
Moreover, as the inequality shows, 
 it is $C^r$-small in the direction of the fibers.
On the other hand, we have, by definition of $g$, 
$$
g^*(g^*_{i_1, \nu} \Xi_{\nu}) = g^*_{i_2, \nu} \Xi_{\nu}.
$$
The part having differential forms of the horizontal
 directions of this equation can be written in the form that  
$$
g^{-1} \del_S g dS + g^{-1} \del_T g dT = - g^{-1} \Phi'_{\nu, i_1} g dS 
  - g^{-1} \Psi'_{\nu, i_1} g dT +
  \Phi'_{\nu, i_2} dS +  \Psi'_{\nu, i_2} dT.
$$
Since $g^{-1} \Phi'_{\nu, i_1} g$ satisfies
$$
\begin{array}{l}
\pll g^{-1} \Phi'_{\nu, i_1} g\pll_{C^r((B_{i_1} \cap B_{i_2}) \times T^2)} \\
  \hs{.5in} < \sum_{0 \leq i, j, i+j \leq r}
  \pll g\pll_{C^i((B_{i_1} \cap B_{i_2}) \times T^2)}
 \pll g^{-1}\pll_{C^j((B_{i_1} \cap B_{i_2}) \times T^2)}
 \pll \Phi'_{\nu, i_1}\pll_{C^{r-i-j}(B_{i_1} \times T^2)},
\end{array}
$$
 (and similarly for $g^{-1} \Psi'_{\nu, i_1} g$) we have
$$
\begin{array}{l}
\pll g^{-1} \del_S g dS
      + g^{-1} \del_T g dT\pll_{C^r((B_{i_1} \cap B_{i_2}) \times T^2)} \\
 \hs{.5in} < 2 \ep_0^{\lam} (\sum_{0 \leq i, j \leq r}
  \pll g\pll_{C^i((B_{i_1} \cap B_{i_2}) \times T^2)}
 \pll g^{-1}\pll_{C^j((B_{i_1} \cap B_{i_2}) \times T^2)} + 2).
\end{array}
$$
From this, we can inductively deduce the estimate
$$
\pll g - 1\pll_{C^r((B_{i_1} \cap B_{i_2}) \times T^2)}
   < C \frac{\ep_0^{\lam}}{\eta}
$$
 (this time the norm is not restricted to the fiber direction).

Now we can glue the gauges on $B_{i_1}$ and $B_{i_2}$
 using a cut-off function.
First, extend the gauge transformation $g$ to $B_{i_1} \times T^2$
 in a way that it satisfies
$$
\pll g - 1\pll_{C^r(B_{i_1} \times T^2)} < C \frac{\ep_0^{\lam}}{\eta}.
$$
Introduce a smooth function $\rho$ on $B_{i_1}$ which satisfies
$$
\begin{array}{l}
1. 0 \leq \rho \leq 1, \\
2. \pll \rho \pll_{C^{\infty}} < 10,\\
3. \rho = 1 \; on \; B_{i_1} \cap B_{i_2}, \\
4. \rho = 0 \; on \; B'_{i_1},
\end{array}
$$
 here $B'_{i_1}$ is defined by
$$
B'_{i_1} = \{ x \in B_{i_1} \big|  d(x, B_{i_2}) > \frac{1}{5} \}.
$$
Denote $g$ in the form $\exp(e)$.
Then define a gauge on $B_{i_1} \cup B_{i_2}$
 by the following gauge transformation $g'$,
$$
g' = \left\{ \begin{array}{l}
      g_{i_2, \nu} \; on \; B_{i_2}\\
      g_{i_1, \nu} \exp(\rho e) \; on \; B_{i_1}.
       \end{array} \right.
$$
Since the base $\frac{1}{\ep_{\nu}}B_{\iota}(x)$
 has no nontrivial topology,
 by taking a covering of $\frac{1}{\ep_{\nu}}B_{\iota}(x)$
 in a reasonable way (namely, there are not too much overlappings,
 which we will not discuss precisely here, but I hope the meaning is clear)
 and repeating the above procedure,
 we can form a gauge on the whole $\frac{1}{\ep_{\nu}}B_{\iota}(x)$
 satisfying the required estimate of the proposition, with 
 $A_0(S, T)$ defined from $A_i$'s using decomposition of the unit,
 noting that the differences of $A_i$'s on the overlappings of 
 the covering are small (remark 6.11).\\
%

\noindent
(ii) the degenerate case\\

In this case, $A'_{\nu, i_1}, A'_{\nu, i_2}$ are unitary gauge equivalent
 to a connection in $\mathfrak{t}$.
Since we are dealing with connections on
 $T_B^2 \setminus (S_{\eta})$,
 we can transform $A'_{\nu}$ to $\mathfrak{t}$ by a gauge transformation
 satisfying 
$$
\pll g - 1\pll_{C^r(B_i \times T^2)} < C \frac{\ep_0^{\lam}}{\eta}.
$$
So we can assume from the first $A'_{\nu, i}$'s are in $\mathfrak{t}$.

Then, $A'_{\nu, i_1}$ and $A'_{\nu, i_2}$ are necessarily the same.
In particular, the difference of the gauge is only by
 a gauge transformation which does not change the fiber component,
 that is, by a gauge transformation which can be written in the form
$$
\exp(e), \; e = \left( \begin{array}{cc}
                         a & 0 \\
                         0 & -a \end{array} \right), 
$$
 here $a$ is a pure imaginary and constant along the fibers.
On the other hand, by the ASD equation, $\Phi'_{\nu}$ and $\Psi'_{\nu}$ 
 are also of diagonal form and constant in the direction of the
 fibers (see lemma 6.13 below).

Using the equation
$$
\exp(e)^* (\Phi'_{\nu, i_1}dS + \Psi'_{\nu, i_1}dT)
   = \Phi'_{\nu, i_2}dS + \Psi'_{\nu, i_2}dT,
$$
 we conclude as before that $g = \exp(e)$ is of the form that
$$
g = g_0 g_1,
$$
 where $g_0$ is a constant diagonal matrix and $g_1$ satisfies the estimate
$$
\pll g_1 - 1\pll_{C^r(B_i \times T^2)} < C \frac{\ep_0^{\lam}}{\eta}.
$$ 
Since $g_0$ does not transform the connection at all, 
 it can be absorbed to the gauge on one of the 
 open covers whose intersection we are considering.
Then as before, we can construct a global gauge satisfying the required
 estimate of the proposition by gluing. \qed \\

\begin{defn}
Let $\mathcal{K}$ be the space of constant sections of $2 \times 2$,
 diagonal, trace free pure imaginary matrices over $T^2$.
\end{defn}
Let $A_0 \in Rep_{T^2}(SO(3)) \setminus p(U_{\eta})$.
Let $U$ be a small neighbourhood of $A_0$ in $Rep_{T^2}(SO(3))$.
Let $V$, $W$, be neighbourhoods of $0$ in the space of 
 traceless $C^{\infty}$ skew-hermitian and hermitian $2 \times 2$ matrix
 valued functions on $T^2$ which are perpendicular to 
 $\mathcal{K}$ or
 $\left\{ \left( \begin{array}{cc}
        c & 0 \\
        0 & -c \end{array} \right)
      \bigg| c \in \R \right\}$, respectively.
The inner products are defined by $<A, B> = \int_{T^2} Tr AB d\mu$
 in both cases.
This condition is the same as saying the integrals
 of the diagonal components are zero.

Let $\mathcal{A}$ be the space of connections on the rank two
 trivial bundle on $T^2$.
Then the map
$$
\Gamma : U \times V \times W \to \mathcal{A}
$$
 defined by
$$
\Gamma(a, v, w) = \exp(w)^* \exp(v)^* (a)
$$
 is a local diffeomorphism (corollary 4.9).
In particular, we can transform $A_{\nu}'$,
 here $g_{\nu}^* \Xi_{\nu} = A_{\nu}' + \Phi_{\nu}' dS
 + \Psi_{\nu}' dT$ to $\mathfrak{t}$
 by a complex gauge transformation $h$ of this form.
Denote the resulting connection as 
$$
\widetilde{\Xi_{\nu}} = \widetilde{A_{\nu}}
 + \widetilde{\Phi_{\nu}}dS + \widetilde{\Psi_{\nu}}dT.
$$
Since the first one of the ASD equations 
 holds for $\widetilde{\Xi_{\nu}}$, we can see the following
 (in fact we have already used it several times).
\begin{lem}
$\widetilde{\Phi_{\nu}}$ and $\widetilde{\Psi_{\nu}}$
 are diagonal matrices, which are constant along the fibers.
\end{lem}
Proof:
Since $\widetilde{A_{\nu}}$ is flat, we can apply
 the Hodge decomposition. 
$\del_T \widetilde{A_{\nu}}$ and
 $*_{T^2_F} \del_S \widetilde{A_{\nu}}$ belong to the harmonic
 part, so they are perpendicular to
 $d_{\widetilde{A_{\nu}}} \widetilde{\Psi_{\nu}}$.
Moreover $*_{T^2_F} d_{\widetilde{A_{\nu}}} \Phi_{\nu}$,
 and $d_{\widetilde{A_{\nu}}} \widetilde{\Psi_{\nu}}$
 and $*_{T^2_F} d_{\widetilde{A_{\nu}}} \Phi_{\nu}$ are
 mutually perpendicular, too.
So the equation reduces to 
$$
\del_T \widetilde{A_{\nu}} + *_{T^2_F} \del_S \widetilde{A_{\nu}} = 0,
$$
$$
d_{\widetilde{A_{\nu}}} \widetilde{\Psi_{\nu}} = 0,
 \;\; d_{\widetilde{A_{\nu}}} \widetilde{\Phi_{\nu}} = 0.
$$
The kernel of $d_{\widetilde{A}}$ is the set of constant, diagonal matrices.
So the lemma follows. \qed\\

We further modify the gauge for our purpose.
We have seen that we can complex gauge transform
 $g_{\nu}^*\Xi_{\nu}$ to a connection of the form
$$
\widetilde{\Xi_{\nu}} = \widetilde{A_{\nu}} + \widetilde{\Phi_{\nu}}dS 
                     + \widetilde{\Psi_{\nu}}dT,
$$
 where $\widetilde{A_{\nu}}$, $\widetilde{\Phi_{\nu}}$ and
 $\widetilde{\Psi_{\nu}}$ are all diagonal and constant in the
 direction of the fibers (in other words, in $\mathcal{K}$).
Moreover, we can require this complex gauge transformation
 (we write it by $h$) to be of the
 form 
$$
h = \exp(f_2) \exp(- f_1), \; f_1 \in W, \; f_2 \in V.
$$
With this requirement, the transformation $h$ is uniquely determined,
 by corollary 4.9.
It satisfies the estimates
$$
\pll f_i \pll_{C^r(\frac{1}{\ep_{\nu}}B_{\iota}(x) \times T^2)}
 < C \frac{\ep_0^{\lam}}{\eta},
 \; i = 1, 2.
$$

%
\begin{defn}
$\Xi_{\nu, G} = \exp(f_2)^* g_{\nu}^* \Xi_{\nu}$.
\end{defn}

In other words, $\Xi_{\nu, G} = \exp(f_1)^* \widetilde{\Xi_{\nu}}$.
Let us write $\Xi_{\nu, G}$ as $\Xi_{\nu, G} =
     A_{\nu, G}(S, T) + \Phi_{\nu, G}(S, T)dS
   + \Psi_{\nu, G}(S, T)dT$.
\begin{cor}
The estimate
$$
\pll A_{\nu, G}(S, T) + \Phi_{\nu, G}(S, T)dS  + \Psi_{\nu, G}(S, T)dT
  - \widetilde{A_{\nu}}(S, T)\pll_{C^r(\frac{1}{\ep_{\nu}}B(\iota) \times T^2)}
           < C \frac{\ep_0^{\lam}}{\eta}
$$
 holds. 
\end{cor}
\proof
This is obvious from proposition 6.9 and definition of
 $\Xi_{\nu, G}$.\qed\\

Having constructed a moderate gauge, we start the analysis using
 the ASD equation to deduce stronger estimates for the components
 of $\Xi_{\nu, G}$.
We put the following convention to avoid non-essential complexity.
\begin{convention}
Although generally
 the metrics of the fibers of $\widehat M_1$
 vary (so do the complex structures),
 we treat in the following calculation as if they were constant.
This is harmless because we are rescaling a small open 
 subset $B_{\iota}$ and so the metrics of the
 fibers do not really largely vary.
So it does not cause differences on the conclusions.
Let $x + iy$ be the complex coordinate of the fiber.
\end{convention}
We use the 
$$
\del_T \Phi_{\nu} - \del_S \Psi_{\nu} - [\Phi_{\nu}, \Psi_{\nu}]
 - f(S, T)*_{T_F}F_{A_{\nu}} = 0
$$
 part of the equation.
We can decompose it to a part which is perpendicular
 (in $L^2$ sense) to matrices in $\mathcal{K}$
 on each fiber 
 and a part which is proportional to these matrices on each fiber.
We denote the former part of the equation by ($E_1$) and
 the latter part by ($E_2$).
Correspondingly, we write $\Xi_{\nu,G} = \Xi_{1, \nu, G} + \Xi_{2, \nu, G}$
 and $A_{\nu, G} = A_{1, \nu, G} + A_{2, \nu, G}$ etc.,
 which are $L^2$ orthogonal decompositions on each fiber.
Namely,
$$
(E_1) \;\; \del_T \Phi_{1, \nu, G} - \del_S \Psi_{1, \nu, G}
       - [\Phi_{\nu, G}, \Psi_{\nu, G}]_1
       - f(S, T) *_{T_F} (F_{A_{\nu, G}})_1 = 0,
$$
$$
(E_2) \;\; \del_T \Phi_{2, \nu, G} - \del_S \Psi_{2, \nu, G}
       - [\Phi_{1, \nu, G}, \Psi_{1, \nu, G}]_2
       - f(S, T) *_{T_F} (F_{A_{\nu, G}})_2 = 0,
$$
 here $(\cdots)_i$, $i = 1, 2$
 respectively denote the $(\mathcal{K})^{\perp}$ and $\mathcal{K}$
 parts of the corresponding object.

We introduce the following notations.
As before, we denote the hermitian gauge transformation 
 which transforms $\widetilde{\Xi_{\nu}}$ to $\Xi_{\nu,G}$
 by $h_1 = \exp(f_1)$.
We put
$$
e_{a, b} = \frac{\del^{a+b}f_1}{\del S^a \del T^b},
$$
$$
f_{a, b}(S, T) = \pll e_{a, b}(S, T) \pll^2_{L^2(T_F)}.
$$

In terms of $f_1$, $\Xi_{\nu, G}$ is given in the following form.
$$
\begin{array}{lll}
\Xi_{\nu, G} &
  = \exp(f_1)^* (\widetilde{A_{\nu}} + \widetilde{\Phi_{\nu}} dS
                    + \widetilde{\Psi_{\nu}} dT) & \\
  & = \exp (-f_1) \overline{\del} \exp (f_1)
      + \exp (-f_1) \widetilde{\Xi_{\nu}^{0, 1}} \exp (f_1) & \\
  &  \hs{1in}    -  (\exp (-f_1) \overline{\del} \exp (f_1))^{\dag} 
      - (\exp (-f_1) \widetilde{\Xi_{\nu}^{0, 1}} \exp (f_1))^{\dag} & \\
  & = -i \del_S f_1 dT + i \del_T f_1 dS + \frac{1}{2}[\del_S f_1, f_1]dS
                                   + \frac{1}{2}[\del_T f_1, f_1]dT & \\
  &  \;\;\;\; + \frac{i}{2}(- f_1 \del_T f_1 - (\del_T f_1) f_1)dS
                + \frac{i}{2}(f_1 \del_S f_1 + (\del_S f_1) f_1)dT & \\
  &  \;\;\;\;    + i[\widetilde{\Psi_{\nu}}, f_1] dS
           - i[\widetilde{\Phi_{\nu}}, f_1] dT 
           - i\del_x f_1 dy + i \del_y f_1 dx & \\
  &  \;\;\;\;    + \frac{1}{2}[\del_x f_1, f_1]dx
                                   + \frac{1}{2}[\del_y f_1, f_1]dy & \\
  &  \;\;\;\; + \frac{i}{2}(- f_1 \del_y f_1 - (\del_y f_1) f_1)dx
                + \frac{i}{2}(f_1 \del_x f_1 + (\del_x f_1) f_1)dy & \\
  &  \;\;\;\;  - i[\widetilde{A_{\nu, x}}, f_1] dy 
                 + i[\widetilde{A_{\nu, y}}, f_1] dx
           + \widetilde{\Xi_{\nu}}
                + A(\widetilde{\Xi_{\nu}^{0, 1}}, f_1, \overline{\del} f_1)
             - A(\widetilde{\Xi_{\nu}^{0, 1}}, f_1,
                     \overline{\del} f_1)^{\dag}, & \\
\end{array}
$$
 here $\widetilde{A_{\nu}} =
 \widetilde{A_{\nu, x}}dx + \widetilde{A_{\nu, y}}dy$,
 $A(\widetilde{\Xi_{\nu}^{0, 1}}, f_1, \overline{\del} f_1)$ 
 is a one form consisted from 
 $\widetilde{\Xi_{\nu}^{0, 1}}$, $f_1$ and $\overline{\del} f_1$,
 each term of which is quadratic
 (or more) with respect to $f_1$ (with no differential).
 
The relevant part of the curvature,
 that is, purely horizontal or vertical two form parts of the curvature 
 (we write it as $F_{\Xi_{\nu, G}}^{HV}$) is then given as follows.
$$
\begin{array}{lll}
F_{\Xi_{\nu, G}}^{HV}
 & = -i (\del_S^2 + \del_T^2) f_1 dSdT
     + \frac{1}{2}(\del_S [\del_T f_1, f_1] - \del_T [\del_S f_1, f_1])dSdT & \\
 &  \;\;\;\;  + \frac{i}{2} \del_T(f_1 \del_T f_1 + (\del_T f_1) f_1)dSdT
               + \frac{i}{2} \del_S(f_1 \del_S f_1 + (\del_S f_1) f_1)dSdT & \\

 &  \;\;\;\;    - i(\del_S [\widetilde{\Phi_{\nu}}, f_1]
                    + \del_T [\widetilde{\Psi_{\nu}}, f_1]) dSdT 
              + (\del_S \widetilde{\Psi_{\nu}}
                 - \del_T \widetilde{\Phi_{\nu}}) dSdT
              + [\del_T f_1, \del_S f_1] dSdT  & \\
 &  \;\;\;\; + [\widetilde{\Phi_{\nu}},
                 -i \del_T f_1] dSdT 
             + [i \del_S f_1,
                    \widetilde{\Psi_{\nu}}] dSdT & \\
 &  \;\;\;\;  + i d_{\widetilde{A_{\nu}}} *_{T_F} d_{\widetilde{A_{\nu}}} f_1
              + [\del_y f_1, \del_x f_1] dxdy
              + \frac{1}{2}(\del_x [\del_y f_1, f_1] - \del_y [\del_x f_1, f_1])dxdy & \\
 & \;\;\;\;    + \frac{i}{2} \del_y(f_1 \del_y f_1 + (\del_y f_1) f_1)dxdy
               + \frac{i}{2} \del_x(f_1 \del_x f_1 + (\del_x f_1) f_1)dxdy & \\

 &  \;\;\;\;  + Q(\Xi, f_1), & \\ 
\end{array}
$$
 here $Q(\Xi, f_1)$ is the sum of the terms which contain $f_1$ as a factor.

For hermitian matrix valued functions on $T^2$ 
 with small $C^r$-norms, we can deduce 
 estimates for the $L^4_1$ norms by the $L^2$ norms.
\begin{lem}
Let $f$ be a smooth hermitian matrix valued
 function on $T^2$ with
 $\pll f \pll_{C^r(T^2)} < \ep$, $\ep$ a small positive number.
Then the estimate
$$
\pll \nabla f \pll_{L^4(T^2)} < C \ep^{\frac{1}{2}}
        \pll f \pll_{L^2(T^2)}^{\frac{1}{2}},
$$
 holds.
\end{lem}
\proof
Let $x, y$ be the standard (multivalued) affine coordinate
 on $T^2$.
Let $< , >$ denote the $L^2$ inner product of hermitian matrix
 valued functions on $T^2$ defined by
$$
<A, B> = \int_{T^2} tr AB d\mu.
$$
We can prove hermitian matrix valued version of Schwarz's inequality:
$$
| <A, B> | \leq <A, A>^{\frac{1}{2}} <B, B>^{\frac{1}{2}}, 
$$ 
 by a straightforward analysis.
Then,
$$
\begin{array}{ll}
\int_{T^2} tr (\del_x f)^4 d \mu 
   & = <\del_x f, (\del_x f)^3> \\
   & = - <f, 3 (\del^2_x f)(\del_x f)^2> \\
   & < 3 \ep <f, f>^{\frac{1}{2}} (\int_{T^2} tr (\del_x f)^4 d \mu)^{\frac{1}{2}}\\
   & < 9 \ep^{\sum_{i=0}^{\infty} 2^{-i}}
          <f, f>^{\sum_{i=1}^{\infty} 2^{-i}} \\
   & = 9 \ep^2 <f, f>.
\end{array}
$$
The same holds for $\int_{T^2} (\del_y f)^4 d \mu$.
It is a straightforward calculation to see that for hermitian matrix
 valued sections, $\pll A \pll_{L^4}^4$ and $\int tr A^4 d\mu$ gives the 
 equivalent norms.
So we have an estimate
$$
\pll \nabla f \pll_{L^4(T^2)}
 < C \ep^{\frac{1}{2}} \pll f \pll_{L^2(T^2)}^{\frac{1}{2}},
$$
 as required.
 \qed\\

We first consider the equation ($E_1$).
Explicitly, it is given as follows.
$$
\begin{array}{l}
-i (\del_S^2 + \del_T^2) f_1
     + \frac{1}{2}(\del_S [\del_T f_1, f_1] - \del_T [\del_S f_1, f_1])_1 
      - i(\del_S [\widetilde{\Phi_{\nu}}, f_1] 
                    + \del_T [\widetilde{\Psi_{\nu}}, f_1]) \\
    + \frac{i}{2} \del_T(f_1 \del_T f_1 + (\del_T f_1) f_1)_1
               + \frac{i}{2} \del_S(f_1 \del_S f_1 + (\del_S f_1) f_1)_1 
    + [\widetilde{\Phi_{\nu}},
                 -i \del_T f_1]  
             + [i \del_S f_1,
                    \widetilde{\Psi_{\nu}}]   \\
  + [\del_T f_1, \del_S f_1]_1  
   + f(S, T)\{i *_{T_F} d_{\widetilde{A_{\nu}}} *_{T_F} d_{\widetilde{A_{\nu}}} f_1
             + [\del_y f_1, \del_x f_1]_1 \\
 + \frac{1}{2}(\del_x [\del_y f_1, f_1] - \del_y [\del_x f_1, f_1])_1
      + \frac{i}{2} \del_y(f_1 \del_y f_1 + (\del_y f_1) f_1)_1
               + \frac{i}{2} \del_x(f_1 \del_x f_1 + (\del_x f_1) f_1)_1 \} \\

  + (\lfloor dSdT + f(S, T) \lfloor dxdy)Q(\Xi, f_1)_1 \\
  = 0,
\end{array}
$$
 here $\lfloor dSdT$, $\lfloor dxdy$ simply means getting rid of 
 these symbols and $(\cdots)_1$ means taking $\mathcal{K}^{\perp}$ part.
Note that the term
 $(\del_S \widetilde{\Psi_{\nu}}
                 - \del_T \widetilde{\Phi_{\nu}})$
 disappears since it belongs to $\mathcal{K}$. 

We want to estimate $f_{0, 0} (S, T) = < f_1, f_1 >$.
To do this, we will soon study the inequality of the form
$$
-(\del_S^2 + \del_T^2) <f_1, f_1>
   \leq \cdots,
$$
 with the right hand side calculated by the above equation.
About the terms appearing in this expression,
 we want to make several remarks.

First, the term
 $*_{T_F} d_{\widetilde{A_{\nu}}} *_{T_F} d_{\widetilde{A_{\nu}}} f_1$
 satisfies the estimate
$$
<*_{T_F} d_{\widetilde{A_{\nu}}} *_{T_F} d_{\widetilde{A_{\nu}}} f_1, f_1>
   \geq \eta^2 <f_1, f_1>
$$
 because $A_{\nu} \in Rep_{T^2}(SO(3)) \setminus \phi^{-1}(U_{\eta})$.

The terms containing 
 $f_1$ (without differential) as a factor can be absorbed into
 $\eta^2 f_1$ of the right hand side of the above inequality
 (precisely speaking
 by slightly changing $\eta^2$, but we neglect it),
 because $\eta$ is much larger than the coefficients
 of these terms. 
So the essential terms are
 $[\del_S f_1, \del_T f_1]$, $[\widetilde{\Phi_{\nu}}, \del_S f_1]
                    + [\widetilde{\Psi_{\nu}}, \del_T f_1]$,
 $(\del_T f_1)^2$, $(\del_S f_1)^2$, $(\del_x f_1)^2$, $(\del_y f_1)^2$,
 $[\widetilde{\Phi_{\nu}}, - i \del_T f_1]
   + [i \del_S f_1, \widetilde{\Psi_{\nu}}]$
 and
 $[\del_x f_1, \del_y f_1]$.
In particular, the above inequality can be written in the following form.
$$
\begin{array}{lll}
-(\del_S^2 + \del_T^2) f_{0, 0}
   & \leq - 2f(S, T)\eta^2 f_{0, 0} - 2<\del_S f_1, \del_S f_1> - 2<\del_T f_1, \del_T f_1>
       & \\
   & \;\;\;\;  
                + 2|<[\del_S f_1, \del_T f_1], f_1>| 
    + 2|<[\widetilde{\Phi_{\nu}}, \del_S f_1]
                    + 2[\widetilde{\Psi_{\nu}}, \del_T f_1], f_1>| & \\
   & \;\;\;\; + 4|<(\del_S f_1)^2, f_1>| + 2|<(\del_T f_1)^2, f_1>| \\
   & \;\;\;\;
        + 2|<[\widetilde{\Phi_{\nu}}, - i \del_T f_1] 
        + [i \del_S f_1, \widetilde{\Psi_{\nu}}], f_1>| \\
   & \;\;\;\; + f(S, T) \{4|<[\del_x f_1, \del_y f_1], f_1>|
              + 2|<(\del_x f_1)^2, f_1>| & \\
   & \;\;\;\; + 2|<(\del_y f_1)^2, f_1>| \}.
\end{array}
$$
\begin{prop}
We have an estimate
$f_{0, 0}(S,T) < C
    \exp(- \frac{dist(\del\frac{1}{\ep_{\nu}}B_{\iota}(x), (S, T))}{C'})$
 for some positive constants $C$, $C'$ which do not depend
 on $x$.
\end{prop}
\proof
We estimate the terms in the above inequality.
First, by Schwarz's inequality and the estimate
 $\pll f_1 \pll_{C^r(\frac{1}{\ep_{\nu}}B_{\iota}(x) \times T^2)}
     < C \frac{\ep_{0}^{\lam}}{\eta}$,
$$
\begin{array}{lll}
|<[\del_S f_1, \del_T f_1], f_1>| &
   \leq C \frac{\ep_0^{\lam}}{\eta}
       <\del_S f_1, \del_S f_1>^{\frac{1}{2}}
 <f_1, f_1>^{\frac{1}{2}} & \\
   & < C \frac{\ep_0^{\lam}}{ \eta} (<\del_S f_1, \del_S f_1>
             + <f_1, f_1>).
\end{array}
$$
The term $C \frac{\ep_0^{\lam}}{ \eta} <\del_S f_1, \del_S f_1>$
 is absorbed in $- 2<\del_S f_1, \del_S f_1>$
 and $C \frac{\ep_0^{\lam}}{ \eta}<f_1, f_1>$ is absorbed in 
 $- 2f(S, T)\eta^2 f_{0, 0}$.
The same estimates hold for terms $|<(\del_S f_1)^2, f_1>|$, 
  $|<(\del_T f_1)^2, f_1>|$.

The terms $<[\widetilde{\Phi_{\nu}}, \del_S f_1]
                    + [\widetilde{\Psi_{\nu}}, \del_T f_1], f_1>|$
 and $|<[\widetilde{\Phi_{\nu}}, - i \del_T f_1] 
        + [i \del_S f_1, \widetilde{\Psi_{\nu}}], f_1>|$
 are estimated in the same manner.

Finally, 
$$
\begin{array}{lll}
|<[\del_x f_1, \del_y f_1], f_1>|
  & \leq 2 <|\nabla^{T^2} f_1|^2, |\nabla^{T^2} f_1|^2>^{\frac{1}{2}}
           <f_1, f_1>^{\frac{1}{2}} & \\
  & \leq (C \frac{\ep_0^{\lam}}{ \eta}) <f_1, f_1>,
\end{array}
$$
 by lemma 6.17, here $\nabla^{T^2}$ means the covariant differential in the
 direction of fibers.
The terms $|<(\del_x f_1)^2, f_1>|$ and
   $|<(\del_y f_1)^2, f_1>|$ satisfy similar estimates.
Consequently, we have the following inequality.
$$
\begin{array}{lll}
-(\del_S^2 + \del_T^2) f_{0, 0}
   & \leq - 2f(S, T)(\eta^2 - \ep') f_{0, 0}
         - (2 - \ep')<\del_S f_1, \del_S f_1>
         - (2 - \ep')<\del_T f_1, \del_T f_1>
       & \\
   & \leq - 2f(S, T)(\eta^2 - \ep') f_{0, 0},
\end{array}
$$
 here $\ep'$ is a small number.
$f(S, T)$ satisfies
$
f(S, T) > c_0 > 0
$
 on $\frac{1}{\ep_{\nu}} B_{\iota}(x)$.
So we have
$$
-(\del_S^2 + \del_T^2) f_{0, 0} \leq - \eta' f_{0, 0},
$$
 for some positive constant $\eta'$, which does not depend on $x$.
From this inequality, it is a standard fact that
 the estimate of the proposition follows. \qed
\begin{cor}
$\pll \del_x^l \del_y^m f_1(S, T) \pll^2_{L^2(T^2_F)}
  \leq 
    C \exp(- \frac{dist(\del\frac{1}{\ep_{\nu}}B_{\iota}(x), (S, T))}{C'})$.
\end{cor}
\proof
This follows from partial integration.
$$
\begin{array}{lll}
<\del_x^l \del_y^m f_1, \del_x^l \del_y^m f_1>
  & \leq |<\del_x^{2l} \del_y^{2m} f_1, f_1>| & \\
  & < C \ep_0^{\lam}/\eta <f_1, f_1>^{\frac{1}{2}}.
\end{array}  
$$ 
The estimate follows from proposition 6.18. \qed\\

Next we want to estimate $f_{a, b}(S, T)$. 
Namely, we prove the following.
\begin{prop}
$f_{a, b}(S,T) < C
    \exp(- \frac{dist(\del\frac{1}{\ep_{\nu}}B_{\iota}(x), (S, T))}{C'}).$ \\
Moreover,
$$
\pll \del_x^l \del_y^m e_{a, b} (S, T) \pll_{L^2(T_F^2)}^2
  \leq C \exp(- \frac{dist(\del\frac{1}{\ep_{\nu}}B_{\iota}(x), (S, T))}{C'})
$$
 for any $a, b \in \Z_{\geq 0}$ holds.
\end{prop}
\proof
We prove it by induction about $a+b$.
We have proved the case $a = b = 0$.

Now suppose the proposition is shown in the cases $a+b < k$.
We can proceed in the same manner as in the case of $a = b = 0$.
Namely, differentiating the $(E_1)$ equation
 by $\del_S^a \del_T^b$, $a+b = k$, gives an expression for
 $(\del_S^2 + \del_T^2)(\del_S^a \del_T^b f_1)$.
By induction hypothesis,
%
 it is easy to show that
\begin{eqnarray*}
\Delta_{S, T}f_{a, b} (S, T)
    &  = & 2<\Delta_{S, T} e_{a, b}, e_{a, b}>
       - 2<\frac{\del e_{a, b}}{\del S}, \frac{\del e_{a, b}}{\del S}> 
       - 2<\frac{\del e_{a, b}}{\del T}, \frac{\del e_{a, b}}{\del T}> \\
    &  < & 2<(f(S, T)*_{T_F} d_{\widetilde{A_{\nu}}}
                          *_{T_F} d_{\widetilde{A_{\nu}}} - \ep') e_{a, b}, 
           e_{a, b}> \\ 
    &  < & -\eta' f_{a, b}.
\end{eqnarray*}
 by the same calculation as in the proof of proposition 6.18.
Again, by a standard argument
 the first estimate of the proposition follows.
The second estimate follows in the same way as corollary 6.19. \qed \\


\begin{cor}
$\pll f_1 \pll_{C^r(B_1(S, T) \times T^2)}
   < C \exp(- \frac{dist(\del\frac{1}{\ep_{\nu}}B_{\iota}(x), (S, T))}{C'})$.
\end{cor}
\proof
This follows from the Sobolev's inequalities and proposition 6.20. \qed \\

\begin{cor}
$\pll \Phi_{1, \nu, G}(S, T) \pll_{C^r(B_1(S, T) \times T^2)},
 \pll \Psi_{1, \nu, G}(S, T) \pll_{C^r(B_1(S, T) \times T^2)}
    < C \exp(- \frac{dist(\del\frac{1}{\ep_{\nu}}B_{\iota}(x), (S, T))}{C'})$.
\end{cor}
Proof:
This follows directly from the above corollary
 and lemma 6.13,
 because $\widetilde{\Xi_{\nu}}$ has no $(\mathcal{K})^{\perp}$
 part and so this part of $\Xi_{\nu}$ has the form 
 $Y(\widetilde{\Xi_{\nu, G}}, f_1)$ with the estimate 
$$
\pll Y(\widetilde{\Xi_{\nu}}, f_1) \pll_{C^r(B_1(S, T) \times T^2)}
          < C \pll f_1 \pll_{C^{r+1}(B_1(S, T) \times T^2)}. 
$$\qed
\bc
$\pll A_{\nu, G} - \widetilde{A_{\nu}} \pll_{C^r(B_1(S, T) \times T^2)}
  < C \exp(- \frac{dist(\del\frac{1}{\ep_{\nu}}B_{\iota}(x), (S, T))}{C'})$.
\ec
\proof
This is clear from the definition of $\Xi_{\nu, G}$
 and corollary 6.21. \qed\\

Using these results, we can estimate the ($E_2$) part of the equation.
First we introduce a moderate gauge for $\Phi_{2, \nu, G}$
 and $\Psi_{2, \nu, G}$.
\begin{lem}
Let $\Omega = A(S, T) + \Phi(S, T)dS + \Psi(S, T)dT$
 be any connection on our bundle.
Let $\Phi_2$, $\Psi_2$ be the $\mathcal{K}$ parts of 
 $\Phi$, $\Psi$.
Then there is a unitary gauge transformation $j$ of the form $j = \exp(b)$, 
 $b$ is a $\mathcal {K}$ valued section over
 $\frac{1}{\ep_{\nu}} B_{\iota}(x)$,
 such that the $\mathcal{K}$ part of the base component of the
 transformed connection satisfy
$$
\del_S \Psi_2 + \del_T \Phi_2 = 0
$$
 on $\frac{1}{\ep_{\nu}} B_{\iota}(x)$
 and 
 orthogonal to harmonic one forms.
\end{lem}
Proof:
If we apply the transformation $j$ to $\Omega$, the base component becomes
$$
\frac{\del b}{\del S} dS + \frac{\del b}{\del T} dT
 + \exp (-b) (\Phi(S, T)dS + \Psi(S, T)dT) \exp (b).
$$
The $\mathcal{K}$ part of it is simply
$$
\frac{\del b}{\del S} dS + \frac{\del b}{\del T} dT
 + \Phi_2(S, T)dS + \Psi_2(S, T)dT.
$$
Now the space $L^2(\Omega^1(\frac{1}{\ep_{\nu}} B_{\iota}(x)))$
 of $L^2$ one forms on a two dimensional disc $\frac{1}{\ep_{\nu}} B_{\iota}(x)$
 decomposes as
$$
L^2(\Omega^1(\frac{1}{\ep_{\nu}} B_{\iota}(x)))
 = d L^2_1(\Omega_0^0(\frac{1}{\ep_{\nu}} B_{\iota}(x)))
 \oplus d^* L^2_1(\Omega^2_0(\frac{1}{\ep_{\nu}} B_{\iota}(x)))
 \oplus \mathcal{H},
$$ 
 where $L^2_1(\Omega^0_0(\frac{1}{\ep_{\nu}} B_{\iota}(x)))$ is the space of 
 $L^2_1$ completion of the space of the smooth real zero forms
 with compact support on the open disk,
 $L^2_1(\Omega^2_0(\frac{1}{\ep_{\nu}} B_{\iota}(x)))$ is the space of 
 $L^2_1$ completion of the smooth real two forms
 with compact support on the disk and
 $\mathcal{H}$ is the space of $L^2$ harmonic one forms. 

On the other hand, 
$$
d L^2_1(\Omega^0_0(\frac{1}{\ep_{\nu}} B_{\iota}(x)))
 \oplus \mathcal{H} = d L^2_1(\Omega^0(\frac{1}{\ep_{\nu}} B_{\iota}(x)))
$$
 holds.

So, we can take the gauge for $\Phi_2$, $\Psi_2$
 simply by subtracting the
 $ d L^2_1(\Omega^0_0(\frac{1}{\ep_{\nu}} B_{\iota}(x))) \oplus \mathcal{H}$
 part of $\Phi_2 dS + \Psi_2 dT$. \qed \\

Note that we can do this by a
 $C^{\infty}(\frac{1}{\ep_{\nu}} B_{\iota}(x))$-gauge
 transformation.
In particular, it does not disturb the estimates in 
 proposition 6.20 and corollary 6.21.
We assume we have taken $\Phi_{2, \nu, G}$ and $\Psi_{2, \nu, G}$
 in this gauge.

\begin{prop}
$\pll \Phi_{2, \nu, G}(S, T) \pll_{C^r(B_1(S, T) \times T^2)},
 \pll \Psi_{2, \nu, G}(S, T) \pll_{C^r(B_1(S, T) \times T^2)}
     < C \frac{\iota}{\ep_{\nu}}
     \exp(- 2\frac{dist(\del\frac{1}{\ep_{\nu}}B_{\iota}(x), (S, T))}{C'})$.
\end{prop}
Proof:
Note that $f_1$ and
 $[\mathcal{K}, (\mathcal{K})^{\perp}]$ have no $\mathcal{K}$-component.
From this remark, ($E_2$) part of the equations has the form 
$$
\begin{array}{l}
(\del_S \widetilde{\Psi_{\nu}}
                 - \del_T \widetilde{\Phi_{\nu}})
 + \frac{1}{2}(\del_S [\del_T f_1, f_1] - \del_T [\del_S f_1, f_1])_2 \\
 + \frac{i}{2} \del_T (f_1 \del_T f_1 + (\del_T f_1) f_1)_2
 + \frac{i}{2} \del_S (f_1 \del_S f_1 + (\del_S f_1) f_1)_2 
 + [\del_T f_1, \del_S f_1]_2 \\
 + f(S, T) \{[\del_y f_1, \del_x f_1]_2 
 + \frac{1}{2}(\del_x [\del_y f_1, f_1] - \del_y [\del_x f_1, f_1])_2 \\
 + \frac{i}{2} \del_y (f_1 \del_y f_1 + (\del_y f_1) f_1)_2
 + \frac{i}{2} \del_x (f_1 \del_x f_1 + (\del_x f_1) f_1)_2 \} \\  
  + (\lfloor dSdT + f(S, T) \lfloor dxdy) Q(\Xi, f_1)_2 \\
 = 0,
\end{array}
$$
 here $( \cdots)_2$ means taking $\mathcal{K}$ part.
All the terms other than $\del_S \widetilde{\Psi_{\nu}}
                 - \del_T \widetilde{\Phi_{\nu}}$
 contain $f_1$ or its derivative as a factor.
So in terms of $\Psi_{2, \nu, G}$ and $\Phi_{2, \nu, G}$,
 we can write it as
$$
\frac{\del}{\del S} \Psi_{2, \nu, G} - \frac{\del}{\del T} \Phi_{2, \nu, G}
 = P(\widetilde{\Xi_{\nu}}, f_1)
$$
 where $P$ satisfies the estimate
$$
\begin{array}{l}
\pll P \pll_{C^r(B_1(S, T) \times T^2)}
 < C  \exp(- 2\frac{dist(\del\frac{1}{\ep_{\nu}}B_{\iota}(x), (S, T))}{C'}).
\end{array}
$$
On $d^* L^2_1(\Omega^2_0(\frac{1}{\ep_{\nu}} B_{\iota}(x)))$,
 the exterior differential satisfies
$$
\pll d \eta \pll_{C^r(B_1(S, T) \times T^2)}
    > C \frac{\ep_{\nu}}{\iota} \pll \eta \pll_{C^r(B_1(S, T) \times T^2)}. 
$$
By our choice of the gauge of $\Phi_{2, \nu, G}dS + \Psi_{2, \nu, G}dT$, 
 we can apply this inequality. 
Namely, 
$$
\pll \Phi_{2, \nu, G} \pll_{C^r(B_1(S, T) \times T^2)},
 \pll \Psi_{2, \nu, G} \pll_{C^r(B_1(S, T) \times T^2)}
 < C \frac{\iota}{\ep_{\nu}} 
   \exp(- 2\frac{dist(\del\frac{1}{\ep_{\nu}}B_{\iota}(x), (S, T))}{C'}).
$$
This proves the proposition. \qed\\
%
%

Theorem 6.6 (and so theorem 3.1) now easily follows.
Recall $S_{\eta} = S_1 \cup S_2 \cup S_{3, \eta}$.
We have shown that at any point 
 $x \in T^2_B \setminus S_{\eta}$,
 there is a neighbourhood $x \in U \subset T_B^2$
 and a gauge transformation $g_{U, \nu}$
 such that for any $r$, the $C^r$ norm of the connection matrix of
 $g_{U, \nu}^* \Xi_{\ep_{\nu}}$
 is uniformly bounded with respect to $\nu$ in a smaller neighbourhood 
 $\pi^{-1}(U') \subset \pi^{-1}(U)$.
So on $\pi^{-1}(U')$, there is a subsequence of $\Xi_{\ep_{\nu}}$,
 which converges
 in $C^{\infty}$ sense to a limit unitary connection $\Xi$.
Since the local convergence of connections modulo gauge transformation is
 equivalent to global convergence
 (\cite{DK}, corollary 4.4.8), theorem 6.6 follows. \qed\\
\br
We showed the theorem under the assumption $(\bold{B})$.
If we can solve the problem of remark 5.14, this assumption is not necessary.
\er
We remark about 
 the relationship between the gauge theory side and
 the minimal surface side, especially 
 the relationship between the divergence loci of the curvature 
 (with respect to the metric on $\widehat{M_1}$) and
 the behavior of associated holomorphic sections
 around them.
First, we note the following.
\begin{prop}
At the loci $\pi^{-1}(S_1)$, the curvature diverges
 (with respect to the metric on $\widehat{M_1}$).
\end{prop}
\proof
If the curvature does not diverge, the $L^2$ norm of the curvature
 on a unit ball in the rescaled space will be of order $\ep_{\nu}$.
By gluing standard Uhlenbeck's Coulomb gauge for ASD connections,
 it is easy
 to construct a gauge on $\frac{1}{\ep_{\nu}} B(\iota) \times U$,
 here $U$ is an open ball on the fiber torus, in which the 
 $C^r$ norm of the connection matrix is of order $\ep_{\nu}$.
In particular, it is finite in the non-rescaled space.
This contradicts the assumption that the section $\phi_{\nu}$
 develops bubbles on $S_1$. \qed\\

\bpr
Moreover, if $x \in S_1 \setminus S_2$, the point measure 
 $\mu(x)$ is positive.
\epr
\proof 
If $\mu(x) = 0$, then for any small positive $\ep$, 
 there is a neighbourhood $U$ of $x$ and $\nu_0$ such that
 for all $\nu > \nu_0$, $\mu_{\nu} (U) < \ep$ holds.
Then applying the discussion of this section to the rescaled open set
 of $\pi^{-1} (U)$, the convergence of the connections 
 with respect to the metric of 
 $\widehat M_1$ is proved.
This contradicts to the assumption that $x \in S_1$. \qed \\

By definition, it is clear that if $x \in S_{3} = \cup_i S_{3, \eta_i}$
 (or, strictly speaking, $x \in S_3 \setminus S_2$),
 the curvature diverges on $\pi^{-1}(x)$.
Conversely, we have shown that if $x \notin S_1 \cup S_2 \cup S_{3}$
 the connection converges.
In other words, the divergence of the curvature and the positiveness
 of the point measure with respect to $\mu$ are the same, modulo $S_2$ 
 (these happen on the loci $S_1 \cup S_3$).
But we cannot yet relate the divergence of the curvature
 and the bubbling of the corresponding holomorphic sections.

However, it may well be the case that 
 the divergence of the curvature and the bubbles of the holomorphic
 sections are the same and all the bubbles are contained
 in $S_1$, so $S = S_1 \cup S_2$ is a finite set.
To prove it, we need several things.
For example, we have to study the sections $\phi_{\nu}$
 around the first type bubbles and
 also we should have the energy bound of doubly periodic instantons by the
 constant $8 \pi^2$ (and further have to prove the non-constantness
 of the map $\rho$ associated to these instantons),
 and of course further study of the third type bubbles
 is needed.

Another problem is the analysis at the 
 most degenerating loci $S_2$.
To consider the reduced version of \cite{DS}, 
 this type of analysis would be inevitable.\\

As a final remark, we argue about the deduction of
 parts of the results of the paper \cite{DS} of Dostoglou-Salamon,
 as in the last section,
 and also reprove the (stronger) result of Chen \cite{C}.
Let $E$ be an $SO(3)$ bundle over $\Sigma_h \times (\Sigma_g)_{\ep_{\nu}}$
 which is non-trivial when restricted to the fiber $(\Sigma_g)_{\ep_{\nu}}$,
 here $(\Sigma_g)_{\ep_{\nu}}$ means a Riemannian surface
 with the metric which is the $\ep_{\nu}$-rescaling of a fixed metric on 
 $\Sigma_g$, $g \geq 2$.
Let $\Xi_{\ep_{\nu}}$ be a sequence of ASD connections on $E$
 with respect to the Riemannian metrics of 
 $\Sigma_h \times (\Sigma_g)_{\ep_{\nu}}$.
By taking a subsequence if necessary, we can define $S_1$ as the points of
 the base which support the type one or type two bubbles.

By the characterization of the 
 third type bubbles, the curvature $F_{A_{\ep_{\nu}}}$ of the fiber part of the 
 connection is small on $\Sigma_h \setminus U_{\nu}(S_1)$,
 here $ U_{\nu}(S_1)$ is a small neibourhood of
 $S_1$ which sufficiently slowly shrinks to $S_1$ as $\nu \to \infty$.
Then by lemma 2.14 of
 \cite{F1}, the fiber part of the connections
 on $\pi^{-1}(x)$, $x \in \Sigma_h \setminus U_{\nu}(S_1)$,
 are complex gauge equivalent to flat connections.
So we have a sequence of maps
$$
\phi_{\nu}: \Sigma_h \setminus U_{\nu}(S_1) \to \mathcal{F},
$$
 here $\mathcal{F}$ is the space of flat connections on the restriction of $E$ 
 to the fiber $\Sigma_g$.

The energies of these maps are again  bounded
 (in this case, since the domains of the maps
 are not closed, we cannot bound the energy simply by
 topological quantity determined by the bundle, as we have done
 in the main text.
 However, by the method of the proof of proposition 5.4
 (see also proposition 5.15), since we do not have singularities in the 
 moduli here, we can bound the energy of holomorphic maps
 by the Yang-Mills energy of the connections),
 so there is a subsequence of $\phi_{\nu}$
 which converges modulo finite bubbles to a limit holomorphic map
$$
\phi: \Sigma_h \setminus S_1 \to \mathcal{F},
$$
 which can be extended to the whole $\Sigma_h$.
Let $S_2$ be the bubbling points of the sequence.

We construct the limit measure $\mu$ as in the main text and set
$$
S_3 = \{ x \in \Sigma_h \big| \mu(x) \geq \ep_0 \}
$$
 for sufficiently small $\ep_0$.
We define the finite subset $S = S_1 \cup S_2 \cup S_3$
 as in the proof of the theorem 6.6.
Here $U_{\eta}$ appeared in the main text is
 absent since the moduli space is smooth
 and transversal in this case.
\begin{prop}
There is a subsequence of $\Xi_{\ep_{\nu}}$ which 
 converges, modulo gauge transformation, to a limit unitary connection
 over $\pi^{-1}(\Sigma_h \setminus S)$ with respect to the metric on
 $\Sigma_h \times \Sigma_g$ (non rescaled).
\end{prop}
\proof
This follows essentially from the argument of Fukaya \cite{F1}, proof
 of the main lemma 4.42.
Let $y \in \Sigma_h \setminus S$
 and $U_y$ be a small neighbourhood of $y$ in $\Sigma_h \setminus S$.
We rescale $\pi^{-1}(U_y)$
 by the rate $\frac{1}{\ep_{\nu}}$ and define a sequence of connections
 $\{ \Xi_{\nu} \}$ as before. 
That is, 
$$
A_{\nu}(x) = A_{\ep_{\nu}}(x_0 + \ep_{\nu} x),
$$
$$
\Phi_{\nu}(x) = \ep_{\nu} \Phi_{\ep_{\nu}}(x_0 + \ep_{\nu} x), \;\;
\Psi_{\nu}(x) = \ep_{\nu} \Psi_{\ep_{\nu}}(x_0 + \ep_{\nu} x),
$$
 $x_0$ is a fixed point on $\pi^{-1}(U_{y})$.
This is a sequence of ASD connections.

As in the proof of proposition 5.15, we can assume that
 there is a slice by flat connections
 in $\mathcal{A}$ which represent
 $[A_{\nu}]$ (the complex gauge equivalence classes
 of the fiber components over $\frac{1}{\ep_{\nu}}U_y$).
We denote it as $Q_{\nu}$.
Since $[A_{\nu}]$
 is the rescaling of a sequence of holomorphic maps to the 
 moduli of flat connections which converges on $U_y$,
 we can assume
$$
\pll \frac{\del^m Q_{\nu}}{\del S^l \del T^{m-l}} \pll_{C^r(\Sigma_g)}
 < C \ep_{\nu}^m.
$$
Since $[A_{\nu}]$ is a holomorphic map, there are unique
 $\widetilde \Phi_{\nu}$ and $\widetilde \Psi_{\nu}$ which satisfy
$$
\frac{\del Q_{\nu}}{\del S} + *_{\Sigma_g} \frac{\del Q_{\nu}}{\del T}
 = d_{Q_{\nu}} \widetilde \Psi_{\nu}
 + *_{\Sigma_g} d_{Q_{\nu}} \widetilde \Phi_{\nu}.
$$

Since in this case the covariant differential
 in the direction of the fiber $d_{Q_{\nu}}$ has no kernel,
 $\widetilde \Phi_{\nu}$ and $\widetilde \Psi_{\nu}$ satisfy
$$
\pll \frac{\del^m \widetilde \Phi_{\nu}}{\del S^l \del T^{m-l}}
 \pll_{C^r(\Sigma_g)}
 < C \ep^{m+1}_{\nu},
$$
$$
\pll \frac{\del^m \widetilde \Psi_{\nu}}{\del S^l \del T^{m-l}}
 \pll_{C^r(\Sigma_g)}
 < C \ep^{m+1}_{\nu}.
$$

Here again the point is that since the connection in the fiber direction
 does not have isotropy group, fixing the connection on all the
 fibers fixes the connection of $E$ (see the proof of proposition 5.15). 
We unitary gauge transform the ASD connection $\Xi_{\nu}$
 so that it is obtained from
 $Q_{\nu} + \widetilde \Phi_{\nu} dS + \widetilde \Psi_{\nu} dt$
 by a hermitian gauge transformation
 $g_{\nu}$.
We denote this unitary transformed connection also by $\Xi_{\nu}$
 to save letters.
 
Now take a disc $D$ in $\frac{1}{\ep_{\nu}}U_x$
 of radius one.
By theorem 5.1, there is a flat connection $\Xi_D$
 over $E|_{D \times \Sigma_g}$ and we can unitary gauge transform 
 $\Xi_{\nu}$ by $h$ so that
$$
\pll h^*\Xi_{\nu} - \Xi_D \pll_{C^r(D \times \Sigma_g)} < \ep_0^{\lam}.
$$
Moreover, we can assume that the fiber component
 $A_D$ of $\Xi_D$ is close to $Q_{\nu}|_D$.
In particular, the fiber components of
 $(g_{\nu}^{-1})^* \Xi_{\nu}$ and $h^* \Xi_{\nu}$
 are close.
Since $h$ is unitary and $g$ is hermitian,
 both transformations should be
 of small norm. 

It follows that 
$$
\pll g_{\nu} - 1 \pll_{C^r(D \times \Sigma_g)} < C \ep^{\lam}_0
$$
 holds.
Note this estimate is global, that is, it is valid over
 $\pi^{-1}(\frac{1}{\ep_{\nu}}U_x)$:
$$
\pll g_{\nu} - 1 \pll_{C^r(\pi^{-1}(\frac{1}{\ep_{\nu}} U_x))} < C \ep^{\lam}_0.
$$
 This estimate corresponds to Fukaya \cite{F1}, lemma 6.1
 (this claim can, in fact, be proved
 without the assumption that $\Xi_{\nu}$ is ASD).

What we want is the stronger estimate
$$
\pll \frac{\del^m (g_{\nu} - 1)}{\del S^l \del T^{m-l}}
 \pll_{C^r(\Sigma_g)} < C \ep_{\nu}^{m+1}.
$$
Let us write $g_{\nu}$ as 
 $\exp (e_{\nu})$, $e_{\nu}$ is a hermitian matrix valued section.
It satisfies the equation
$$
\Delta_{S, T} e_{\nu} + \Delta_{A_{\nu}} e_{\nu}
 = Q(e; S, T),
$$
 here $Q$ is the sum of terms alike those of $F^{HV}$
 appeared before lemma 6.17.
However, here we cannot delete the term
 $\del_S \widetilde \Psi_{\nu} - \del_T \widetilde \Phi_{\nu}$.
By this reason, the differential inequality we obtain
 is different from the main text (see proposition 6.18).
Namely, it becomes the one which appeared in \cite{F1}, lemma 6.46:
$$
\Delta_{S, T} f(S, T) \leq - \tau f(S, T) +
 C \ep_{\nu}^2 \sqrt{f(S, T)},
$$
 here $\tau$ is a positive constant
 and $f = \pll e_{\nu}(S, T) \pll^2_{L^2(\Sigma_g)}$
 (note this equation is a little different from the one actually
 appeared in \cite{F1}, lemma 6.46. 
Namely, the third term of the right hand side there is dropped.
This is safely done by the same calculation as in
 the proof of proposition 6.18.
The reason to do this is that there is an error in the estimate
 of \cite{F1}, lemma 6.10.
This point can also be fixed by the argument as
 in the proof of proposition 6.18).    
Then we can apply lemma 6.46 of \cite{F1} and
 deduce the estimate
$$
f(S, T) \leq C \ep_{\nu}^4
 + C \exp(- \frac{dist(\del \frac{1}{\ep_{\nu}}U_x, (S, T))}{C}) \sup f.
$$
The estimates of 
 $f_{a, b} = \pll \del_S^a \del_T^b e_{\nu}(S, T) \pll^2_{L^2(\Sigma_g)}$
 and of terms including differentials in the fiber direction
 can be done in the same way (see \cite{F1}, pages 526 and 527
 for full details).\qed\\


In the non-reducible case, we can make the 
 relationship between
 the bubbling of the sequence of the holomorphic maps
 $\phi_{\nu}$ and the bubbling of the connections
 clearer than in the reducible case.
Namely, we can prove the following.
\bpr
$S_3 = S_1 \cup S_2$ holds.
Moreover, the loci $S_2$ is precisely the same as the loci
 at which the sequence of HYM connections
 $\Xi_{\ep_{\nu}}$ develops the third type bubbles. 
\epr 
\proof
In this case, we have, as noted before, the relation between 
 Yang-Mills energy of the
 connections and the energy of the holomorphic maps defined
 on a small disc $U_x$ on $\Sigma_h \setminus U_{\nu}(S_1)$ 
 as follows.
$$
E(\phi_{\ep_{\nu}}|_{U_x}) - C (YM_{\ep_{\nu}}(\Xi_{\ep_{\nu}}|_{U_x}))
  < YM_{\ep_{\nu}}(\Xi_{\ep_{\nu}}|_{U_x})
  < E(\phi_{\ep_{\nu}}|_{U_x}) + C (YM_{\ep_{\nu}}(\Xi_{\ep_{\nu}}|_{U_x})),
$$
 for some positive constant $C$, which does not depend on $U_x$
 and $YM_{\ep_{\nu}}$ is the Yang-Mills energy defined using the 
 metric on $\Sigma_h \times (\Sigma_g)_{\ep_{\nu}}$.
So, if $x \in S_2$, that is, $E(\phi_{\ep_{\nu}}|_{U_x}) > C$
 for any small $U_x$ and large $\nu$, 
 we have 
 $YM_{\ep_{\nu}}(\Xi_{\ep_{\nu}}|_{U_x}) > C$.
In particular, the curvature blows up at $x$.
Since the first and the second type bubbles are already 
 contained in the loci $S_1$, 
 this must be the third type bubble.
This also shows $S_1 \cup S_2 \subset S_3$.
 
Conversely, let $x \in S_3 \setminus S_1 \cup S_2$.
Then the situation is similar to the one in the last proposition.
Around this point, the rescaled ASD connection $\Xi_{\nu}$
 can be unitary gauge transformed
 so that it is obtained by a hermitian gauge transformation
 from the connection $\Xi_{0, \nu}$ 
 whose fiber component is in the form $Q_{\nu}$
 in the proof of the previous proposition
 and 
 complex gauge equivalent to those of $\Xi_{\nu}$.
And the whole $\Xi_{0, \nu}$ is constructed from this component
 as before.
The connection $\Xi_{0, \nu} = Q_{\nu} + \widetilde{\Phi_{0, \nu}} dS
 + \widetilde{\Psi_{0, \nu}} dT$ satisfies the estimates
$$
\pll \nabla^m_B Q_{\nu} \pll_{C^r(\Sigma_g)} < C\ep_{\nu}^m 
$$  
 and
$$
\pll \nabla^m_B \widetilde{\Phi_{0, \nu}} \pll_{C^r(\Sigma_g)},
  \pll \nabla^m_B \widetilde{\Psi_{0, \nu}} \pll_{C^r(\Sigma_g)}
    < C\ep_{\nu}^{m+1}, 
$$
 where $\nabla_B$ is the covariant differential in the direction of the base.

The difference from the situation of proposition 6.29 is that there 
 is no bound, by a small constant, of the energy of the connection,
 by definition of $S_3$.
However, this does not really cause a big problem here
 by the fact that $\Xi_{\nu}$ is ASD (see the proof of the previous
 proposition). 
Let us write by $g_{\nu} = \exp(e_{\nu})$ the hermitian
 gauge transformation which maps  $\Xi_{0, \nu}$ to $\Xi_{\nu}$.
Because $x \notin S_1 \cup S_2$, the $C^0$ norm of
 $e_{\nu}$ is bounded by a small constant, say $\ep$.
Moreover, because $\Xi_{\nu}$ is ASD, the $C^r$ norms of $e_{\nu}$
 are also bounded by some constants.
These two conditions imply that the $C^r$ norms of $e_{\nu}$ are also
 bounded in terms of $\ep$.
Then the argument of the proof of proposition 6.29 
 applies
 and so $\Xi_{\ep_{\nu}}$ converges around $x$ in $C^{\infty}$ sense.
This contradicts the assumption that $x \in S_3$
 and so $S_3 = S_1 \cup S_2$ \qed

\begin{note}
This corresponds to the statement required in section 3 of
 the erratum of \cite{DS}, namely, the nontrivialness
 of the third type bubbles as holomorphic maps.
Note also that
 in the erratum they used the divergence of the function 
 $\ep^{-1} \pll F_{A_{\ep_{\nu}}} \pll_{L^{\infty}}
  + \pll \del_s A_{\ep_{\nu}} - d_{A_{\ep_{\nu}}} \Phi_{\ep_{\nu}} \pll_{L^{\infty}}$
 for definition of bubbles.
However, we have shown that we can instead use the function
 $\ep^{-1} \pll F_{A_{\ep_{\nu}}} \pll_{L^{\infty}}^{\frac{1}{2}}
  + \pll \del_s A_{\ep_{\nu}} - d_{A_{\ep_{\nu}}} \Phi_{\ep_{\nu}} \pll_{L^{\infty}}$
 as in the original paper of \cite{DS}.
\end{note}

Now we mention about the relation between our method and
 Chen's result \cite{C}.
Let $E$ be an $SU(N)$ (or $SO(N)$) bundle over the metrically degenerating
 family of product of Riemannian surfaces 
 $\Sigma_h \times (\Sigma_g)_{\ep_{\nu}}$, $g \geq 2$, 
 and consider a smooth component of the moduli space of
 flat connections on the bundle over $\Sigma_g$ obtained by
 restricting $E$.
Let $\Xi_{\ep_{\nu}}$ be a sequence of ASD connections on $E$
 with respect to the degenerating metrics.
Then he proved the $C^0$ convergence of the fiber component
 of the connections to a holomorphic map from the union of $\Sigma_h$
 and finitely many spheres to the moduli space of flat connections on
 $\Sigma_g$ associated to the restriction of the given bundle,
 modulo gauge fixing and taking subsequences
(actually the original statement is about general $SU(N)$ bundles, 
 however as noted in the introduction, Chen seems to assume
 the moduli is smooth.
For example, the estimate (4.3) in \cite{C} does not seem to hold
 in the reducible cases).

Our version of the theorem is the following.
The sets $S_1$, $S_2$ and $S = S_1 \cup S_2$ are defined in the same way
 as proposition 6.29.
The proof of it is already clear.
\bt
There is a subsequence of connections $\Xi_{\ep_{\nu}}$
 which converges in $C^{\infty}$ sense 
 on $\pi^{-1}(\Sigma_h \setminus S)$.
The fiber part of the limit connection is flat
 and describes a holomorphic map $\phi$ from
 $\Sigma_h \setminus S$ to $\mathcal{F}$,
 which can be extended to a map from the whole
 $\Sigma_h$.
 
Moreover, we can define a sequence of holomorphic maps 
$$
\phi_{\nu} : \Sigma_h \setminus U_{\nu} \to \mathcal{F},
$$
 associated to $\Xi_{\ep_{\nu}}$,
 where $U_{\nu}$ is a neighbourhood of $S_1$ which shrinks to 
 $S_1$ as $\nu \to \infty$.
At the locus $S_2$, this sequence develops sphere bubbles
 and on any compact subset of
 $\Sigma_h \setminus S$, it converges (in $C^{\infty}$ sense)
 to $\phi$.
  \qed
\et
\begin{note}
In our case, around the locus $S_1$, the holomorphic maps are not defined,
 before going to the limit.
On the other hand, Chen \cite{C} proved his theorem by 
 different methods and in particular constructed
 a sequence of maps defined on small discs covering
 $\Sigma_h$ to $\mathcal{F}$, and proved that
 this converges (modulo patching by gauge transformations)
 to a map from the union of $\Sigma_h$
 and finitely many spheres to the moduli space.
In particular, he proved the bubblings at $S_1$ (or at least
 at the first type bubbling locus, the analysis of the second type bubbles
 seems to be absent there).
From this point, it is desirable to prove the bubblings at $S_1$
 in our case too.
\end{note}

\section{Proof of the main theorem}
In this section, we prove the main theorem 2.6.
Here we get back the $i \R$ part of the connection
 and consider HYM (not necessarily ASD) connections (see section 3).
So in this section, $E_{\ep_{\nu}}$ are $U(2)$-bundles as originally defined,
 contrary to convention 6.1.
Let us take $(\widehat{M_{\ep_{\nu}}}, E_{\ep_{\nu}}, \Xi_{\ep_{\nu}})$
 as before.
First we have to treat the $i\R$-part of the connections.
We do not prove the convergence of the whole of the connections, 
 but only treat the parts which are
 relevant to the construction of the Lagrangians.
During this discussion, we again hyperK\"ahler
 rotate the manifolds and take isothermal coordinates ($S, T$)
 on the base of the fibration.

We write the $i \R$-part of the connection $\Xi_{\ep_{\nu}}$ as
 $\Xi'_{\ep_{\nu}} = A'_{\ep_{\nu}} + \Phi'_{\ep_{\nu}}dS
           + \Psi'_{\ep_{\nu}}dT$.
Then the equations which $\Xi'_{\ep_{\nu}}$ satisfies become
 (they are the equations(4) and (6) in section 3)
$$
\begin{array}{l}
(h_1): \del_T \Phi'_{\ep_{\nu}} - \del_S \Psi'_{\ep_{\nu}}
         - f(S, T) \ep_{\nu}^{-2} F_{A'_{\ep_{\nu}}} = 0,\\ 
(h_2): *_{T_F^2} (-(d_F \Phi'_{\ep_{\nu}} - \del_S A'_{\ep_{\nu}}) dT
         + (d_F \Psi'_{\ep_{\nu}} - \del_T A'_{\ep_{\nu}}) dS) \\
   \hs{1in}
       + (d_F \Phi'_{\ep_{\nu}} - \del_S A'_{\ep_{\nu}}) dS
         + (d_F \Psi'_{\ep_{\nu}} - \del_T A'_{\ep_{\nu}}) dT
        = 2c_0 E_2 \omega,
\end{array}
$$
 here $d_F$ is the exterior differential in the direction of the fibers
 and $f(S, T)$ is a positive function as in the previous section.
We divide the equations into two parts:
 one for fiberwise constant part and the other for its perpendicular part
 (in fiberwise $L^2$ sense. This is the same as saying the 
 integrals along the fibers are zero).
We write, as in the previous section, the fiberwise constant part of some object
 as $(\cdots)_2$ and the perpendicular part as  $(\cdots)_1$.
The equations become
$$
\begin{array}{l}
(h_1)_1: \del_T(\Phi'_{\ep_{\nu}})_1 - \del_S(\Psi'_{\ep_{\nu}})_1
           - f(S, T) \ep_{\nu}^{-2} F_{A'_{\ep_{\nu}}}
           = 0,\\
(h_1)_2: \del_T(\Phi'_{\ep_{\nu}})_2
           - \del_S (\Psi'_{\ep_{\nu}})_2 = 0,\\
(h_2)_1: \ast_{T^2_F} (-(d_F \Phi'_{\ep_{\nu}} - \del_S A'_{\ep_{\nu}})_1 dT
             + (d_F \Psi'_{\ep_{\nu}} - \del_T A'_{\ep_{\nu}})_1 dS) \\
          \hs{1in}
       + (d_F \Phi'_{\ep_{\nu}} - \del_S A'_{\ep_{\nu}})_1 dS
         + (d_F \Psi'_{\ep_{\nu}} - \del_T A'_{\ep_{\nu}})_1 dT
        = 0,\\
(h_2)_2: \ast_{T^2_F} (\del_S(A'_{\ep_{\nu}})_2 dT
                - \del_T(A'_{\ep_{\nu}})_2 dS)
    -  \del_S(A'_{\ep_{\nu}})_2 dS
                - \del_T(A'_{\ep_{\nu}})_2 dT
  = 2c_0 E_2 \omega.        
\end{array}
$$
We treat only $(\cdots)_2$-part of the connections, since 
 $(A'_{\ep_{\nu}})_1$-part is complex gauge equivalent to zero
 on each fiber
 and irrelevant to our construction.
Note that the dependence of $(A'_{\ep_{\nu}})_2$
 to the gauge is only by the ambiguity of $\Z^2$, 
 coming from the gauge transformations of the type (see section 4)
$$
g = \left( \begin{array}{cc}
         e^{2\pi i(nx - my)} & 0 \\
                0            & e^{-2\pi i(nx - my)}
          \end{array} \right).
$$
\bpr
There is a finite open covering 
 $\{U_i\}$ of $T^2_B$ and on each $U_i$
 there is a subsequence of  $(A'_{\ep_{\nu}}
 + \Phi'_{\ep_{\nu}} dS + \Psi'_{\ep_{\nu}} dT)_2$
 such that it converges in $C^{\infty}$ sense on $\pi^{-1}(U_i)$
 in the metric of $\widehat M_1$
 modulo gauge transformations.
\epr
\proof
The equations relevant to this part are $(h_1)_2$ and $(h_2)_2$.
We first consider $(A'_{\ep_{\nu}})_2$.

By the remark before the proposition,
 we can think of $(A'_{\ep_{\nu}})_2$
 as a section over $T_B^2$ of the bundle, whose fiber
 over $x \in T_B^2$ is $Rep_{T^2_{F, x}}(U(1))$,
 which is the quotient of the space $T^*_0 T_{F, x}^2$,
 where $0$ means the intersection of the fiber $T_{F, x}$ over $x$
 and the zero section, with the complex structure given by
 $*_{T^2_F}$. 
This is a holomorphic torus bundle.
In particular, the equation $(h_2)_2$ says that $(A'_{\ep_{\nu}})_2$
 is a solution of a holomorphic differential equation
 which does not depend on $\nu$.
The difference of two solutions can be locally seen as
 a holomorphic function.
Moreover, locally, by gauge transformation of the type
$$
g = \left( \begin{array}{cc}
         e^{2\pi i(nx - my)} & 0 \\
                0            & e^{-2\pi i(nx - my)}
          \end{array} \right),
$$
 we can assume
$$
inf_{U_i} |(A'_{\ep_{\nu}})_2| < 2 \pi.
$$

On the other hand, the cohomology classes of the curvatures 
 $[Tr(F_{\Xi'_{\ep_{\nu}}})]$ do not depend on $\nu$, since
 we have fixed the topological type.
Moreover, $2 c_0 \omega$, the self dual part of the curvature,
 does not depend on $\nu$, either.
So the cohomology classes of the ASD part (with respect to the rescaled
 metrics) of the curvature 
 $[Tr(F_{\Xi'_{\ep_{\nu}}}) - 2c_0 \omega]$
 do not depend on $\nu$.
Note that the self duality of $\omega$ is preserved in the rescaling.
In particular, the integral
 $\int |F_{\Xi'_{\ep_{\nu}}}|^2 d\mu_{\widehat M_{\ep_{\nu}}}$
 is uniformly bounded by some constant 
 not depending on $\nu$.

On the other hand, the curvature associated to $(A'_{\ep_{\nu}})_2$
 is of mixed type, that is, the terms of it are the wedge product of
 the fiber and the base direction differential forms.
The $L^2$-norm of these forms are not affected by rescaling.
So the $L^2$-norms of the curvatures related to $(A'_{\ep_{\nu}})_2$,
 namely, the $L^2$-norms of
 $\del_S (A'_{\ep_{\nu}})_2 dS$ and $\del_T (A'_{\ep_{\nu}})_2 dT$
 are also bounded in the metric of $\widehat M_1$.

Summarizing, there is locally a gauge in which $(A_{\ep_{\nu}})_2$
 can be seen as a holomorphic function with
 $inf_{U_i} |(A'_{\ep_{\nu}})_2| < 2 \pi$, and
 $\int_U |\del_S(A'_{\ep_{\nu}})_2|^2 d \mu_{T_B^2}$
 and
 $\int_U |\del_T(A'_{\ep_{\nu}})_2|^2 d \mu_{T_B^2}$
 are also bounded.
From this, it is easy to deduce that the $C^r$-norms of 
 $(A_{\ep_{\nu}})_2$ are bounded by constants independent of $\nu$.

On the other hand, the 
 $(\Phi'_{\ep_{\nu}} dS + \Psi'_{\ep_{\nu}} dT)_2$ part
 of the connections satisfy the equation $(h_1)_2$.
Moreover, we can locally choose the gauge 
 (without disturbing the gauge for $(A_{\ep_{\nu}})_2$)
 as in lemma 6.24.
In this gauge, this part becomes zero.
The proposition follows from these.\qed\\

We have now shown the existence of a subsequence of the traceless part of
 $\Xi_{\ep_{\nu}}$, which converges on 
 $\pi^{-1}(T_B^2 \setminus S)$ to a limit connection 
 on $E \to \widehat{M_1}$
 in a $C^{\infty}$ manner,
 and the existence of a locally converging subsequence of the fiberwise 
 constant part 
 of the
 $i \R$-part of the  connections.
We write (on one of the open disks $U_i$ of proposition 7.3)
 the limit of the latter as $\Xi' = A' + \Phi' ds + \Psi' dt$
 (we have returned to the original coordinates).
The fiber component $A^{\circ}$ of the traceless limit connection, which we write as
 $\Xi^{\circ} = A^{\circ} + \Phi^{\circ} ds + \Psi^{\circ} dt$.
 was shown to be flat, and so
 we can associate to it
 a family of points
 on the fibers of the dual torus fibration (that is, the mirror
 symplectic manifold) $M$ of the K\"ahler
 manifold $\widehat M_1$ with a K\"ahler $T^2$ structure. 
The part $A'$ plays the role of fiberwise parallel transportation
 of these points.
Our final task is to show that these points constitute
 a Lagrangian subvariety of $M$.
We write this subspace of $M$ as $L$.
Note that on the overwrapping regions of the covering 
 $\{ U_i \}$, the parallel transportation by $A'$
 can be assumed to be compatible, by diagonal argument. 

We denote by
$$
u : T^2_B \setminus S \to X,
$$
 the map defined by taking the 
 complex gauge equivalence classes
 of the fiber component of $\Xi$.
Here $X$ is the $S^2$ bundle over $T_B^2$ introduced in the last section,
 whose 
 fiber over $x \in T_B^2$ is $Rep_{T^2_F}(SO(3))$.

%
%
%
%
%
On the other hand, we have a limit map
$$
\phi : T^2_B \to X
$$
 which is the limit of the holomorphic maps $\phi_{\nu}$
 given by taking the complex gauge equivalence classes
 of the traceless parts of the
 connections $\Xi_{\ep_{\nu}}$ restricted to the fibers,
 modulo finite bubbles (proposition 6.3).
Recall that by fixing the natural
 complex structure on $T^2_B$,
 $\phi_{\nu}$ and so $\phi$ become holomorphic maps (remark 6.2).
The map $\phi$ coincides with $u$ on $T_B^2 \setminus S$.
In particular, the map $u : T^2_B \setminus S \to X$
 naturally extends to a holomorphic map from the whole $T^2_B$
 to $X$.
As mentioned, we interpret this as a double valued section over
 the dual torus fibration (that is, the mirror of $\widehat M_1$).
The fiberwise parallel transportations caused by
 the fiberwise constant part of  the $i \R$ part of the limit connection
 do not affect the diffeomorphism class of the section.

Now we recall the statement of the main theorem (theorem 2.6).
\bt
Let $M$ be the mirror symplectic manifold of $\widehat{M}_1$.
There is a double valued Lagrangian multisection, possibly non-reduced
 and possibly with
 ramifications
 for $M \ra T_B^2$ determined by the family
 $(E_{\ep_{\nu}}, \Xi_{\ep_{\nu}})$.
The ramifications occur at finitely many points.
%
Moreover, if 
 the first Chern class of the bundle $E$ is $0$, the
 multisection satisfies the special Lagrangian condition
 on the smooth part. 
\et
\br
The non-reduced case only occurs when the
 map $\phi$ above is a constant map to a singular point.
\er
\br
As we can see from the above argument, the map $\phi$ can be
 defined without referring to the convergence of the connections.
However, the fact that it really comes from the convergence of 
 the connections is important when one wants to compare
 gauge theory and holomorphic curve theory, as required for example
 in Atiyah-Floer conjecture (see also introduction).

Only in the case when the map $\phi$ is a constant map
 to a singular point of $Rep_{T^2}(SO(3))$, this
 convergence of the connections is not proved yet (see remark 6.24).
\er
The rest of this section is devoted to the proof of this theorem.
We assume the map $\phi$ is not a constant map 
 (in the case it is constant the following argument becomes simpler).
We already have a double valued section,
 which is locally diffeomorphic to a complex curve (so remark 2.7 follows).
The statement in the case of $c_1(E) = 0$ immediately follows from the
 hyperK\"ahler rotation.
We calculate the effect of the $i \R$ part of the connection 
 in the following.

Since at the (countable) loci $S$, the map defined by
 the limit connection extends continuously and
 in fact smoothly outside finite ramification points, 
 it suffices to prove the (special) Lagrangian
 property only over $T_B^2 \setminus S$. 
Moreover, note that the loci where
 the fiber part of the limit connection 
 have holonomy $U(2)$ are contained in $S$.
So the isotropy groups of the flat connections
 on the fibers over $T_B^2 \setminus S$ are all isomorphic.
The following is clear.
\bl
Let $U \subset T^2_B \setminus S$ be an open disc.
Then, there is a gauge on $E|_{\pi^{-1}(U)}$ 
 in which the fiber part of the limit 
 connection diagonalizes
 to constant matrices.
\qed
\el

Using a local frame of $E$ on $\pi^{-1}(U)$, $U$ is a small open subset of
 the base (recall $E$ is trivial when restricted to fibers and so
 we can take a frame globally in the fiber direction),
 denote $A = A^{\circ} + A'$ as $A = A_x dx + A_y dy$, 
 here $A_x$ and $A_y$ are the sums of 
 skew hermitian matrix-valued functions on
 $\pi^{-1}(U)$ and pure imaginary multiples of the identity matrix-valued
 functions on $\pi^{-1}(U)$ and both are fiberwise constant.

On a disc of the above lemma, 
%
 we can write $A$ as 
\begin{equation*}
\begin{array}{ll}
A = A_x dx + A_y dy &  \\
   \;\;\; = \left(  \left(\begin{array}{ll}
                  a(s, t) & 0\\
                     0      & -a(s, t) \end{array} \right) +
         \left(\begin{array}{ll}
                  c_1(s, t) & 0\\
                     0      & c_1(s, t) \end{array} \right) \right)dx & \\ 
 \hs{1in} + \left(  \left(\begin{array}{ll}
                  b(s, t) & 0\\
                     0      & -b(s, t) \end{array} \right) +
         \left(\begin{array}{ll}
                  c_2(s, t) & 0\\
                     0      & c_2(s, t) \end{array} \right) \right)dy & \\
  \;\;\; = \left(\begin{array}{ll}
                  a_1(s, t) & 0\\
                     0      & a_2(s, t) \end{array} \right)dx +
           \left(\begin{array}{ll}
                  b_1(s, t) & 0\\
                     0      & b_2(s, t) \end{array} \right) dy.
\end{array} 
\end{equation*}
We assume this disc is contained in some $U_i$ of the covering
 of $T^2_B$ in proposition 7.1.
The limit `connection'
 (in fact this is not a true connection on $\widehat M$.
But it suffices for our local argument below and the result
 does not depend on the gauge.)
 $\Xi = \Xi^{\circ} + \Xi' = A + \Phi ds + \Psi dt$ 
 satisfies the following equation:
$$
\ast_{\widehat M} (F_{mix} - c_0 \omega) = - (F_{mix} - c_0 \omega),
$$
where $F_{mix}$ is given by
 $F_{mix} = (d_A \Phi - \del_s A) ds + (d_A \Psi - \del_t A) dt$.
%
In this gauge, we can decompose the equation
 $\ast_{\widehat M} (F_{mix} - c_0 \omega) = - (F_{mix} - c_0 \omega)$,
 to the harmonic part and its fiberwise perpendicular part
 with respect to the fiberwise differentials
 $d_A$. 
The harmonic part
 (that is, the fiberwise constant part)
 is then of the following form (see below for the
 calculation):
\begin{equation*}
(\ast) \hs{1in}  \begin{array}{cc}
\del_s a_i + \del_t b_i - c_0 = 0 \\
(g^{ss}\del_s + 2g^{st}\del_t) b_i - g^{tt}\del_t a_i - c_0 g^{st} = 0.
\end{array}
\end{equation*}

%

$a_i(s, t)$ and $b_i(s, t)$ are points of the dual torus of
 the fiber of $\widehat M_1$, that is, the fiber of the mirror $M$.
These points consist a multi-section on $T_B^2$ of the 
 mirror torus fibration $M$. 
This is the locus $L$ we defined at the beginning of this section.

What is left to prove is that the smooth parts
 of our multisections are  
 Lagrangians, and when $c_1 = 0$, they are indeed special.
This is a partial inverse of the Fourier transformation from
 A-cycles to B-cycles in \cite{AP} 
 (partial means we do not get bundles on Lagrangians.
 See the next section).
This is also discussed in the paper \cite{LYZ}
 in the $U(1)$-bundle case.
In $\cite{LYZ}$, they treat so called $MMMS \; equations$
$$ \begin{array}{cc} F_A^{0,2} = 0, \\
                     Im e^{i \theta}(\omega + F_A)^n = 0,
\end{array} $$
 here $n$ is the complex dimension of the manifold.
These are presumed to be a deformation of the HYM equation.
Although the following argument should be known to experts, we
 give an explicit calculation for the sake of completeness and
 also for clarifying the relation between the specialness
 and the Chern class of the mirror bundle.

Now we begin to prove that the locus $L$ is a Lagrangian
 subvariety of $M$.
Recall that the mirror symplectic manifold $(M, \Check{\omega})$
 of $(\widehat M_1, \omega_1, J_1)$ has Darboux coordinates
 $(\Check s, \Check t, x^*, y^*)$ and 
 $\Check{\omega} = d\Check s \wedge dx^* + d\Check t \wedge dy^*$
 in this coordinates.

Taking the connection $A$ in the above form,
 the tangent space of $L$ at $(s, t, a_i(s, t), b_i(s, t))$
 is spanned by vectors of the form
\begin{eqnarray}
 \ell_1 = \del_{\Check s} + \del_{\Check s} a_i \del_{x^*}
    + \del_{\Check s} b_i \del_{y^*},\\
 \ell_2 = \del_{\Check t} + \del_{\Check t} a_i \del_{x^*}
    + \del_{\Check t} b_i \del_{y^*}.
\end{eqnarray}

Substituting these vectors to the symplectic form $\Check{\omega}$,
 the condition that the submanifold $L$ is Lagrangian
 is given by the following equation
\begin{equation}
g^{st} \del_s A_x + g^{tt} \del_t A_x - g^{ss} \del_s A_y - g^{st} \del_t A_y
 = 0.
\end{equation} 

Now, looking at
$$F_{mix} - c_0 \omega =  
       (d_A \Phi - \del_s A + c_0 dx) ds +
                     (d_A \Psi - \del_t A + c_0 dy) dt,$$
 we calculate
 $\ast_{\widehat M}\{(- \del_s A + c_0 dx) ds + (- \del_t A + c_0 dy) dt\}$:
\begin{eqnarray}
\ast_{\widehat M} (- \del_s A ds - \del_t A dt) \hs{3in} \\\nonumber
 = (g^{ss} g_{ss} \del_s A_x + g^{ss} g_{st} \del_s A_y
        + g^{st} g_{ss} \del_t A_x + g^{st} g_{st} \del_t A_y) dt \wedge dy
\\\nonumber
 - (g^{st} g_{ss} \del_s A_x + g^{st} g_{st} \del_s A_y
        + g^{tt} g_{ss} \del_t A_x + g^{tt} g_{st} \del_t A_y) ds \wedge dy
\\\nonumber
 - (g^{ss} g_{st} \del_s A_x + g^{ss} g_{tt} \del_s A_y
        + g^{st} g_{st} \del_t A_x + g^{st} g_{tt} \del_t A_y) dt \wedge dx
\\\nonumber
 + (g^{st} g_{st} \del_s A_x + g^{st} g_{tt} \del_s A_y
        + g^{tt} g_{st} \del_t A_x + g^{tt} g_{tt} \del_t A_y) ds \wedge dx.
\end{eqnarray}

Since $\omega$ is a self dual form, we have
\begin{equation}
\ast_{\widehat M} c_0 \omega = c_0 \omega
\end{equation}

Using these and fiberwise harmonic part of the ASD equation,
 we have the following identities:
\begin{eqnarray}
(i) \; g^{ss} g_{ss} \del_s A_x + g^{ss} g_{st} \del_s A_y
      + g^{st} g_{ss} \del_t A_x + (g^{st} g_{st} + 1) \del_t A_y - c_0 = 0\\
(ii)\; g^{st} g_{ss} \del_s A_x + (g^{st} g_{st} - 1) \del_s A_y
        + g^{tt} g_{ss} \del_t A_x + g^{tt} g_{st} \del_t A_y = 0 \hs{.7cm}\\
(iii)\; g^{ss} g_{st} \del_s A_x + g^{ss} g_{tt} \del_s A_y
        + (g^{st} g_{st} - 1) \del_t A_x + g^{st} g_{tt} \del_t A_y = 0
 \hs{.7cm}\\
(iv)\; (g^{st} g_{st} + 1) \del_s A_x + g^{st} g_{tt} \del_s A_y
        + g^{tt} g_{st} \del_t A_x + g^{tt} g_{tt} \del_t A_y - c_0 = 0.
\end{eqnarray}

Now $(i) + (iv)$ gives
$$
(A) : \del_s A_x + \del_t A_y - c_0 = 0.
$$

Further,
$$
(B): (i) \times g^{st} - (ii) \times g^{ss}
                = g^{ss}\del_s A_y - g^{tt} \del_t A_x 
                                + 2g^{st} \del_t A_y - c_0 g^{st} = 0.
$$
The identity $(A) \times g^{st} - (B) = 0$ gives precisely the 
 Lagrangian condition for $L$ given above.

\begin{note}
In the above calculation, we do not need the Calabi-Yau
 condition $det g = 1$.
So if we can prove the convergence of the connection in non Calabi-Yau
 situation, we can construct the corresponding Lagrangian on the mirror.
\end{note}
Next, we prove the specialness of $L$ 
 when the 
 $c_1$ of the bundles vanish
 (although this is immediate from hyperK\"ahler rotation,
 we prove it by calculation to see the role of the constant $c_0$).
Since we have already proved that $L$ is Lagrangian,
 to prove it is special it suffices to show the imaginary part of the
 holomorphic volume form restricts to zero on $L$.

The imaginary part of $dz_1 \wedge dz_2$
 is given by
$$
Im(dz_1 \wedge dz_2) = ds\wedge dy^* - dt\wedge dx^*.
$$
Substituting the tangent vectors $\ell_1, \ell_2$ of $L$ given above, 
 we have
$$
\begin{array}{ll}
Im(dz_1 \wedge dz_2)(\ell_1, \ell_2)
 & = g^{ss}\del_{\Check t} b_i - g^{st}\del_{\Check s}b_i
   + g^{tt}\del_{\Check s}a_i - g^{st}\del_{\Check t}a_i \\
 & = g^{ss}g^{st}\del_s b_i + g^{ss}g^{tt}\del_t b_i 
     - g^{st}g^{ss}\del_s b_i - (g^{st})^2\del_t b_i \\

   & \;\;\; + g^{tt}g^{ss}\del_s a_i + g^{tt}g^{st}\del_t a_i 
     - (g^{st})^2\del_s a_i - g^{st}g^{tt}\del_t a_i \\
 & = \del_t b_i + \del_s a_i
\end{array}
$$
 using $det g = 1$.

Since we have assumed $c_1 = 0$, the constant $c_0$
 becomes zero
 in the identities $(i)$ and $(iv)$ above. 
So we have
$$
(i) + (iv) = 2(\del_t b_i + \del_s a_i) = 0.
$$
Now our theorem is proved. \qed

\section{Miscellaneous}

Here we make some remarks.\\

\nnn
1. About flat bundles \\

We mention about the flat connection on $L$.
From the usual mirror symmetry point of view, it is
 desirable to attach to our Lagrangian subvariety 
 a flat vector bundle in a natural manner.
In fact, there is an obvious candidate for it as mentioned below.

Here again we talk about ASD connections on $SO(3)$ bundle
 and hyperK\"ahler rotate the base as in section 6.
In a gauge which diagonalize $A$ with constant components on each fiber,
 the $d_A$ exact part and the $d_A$ coexact parts of the
 $(H_1)$ part (in the notation of the beginning
 of section 6) of the ASD equations are equivalent to
$$
d_A \Phi = 0, d_A \Psi = 0.
$$
 as we mentioned several times.

That is, $\Phi$ and $\Psi$ are sections which are constant and diagonal
 on each fiber (we assume $A \neq 0$).
Moreover, note that the estimate of corollary 6.23 says that the curvature
 of the fiber part $\pll F_{A_{\ep_{\nu}}} \pll_{L^{\infty}}$
 converges exponentially fast to zero on $T_B^2 \setminus S$.
So the equation 
$$
\del_t \Phi_{\ep_{\nu}} - \del_s \Psi_{\ep_{\nu}}
 - [\Phi_{\ep_{\nu}}, \Psi_{\ep_{\nu}}] +
 \ep_{\nu}^{-2} F_{A_{\ep_{\nu}}} 
 = 0
$$
 yields the equation for the limit connection
$$
\del_t \Phi - \del_s \Psi - [\Phi, \Psi] = 0.
$$
That is, the 1-form $\Phi ds + \Psi dt$
 defines a flat connection on the the restriction of 
 $E$ to the smooth part of $L$.

This bundle is rank two and so we want to split
 it to line bundles.
However, to do this in a well-defined way, 
 the ramification of the Lagrangian and that of the
 bundle (roughly speaking, the configuration of points
 at which $\Phi$ and $\Psi$ becomes zero)
 must be compatible, and we cannot prove it at present.\\


\noindent
2. Gauge theory on Ooguri-Vafa spaces.\\

In the paper \cite{OV}, physicists
 H. Ooguri and C. Vafa constructed a
 Ricci-flat K\"ahler metric on a neibourhood of
 an $I_1$ fiber of an elliptically fibered K3 surface.
Gross and Wilson \cite{GW} used this to
 construct a Ricci-flat K\"ahler metric on
 a K3 surface which is near the large structure limit.

Then it is natural to consider gauge theory on K3
 with these metric.
It is the extension of this paper, where we treated
 $T^4$.

Obviously, gauge theory on Ooguri-Vafa spaces
 will be the key ingredient for it.
The gauge theory on the whole K3 should be 
 studied by some appropriate gluing of this and
 the gauge theory on the nonsingular part, which
 is studied in this paper.
In fact, the metrics on Ooguri-Vafa spaces
 are known to converge exponentially fast to semi-flat
 metrics at infinity.
So it would be possible to apply our methods to these spaces.
In particular, it is quite plausible that
 we can prove the energy quantization of ASD connections
 on Ooguri-Vafa spaces, 
 whether the bundle is trivial when restricted to the fiber or not.\\

\noindent
3. Donaldson-Floer theory with degenerate boundary condition.\\

There is an expectation that Donaldson invariants of a closed manifold
 $M = M_1 \cup_Y M_2$ can be computed from some pairing of the Floer-homology
 valued invariants of the components $M_1, M_2$.
This has been realized when $Y$ is a homology 3-sphere \cite{D}.
There has been attempts to extend it to more general situations \cite{MMR},
 \cite{T}.
It will be interesting to try to apply our methods to these
 cylinder end cases, while we treated cone end cases in 
 this paper.
I hope that our method can be applied to some of these situations.
In particular, it may be suited for the analysis of bubbles,
 as is done in this paper.

\nnn
Takeo Nishinou, Department of Mathematics, Faculty of Science,
 Kyoto University, Kitashirakawa, Kyoto 606-8502, Japan\\
Current address: Mathematisches Institut, Albert-Ludwigs-Universit\"at,
 Eckerstrasse 1. 79104, Freiburg, Germany

  nishinou@kusm.kyoto-u.ac.jp
\end{document}